\documentclass[a4paper]{amsart} 
\usepackage[ascii]{inputenc} 

\usepackage{amssymb}
\usepackage{mathrsfs}
\usepackage[hidelinks]{hyperref}

\usepackage{xcolor}

\newcommand{\cyan}[1]{\textcolor{cyan}{#1}}

\usepackage[shortlabels]{enumitem}
\setlist[enumerate,1]{label={(\Alph*)}}
\setlist[enumerate,2]{label={(\alph*)}}
\setlist[enumerate,3]{label={$\bullet_{\arabic*}$}}

\newenvironment{PROOF}[2][\proofname.]
   {\begin{proof}[#1]\renewcommand{\qedsymbol}{\textsquare$_{\rm #2}$}}
   {\end{proof}}

\newtheorem{theorem}{Theorem}[section] 
\newtheorem{claim}[theorem]{Claim}
\newtheorem{lemma}[theorem]{Lemma} 
\newtheorem{thesimplelemma}[theorem]{The Simple B.B. Lemma} 
 
\newtheorem{conclusion}[theorem]{Conclusion}

\theoremstyle{definition}
\newtheorem{definition}[theorem]{Definition}
\newtheorem{example}[theorem]{Example}
\newtheorem{problem}[theorem]{Problem}
\newtheorem{fact}[theorem]{Fact}
\newtheorem{subfact}[theorem]{Subfact}
\newtheorem{exercise}[theorem]{Exercise}
\newtheorem{discussion}[theorem]{Discussion}

\theoremstyle{remark}
\newtheorem{remark}[theorem]{Remark}

\newtheorem{notation}[theorem]{Notation}

\newcommand{\cd}{\mathrm{cd}}
\newcommand{\rel}{\mathrm{rel}}
\newcommand{\St}{\mathsf{St}}

\newcommand{\closure}{\mathrm{closure}}

\newcommand{\rank}{\mathrm{rank}}

\newcommand{\tthh}{\mathrm{th}}

\newcommand{\acc}{\mathrm{acc}}

\newcommand{\cf}{\mathrm{cf}}

\newcommand{\Hom}{\operatorname{Hom}}
\newcommand{\id}{\mathrm{id}}

\newcommand{\Ord}{\mathrm{Ord}}
\newcommand{\otp}{\mathrm{otp}}

\newcommand{\Rang}{\mathrm{Rang}}
\newcommand{\rang}{\mathrm{rang}}
\newcommand{\seq}{\mathrm{seq}}

\newcommand{\bd}{\mathrm{bd}}

\newcommand{\Dom}{\mathrm{Dom}}

\newcommand{\EM}{\mathrm{EM}}

\newcommand{\tcf}{\mathrm{tcf}}

\newcommand{\tr}{\mathrm{tr}}

\newcommand{\pp}{\mathrm{pp}}

\newcommand{\bfB}{\mathbf{B}}

\newcommand{\bfC}{\mathbf{C}}
\newcommand{\bfc}{\mathbf{c}}

\newcommand{\bfF}{\mathbf{F}}
\newcommand{\bff}{\mathbf{f}}

\newcommand{\bfL}{\mathbf{L}}

\newcommand{\bft}{\mathbf{t}}

\newcommand{\bfV}{\mathbf{V}}
\newcommand{\bfW}{\mathbf{W}}

\newcommand{\bfx}{\mathbf{x}}

\newcommand{\bbI}{\mathbb{I}}

\newcommand{\bbP}{\mathbb{P}}
\newcommand{\bbQ}{\mathbb{Q}}

\newcommand{\bbZ}{\mathbb{Z}}

\newcommand{\mn}{\medskip\noindent}
\newcommand{\sn}{\smallskip\noindent}
\newcommand{\bn}{\bigskip\noindent}

\newcommand{\cA}{\mathscr{A}}

\newcommand{\cE}{\mathscr{E}}
\newcommand{\cF}{\mathscr{F}}

\newcommand{\cL}{\mathscr{L}}
\newcommand{\cM}{\mathscr{M}}

\newcommand{\cU}{\mathscr{U}}

\newcommand{\clH}{\mathcal{H}}

\newcommand{\clP}{\mathcal{P}}

\newcommand{\clS}{\mathcal{S}}
\newcommand{\clT}{\mathcal{T}}

\newcommand{\gA}{\mathfrak{A}}

\newcommand{\gC}{\mathfrak{C}}

\newcommand{\eps}{\varepsilon}
\newcommand{\lh}{{\ell g}}
\newcommand{\rest}{\restriction}

\newcommand{\caret}{{\char 94}}
\newcommand{\LL}{\langle}
\newcommand{\RR}{\rangle}

\newcommand{\olsi}[1]{\,\overline{\!{#1}}} 

\newcommand*{\defeq}{\mathrel{\vcenter{\baselineskip0.5ex \lineskiplimit0pt\hbox{\scriptsize.}\hbox{\scriptsize.}}}=}

\newcommand{\mgrimes}[1]{\cyan{\textbf{[#1]}}}

\newcount\skewfactor
\def\mathunderaccent#1#2 {\let\theaccent#1\skewfactor#2
\mathpalette\putaccentunder}
\def\putaccentunder#1#2{\oalign{$#1#2$\crcr\hidewidth
\vbox to.2ex{\hbox{$#1\skew\skewfactor\theaccent{}$}\vss}\hidewidth}}
\def\name{\mathunderaccent\tilde-3 }
\def\Name{\mathunderaccent\widetilde-3 }

\newbox\noforkbox \newdimen\forklinewidth
\setbox0\hbox{$\textstyle\bigcup$}
\setbox1\hbox to \wd0{\hfil\vrule width \forklinewidth depth \dp0
                        height \ht0 \hfil}
\setbox\noforkbox\hbox{\box1\box0\relax}
\def\unionstick{\mathop{\copy\noforkbox}\limits}
\def\nonfork#1#2_#3{#1\unionstick_{\textstyle #3}#2}
\def\nonforkin#1#2_#3^#4{#1\unionstick_{\textstyle #3}^{\textstyle
    #4}#2}
\setbox0\hbox{$\textstyle\bigcup$}
\setbox1\hbox to \wd0{\hfil{\sl /\/}\hfil}
\setbox2\hbox to \wd0{\hfil\vrule height \ht0 depth \dp0 width
                                \forklinewidth\hfil}
\newbox\doesforkbox
\setbox\doesforkbox\hbox{\box1\box0\relax}
\def\nunionstick{\mathop{\copy\doesforkbox}\limits}

\def\fork#1#2_#3{#1\nunionstick_{\textstyle #3}#2}
\def\forkin#1#2_#3^#4{#1\nunionstick_{\textstyle #3}^{\textstyle
    #4}#2}

\newcommand{\stickT}{%
\setbox255=\hbox{\raise1ex\hbox{$\hspace{0.2pt}\,\bullet\,$}}
\mathord{\rlap{\hbox to\wd255{\hss\hbox{$|$}\hss}}
\box255}
}
\newcommand{\stickS}{%
\setbox255=\hbox{\raise0.6ex\hbox{$\scriptstyle\bullet$}}
\mathord{\rlap{\hbox to\wd255{\hss\hbox{$\scriptstyle|$}\hss}}
\box255}
}

\begin{document}
\makeatletter\def\shfiuwefootnote{\gdef\@thefnmark{}\@footnotetext}\makeatother\shfiuwefootnote{Version 2023-04-26. See \url{https://shelah.logic.at/papers/309/} for possible updates.}

\title {Black Boxes \\
Sh309}
\author {Saharon Shelah}
\address{Einstein Institute of Mathematics\\
Edmond J. Safra Campus, Givat Ram\\
The Hebrew University of Jerusalem\\
Jerusalem, 91904, Israel\\
 and \\
 Department of Mathematics\\
 Hill Center - Busch Campus \\ 
 Rutgers, The State University of New Jersey \\
 110 Frelinghuysen Road \\
 Piscataway, NJ 08854-8019 USA}
\email{shelah@math.huji.ac.il}
\urladdr{http://shelah.logic.at}
\thanks{For versions up to 2019, the author thanks Alice Leonhardt for the beautiful typing. In the latest version, the author thanks an individual who wishes to remain anonymous for generously funding typing services, and thanks Matt Grimes for the careful and beautiful typing. For their partial support of this research, the author would like to thank: an NSF-BSF 2021 grant with M. Malliaris, NSF 2051825, BSF 3013005232 (2021/10-2026/09); and for various grants from the BSF (United States Israel Binational Foundation), the Israel Academy of Sciences and the NSF via Rutgers University.
This paper is number 309 in the author's publication list.\\
This is a revised version of \cite[Ch.III,\S4,\S5]{Sh:300}; it has existed 
(and been occasionally revised) for many years. It was mostly ready in the early nineties, and was made public to some extent. This was written as Chapter IV of 
the book \cite{Sh:e}, which hopefully will materialize some day, but in the 
meantime it is \cite{Sh:309}. The intention was to have \cite{Sh:E58} 
(revising \cite{Sh:229}) for Ch.I, \cite{Sh:421} for Ch.II, \cite{Sh:E59} for 
Ch.III, \cite{Sh:309} for Ch.IV, \cite{Sh:363} for Ch.V, \cite{Sh:331} for Ch.VI,
\cite{Sh:511} for Ch.VII, \cite{Sh:E60} (a revision of \cite{Sh:128}) for Ch.VIII,
\cite{Sh:E62} for the appendix, and probably \cite{Sh:757}, \cite{Sh:384}, \cite{Sh:482}, and \cite{Sh:800}. References like \cite[3.26=L6.12]{Sh:309}
means that 6.12 is the label of Lemma \ref{6.12} in \cite{Sh:309}.\\
The reader should note that the version in my website is usually more
up-to-date than the one in the mathematical archive.}

\makeatletter
\@namedef{subjclassname@2020}{\textup{2020} Mathematics Subject Classification}
\makeatother



\date {April 24, 2023}

\begin{abstract}
We shall deal comprehensively with Black Boxes, the intention being
that \underline{provably in ZFC}  we have a sequence of guesses of 
extra structure on small subsets, {where} the guesses are 
pairwise ``almost disjoint;'' by this we mean they have quite 
little interaction, and are far {apart} but together are ``dense."
We first deal with the simplest case, where the existence comes
from winning a game by just writing down the opponent's
moves. We show how it helps when instead of orders
we have trees with boundedly many levels, having freedom
in the last. After this we quite systematically look at existence 
of black boxes, and make connection to non-saturation of 
natural ideals and diamonds on them.
\end{abstract}

\maketitle
\numberwithin{equation}{section}
\setcounter{section}{-1}
\newpage

\section{Introduction}

The non-structure theorems we have discussed in \cite{Sh:E59} 
usually rest on some freedom on finite sequences and on a kind of order.  
When our freedom is related to infinite sequences, and to trees, 
our work is sometimes harder. In particular, we may consider, 
for $\lambda \ge \chi$, $\chi$ regular, and 
$\varphi = \varphi(\bar x_0,\ldots,\bar x_\alpha,\ldots)_{\alpha<\chi}$ 
in a vocabulary $\tau$:
\sn
\begin{enumerate}
    \item[$(*)$]   For any $I \subseteq {}^{\chi\ge}\lambda$ we have a $\tau$-model $M_I$ and sequences $\bar a_\eta$ (for $\eta\in{}^{\chi>}\!\lambda$), where
\[
[\eta\lhd\nu\ \Rightarrow\ {\bar a}_\eta \ne {\bar a}_\nu],\qquad
\lh({\bar a}_\eta) = \lh({\bar x}_{\lh(\eta)}),
\]
such that for $\eta\in {}^\chi\lambda$ we have:
\[
M_I\models\varphi (\ldots,{\bar a}_{\eta\rest\alpha},
\ldots)_{\alpha<\chi } \text{ if and only if } \eta\in I.
\]
\end{enumerate}
\mn
(Usually, $M_I$ is to some extent ``simply defined'' from $I$). Of course,
if we do not ask more from $M_I$, we can get nowhere: we certainly restrict
its cardinality  and/or usually demand it is $\varphi$-representable\footnote{see
Definition \cite[2.7=Lf4]{Sh:E59} clauses (c),(d).} in (a variant of) 
$\cM_{\mu,\kappa }(I)$ (for suitable $\mu,\kappa$). 
Certainly for $T$ un-superstable we have such a formula $\varphi$:
\[
\varphi(\ldots,\bar{a}_{\eta\rest n},\ldots) = \big(\exists{\bar x} \big)
\textstyle\bigwedge\limits_n\varphi_n({\bar x},{\bar a}_{\eta\rest n}).
\]

\mn
There are many natural examples.

Formulated in terms of the existence of $I$ for which our favorite
``anti-isomorphism'' player has a winning strategy, we proved this in 1969/70
(in proofs of lower bounds of 
$\dot{\bbI}(\lambda,T_1,T)$, $T$ un-superstable), but
it was shortly superseded. However, eventually the method was used in one of
the cases in \cite[Ch.VIII,\S2]{Sh:a} --- for strong limit singular
\cite[Ch.VIII,2.6]{Sh:a}, which comes from \cite{Sh:31}. It was developed in \cite{Sh:172}, \cite{Sh:227} for
constructing Abelian groups with prescribed endomorphism groups. See further
a representation  of one of the results here in Eklof-Mekler \cite{EM},
\cite{EM02} a version which was developed for a proof of the existence 
of an Abelian (torsion-free $\aleph_1$-free) group $G$ with
\[
G^{***} = G^* \oplus A \qquad \big(G^* \defeq \Hom(G,\bbZ)\big)
\]

\mn
in a work by Mekler and Shelah. A preliminary version of this paper
appeared in \cite[Ch.III,\S4,\S5]{Sh:300}, but \S3 here was just 
almost ready and \S4 (on partitions of stationary sets and $\diamondsuit_D$)
was written up as a letter to Foreman in the late nineties.

The saturation of ideals was continued much later in Gitik-Shelah \cite{Sh:898} and more recently in \cite{Sh:1028} and Asgarzadeh-Golshani-Shelah \cite{Sh:1232}. 

\newpage

\section {The Easy Black Box and an Easy Application}

In this section we do not try to get the strongest results, but just provide
some examples (e.g. we do not present the results when $\lambda=
\lambda^\chi$ is replaced by $\lambda=\lambda^{<\chi})$. By the proof of
\cite[Ch.VIII,2.5]{Sh:a} (see later for a complete proof):

\begin{theorem}\label{4.2}
Suppose that
\mn
\begin{enumerate}
    \item[$(*)$]
    \begin{enumerate}
        \item[$(a)$]  $\lambda=\lambda^\chi$
\sn
        \item[$(b)$]  $\tau$ is a vocabulary and $\varphi = 
        \varphi(\bar x_0,\bar x_1,\ldots,\bar x_{\alpha}\ldots)_{\alpha<\chi}$ 
        is a formula in $\cL(\tau)$ for some logic $\cL$.
\sn
        \item[$(c)_{\tau,\varphi}$]   For any $I$ such that  
        ${}^{\chi>}\!\lambda \subseteq I \subseteq {}^{\chi\ge}\lambda$, 
        we have a $\tau$-model $M_I$ and sequences 
        $\bar a_\eta$ (for $\eta\in {}^{\chi>}\!\lambda$), where
        \[
        [\eta \lhd \nu \Rightarrow \bar a_\eta \ne \bar a_\nu],\qquad
        \lh(\bar a_\eta) = \lh(\bar x_{\lh(\eta)}),
        \]
        such that for $\eta\in {}^\chi\lambda$ we have:
        \[
        M_I \models \varphi (\ldots,{\bar a}_{\eta\rest\alpha},\ldots)_{\alpha<\chi} \text{ if and only if } \eta\in I.
        \]

        \item[$(c)$]  $\|M_I\| = \lambda$ for every $I$ satisfying ${}^{\chi>}\!\lambda \subseteq I\subseteq {}^{\chi\le}\lambda$, and $\lh({\bar a}_\eta) \le \chi$ or just $\lambda^{\lh({\bar a}_\eta)}= \lambda$.
    \end{enumerate}
\end{enumerate}
\mn
\underline{Then} (using ${}^{\chi>}\!\lambda\subseteq I\subseteq
{}^{\chi\ge}\lambda)$:

\sn
1) There is no model $M$ of cardinality $\lambda$ into which every $M_I$ can be
$(\pm\varphi)$-embedded (i.e., by a function preserving $\varphi$ and
$\neg\varphi)$.

\sn
2) For any $M_i$ (for $i<\lambda$), $\|M_i\| \le \lambda$, for some $I$
satisfying ${}^{\chi>}\!\lambda\subseteq I\subseteq {}^{\chi\ge}\lambda$, the
model $M_I$ cannot be $(\pm\varphi)$-embedded into any $M_i$.
\end{theorem}

\begin{example}
\label{4.3}
Consider the class of Boolean algebras and the formula
\[
\varphi(\ldots,x_n,\ldots) \defeq \big(\bigcup\limits_{n} x_n\big) = 1
\]

\sn
(i.e., there is no $x \ne 0$ such that $x \cap x_n=0$ for each $n$).

\noindent
For ${}^{\omega>}\!\lambda\subseteq I\subseteq {}^{\omega\ge}\lambda$, let
$M_I$ be the Boolean algebra generated freely by $x_\eta$ (for $\eta\in I$)
except the relations:  for $\eta \in I$, if $n < \lh(\eta) =
\omega$ then $x_\eta\cap x_{\eta\rest n}=0$.

So\footnote{With more work, we can demand that $M_I$ satisfies the c.c.c.} $\|M_I\| = |I| \in [\lambda,\lambda^{\aleph_0}]$ and in $M_I$ for 
$\eta \in {}^\omega\!\lambda$ we have: $M_I \models 
``\big(\bigcup\limits_n x_{\eta\rest n} \big) = 1$" if and only if 
$\eta \notin I$ (work a little in Boolean algebras). 
\end{example}

\noindent
So

\begin{conclusion}\label{4.4}
If $\lambda = \lambda^{\aleph_0}$, \underline{then} there is no 
Boolean algebra $\bfB$ of cardinality $\lambda$ universal 
under $\sigma$-embeddings (i.e., ones preserving countable unions). 
\end{conclusion}

\begin{remark} 
\label{4.4r}
This is from \cite[Ch.VIII,Ex.2.5,pg.464]{Sh:a}.
\end{remark}

\bn
\textbf{Proof of the Theorem \ref{4.2}.}
First we recall the simple black box (and a variant) in \ref{4.5A}, \ref{4.5D} below:

\begin{thesimplelemma}\label{4.5A}
There are functions $f_\eta$ (for $\eta\in {}^\chi\lambda$) such that:
\mn
\begin{enumerate}
    \item[$(i)$]   $\Dom(f_\eta)=\{\eta\rest\alpha:\alpha<\chi\}$,
\sn
    \item[$(ii)$]  $\Rang(f_\eta)\subseteq\lambda$,
\sn
    \item[$(iii)$]  If $f:{}^{\chi>}\!\lambda\to\lambda$, then for some $\eta\in {}^\chi\lambda$ we have $f_\eta\subseteq f$.
\end{enumerate}
\end{thesimplelemma}

\begin{PROOF}{\ref{4.5A}}
For $\eta\in {}^\chi\lambda$, let $f_\eta$ be the function (with
domain $\{\eta \rest \alpha:\alpha<\chi\})$ such that
\[
f_\eta(\eta \rest \alpha)=\eta(\alpha).
\]

\sn
So $\LL f_\eta : \eta\in {}^\chi\lambda\RR$ is well
defined. Properties (i) and (ii) are straightforward, so let us prove
(iii).  Let $f:{}^{\chi >}\!\lambda\to\lambda$. We define
$\eta_\alpha = \LL\beta_i : i < \alpha\RR$ by induction on $\alpha$.

For $\alpha=0$ or $\alpha$ limit --- no problem.

For $\alpha+1$: let $\beta_\alpha$ be the ordinal such that $\beta_\alpha
=f(\eta_\alpha)$.

\sn
So $\eta\defeq\LL\beta_i:i<\chi\RR $ is as required.
\end{PROOF}

\begin{fact}\label{4.5D}
In \ref{4.5A}:
\mn
\begin{enumerate}
    \item We can replace the range of $f,f_{\eta}$ by any fixed set of power $\lambda$.
\sn
    \item We can replace the domains of $f,f_{\eta}$ by 
    $\{{\bar a}_\eta : \eta\in {}^{\chi>}\!\lambda\}$, 
    $\{{\bar a}_{\eta \rest \alpha} : \alpha < \chi\}$, respectively, as long as
    \[
        \alpha < \beta < \chi \wedge \eta\in {}^\chi\lambda\ \Rightarrow 
        \ {\bar a}_{\eta \rest \alpha} \ne {\bar a}_{\eta\rest \beta}.
    \]
\end{enumerate}
\end{fact}

\begin{remark}\label{4.5E}
We can present it as a game. (As in the book \cite[Ch.VIII,2.5]{Sh:a}).
\end{remark}

\mn
\textbf{Continuation of the Proof of Theorem \ref{4.2}.}

It suffices to prove \ref{4.2}(2). Without loss of generality $\LL
|M_i|:i<\lambda\RR$ are pairwise disjoint. Now we use \ref{4.5D}; for
the domain we use $\LL{\bar a}_\eta:\eta\in {}^{\chi>}\!\lambda\RR$
from the assumption of \ref{4.2}, and for the range: $\bigcup\limits_{i<
\lambda} {}^{\chi\ge}|M_{i}|$ (it has cardinality $\le \lambda$ as
$\|M_i\|\leq\lambda=\lambda^\chi$). 

We define
\begin{align*}
I = ({}^{\chi>}\!\lambda) \cup \big\{\eta\in {}^\chi\lambda:
&\text{ for some } i< \lambda,\ \Rang(f_\eta) \text{ is a set of sequences}\\
  &\text{ from } |M_i| \text{ and } M_i \models \neg 
\varphi(\ldots,f_\eta({\bar a}_{\eta \rest \alpha}),
\ldots)_{\alpha<\chi} \big\}.
\end{align*}

\mn
Look at $M_I$. It suffices to show:
\mn
\begin{enumerate}
    \item[$\otimes$]   There is no $(\pm\varphi)$-embedding of 
    $M_I$ into $M_i$ for $i<\lambda$.
\end{enumerate}
\mn
Why does $\otimes$ hold?

If $f:M_I\to M_i$ is a $(\pm\varphi)$-embedding, then by Fact
\ref{4.5D}, for some $\eta\in {}^\chi\lambda$ we have
\[
f \rest \{{\bar a}_{\eta \rest \alpha} : \alpha < \chi\} = f_\eta.
\]

\mn
By the choice of $f$,
\[
M_{I}\models\varphi\left[\ldots,{\bar a}_{\eta\rest\alpha},\ldots\right]_{
\alpha<\chi} \iff M_i\models\varphi\left[\ldots,f({\bar a}_{\eta \rest
\alpha}),\ldots\right]_{\alpha<\chi},
\]

\mn
but by the choice of $I$ and $M_I$ we have
\[
M_I\models\varphi\left[\ldots,{\bar a}_{\eta\rest\alpha},\ldots\right]_{
\alpha<\chi}\iff M_i\models\neg\varphi\left[\ldots, f_\eta({\bar a}_{\eta
\rest \alpha}),\ldots\right]_{\alpha<\chi}.
\]

\mn
This is a contradiction, as by the choice of $\eta$,
\[
\bigwedge\limits_{\alpha<\chi} f({\bar a}_{\eta \rest \alpha})
= f_\eta(\bar{a}_{\eta \rest \alpha}). 
\]
QED. \hfill\textsquare$_{\ref{4.2}}$

\bn
\begin{discussion}\label{4.5G}
We may be interested whether, in \ref{4.2}, when 
$\lambda^+ < 2^\lambda$ we may
\mn
\begin{enumerate}
    \item   in \ref{4.2}(1), allow $\|M\|=\lambda^+$, and/or
\sn
    \item  get $\ge \lambda^{++}$ non-isomorphic models of the form $M_I$, assuming $2^\lambda>\lambda^+$.
\end{enumerate}
\mn
The following lemma shows that we cannot prove those better statements in
ZFC, though (see \ref{4.7}) in some universes of set theory we can. So this
requires (elementary) knowledge of forcing, but is not used later. 
It is here just to justify the limitations of what we can prove, and
the reader can skip it.
\end{discussion}

\begin{lemma}\label{4.6}
Suppose that in the universe $\bfV$ we have
$\kappa<\lambda=\cf(\lambda)=\lambda^{<\lambda}$, 
$(\forall\lambda_1<\lambda) \big[(\lambda_1)^\kappa<\lambda \big]$, 
and $\lambda<\mu=\mu^\lambda$.

\underline{Then}, for some notion forcing $\bbP$:
\begin{enumerate}
    \item[$(a)$]  $\bbP$ is $\lambda$-complete and satisfies the 
    $\lambda^+$-c.c., and $|\bbP|=\mu$, $\Vdash_\bbP ``2^\lambda=\mu"$ 
    (so forcing with $\bbP$ collapses no cardinals, changes no 
    cofinalities, adds no new sequences of ordinals of length $<\lambda$, 
    and $\Vdash_{\bbP} ``\lambda^{<\lambda}=\lambda"$).
\sn
    \item[$(b)$]  We can find $\varphi,M_I$ (for ${}^{\kappa>}\!\lambda \subseteq I\subseteq {}^{\kappa\ge}\lambda$) as in \ref{4.2}$(*)$, so with $\|M_I\| = \lambda$, ($\tau$-models with $|\tau| = \kappa$ for simplicity) such that:
    \begin{enumerate}
        \item[$\oplus$]  There are, up to isomorphism, exactly 
        $\lambda^+$ models of the form $M_I$ \\
        (for ${}^{\kappa>}\!\lambda\subseteq I\subseteq {}^{\lambda\ge}\lambda$).
    \end{enumerate}
\sn
    \item[$(c)$]  In (b), there is a model $M$ such that $\|M\|=\lambda^+$ 
    and every model $M_I$ can be $(\pm\varphi)$-embedded into $M$.
\end{enumerate}
\end{lemma}

\sn
\begin{remark}\label{4.6A}
1) $M_{I}$ is essentially $(I^+,\lhd)$: the addition of level
predicates is immaterial, where $I^+$ extends $I$ ``nicely'' so that we can
let $a_\eta=\eta$ for $\eta\in I$.

\sn
2)  Clearly clause (c) also shows that weakening $\|M\| = \lambda$, even
when $\lambda^+ < 2^\lambda$, may make \ref{4.2} false.

\sn
3) In the proof of Lemma \ref{4.6}, the class of models isomorphic to some $N_j^*$ with $j < \lambda^+$ is not so nice. But the following class of models, which is reasonably well defined, will fail to satisfy the statement in \ref{4.2}(2) (in $\bfV^\bbP$).
\begin{enumerate}
    \item[$\boxplus$] $N \in K$ \underline{iff}
    \begin{enumerate}
        \item $N$ is a $\tau$-model.

        \item For some ordinal $\alpha$ and $S \subseteq {}^\kappa\alpha$, $N$ is isomorphic to $N_{I[\delta]}$, where $I = \{\eta \rest \zeta : \eta \in S,\ \zeta \leq \kappa\}$ and $N_{I[\delta]}$ is defined as in the proof below.
    \end{enumerate}
\end{enumerate}
\end{remark}

\sn
\begin{PROOF}{\ref{4.6}}
Let $\tau=\{R_{\zeta}:\zeta \le \kappa\} \cup \{<\}$ with $R_\zeta$
being a monadic predicate, and $<$ being a binary predicate. For a set $I$,
${}^{\kappa>}\!\lambda\subseteq I\subseteq {}^{\kappa\geq}\lambda$ let $N_I$
be the $\tau $-model:
\[
|N_I|=I,\quad R^{N_I}_\zeta=I\cap {}^\zeta\lambda,\quad
<^{N_I}=\{(\eta,\nu):
\eta,\nu\in I,\ \eta\lhd\nu\}, 
\]

and
\[
\varphi(\ldots,x_\zeta,\ldots)_{\zeta<\kappa} = \bigwedge\limits_{\zeta<
\xi<\kappa}\big(x_\zeta<x_\xi\ \wedge\ R_\zeta(x_\zeta)\big)\wedge (\exists
y) \big[R_\kappa(y)\ \wedge\ \textstyle\bigwedge\limits_{\zeta<\kappa} x_\zeta < y \big].
\]

\mn
Now we define the forcing notion $\bbP$. It is $\bbP_{\lambda^+}$, where
\[
\big\LL\bbP_i,\Name{\bbQ}_j:i\leq\lambda^+,\ j<\lambda^+ \big\RR
\]

\sn
is an iteration with support $<\lambda $, of $\lambda$-complete forcing
notions, where $\Name{\bbQ}_j$ is defined as follows.

\noindent 
For $j=0$ we add $\mu$ many Cohen subsets to $\lambda$:

\[
\bbQ_0 = \{f:f \text{ is a partial function from } \mu \text{ to }
\{0,1\},\ |\Dom(f)|<\lambda\},
\]

\mn
the order is the inclusion.

\noindent 
For $j>0$, we define $\Name{\bbQ}_j$ in $\bfV^{\bbP_j}$. Let
$\big\LL I(j,\alpha):\alpha<\alpha(j) \big\RR$ list, without 
repetition, all sets  $I \in \bfV^{\bbP_j}$ such that 
$^{\kappa>}\!\lambda \subseteq I \subseteq {}^{\kappa\ge} \lambda$. 
(Note that the interpretation of $^{\kappa\geq}\lambda$ does not 
change from $\bfV$ to $\bfV^{\bbP_j}$ (as $\kappa<\lambda$), 
but the family of such $I$-s increases.)  

Now
\begin{align*}
\Name{\bbQ}_j = \Big\{\bar{f} : &\ f = \LL f_\alpha:\alpha<\alpha(j)\RR,\
f_\alpha \text{ is a partial isomorphism}\\
  & \text{ from } N_{I(j,\alpha)} \text{ into } N_{({}^{\kappa\ge}\lambda)},\\
  &\ w({\bar f}) \defeq \{\alpha:f_\alpha \ne \varnothing\} \text{ has
    cardinality } <\lambda,\\
  &\ \Dom(f_\alpha) \text{ has the form }\bigcup\limits_{\beta<\gamma} 
{}^{\kappa \ge} \beta \cap N_{I(j,\alpha)} \text { for some }\gamma<\lambda;\\
  &\text{ and if } \alpha_1,\alpha_2 < \alpha(j) \text{ and }\eta_1,\eta_2\in
 {}^\kappa \lambda, \text{ and for every } \zeta<\kappa\\
  &\ f_{\alpha_1}(\eta_1 \rest \zeta),f_{\alpha_2}(\eta_2 \rest \zeta)
\text{ are well-defined and equal, then}\\
  &\ \eta_1 \in I(j,\alpha_1) \Leftrightarrow \eta_2 \in I(j,\alpha_2)\Big\}.
\end{align*}

The order is: 
\begin{align*}
{\bar f}^1 \le {\bar f}^2 \quad \text{ if and only if }
&\ \big(\forall\alpha < \alpha(j)\big)\big[f^1_\alpha \subseteq f^2_\alpha\big] \text{ and}\\
   & \text{ for all } \alpha<\beta<\alpha(j),\ f^1_\alpha \ne \varnothing
\wedge f^1_\beta \ne \varnothing \text{ implies}\\
   &\ \Rang(f^2_\alpha) \cap \Rang(f^2_\beta) = \Rang(f^1_\alpha)\cap
\Rang(f^1_\beta).
\end{align*}

\mn
Then, $\Name{\bbQ}_j$ is $\lambda$-complete and it satisfies the `$*_\lambda^\omega$' version of
$\lambda^+$-c.c. from \cite{Sh:80}\footnote{See more in \cite{Sh:546}, and much later in \cite{Sh:1036}.}, hence each
$\bbP_j$ satisfies the $\lambda^+$-c.c. (by \cite{Sh:80}).

Now the $\bbP_{j+1}$-name $\name{I}_j$ (interpreting it in 
$\bfV^{\bbP_{j+1}}$, we get $I^*_j$) is:
\begin{align*}
I^*_j = {}^{\kappa>}\!\lambda \cup \big\{\eta\in {}^\kappa\lambda: & \text{ for
some }{\bar f}\in\Name{G}_{\bbQ_j},\ \alpha<\alpha(j), \text{ and } \nu\in
N_{I(j,\alpha)},\\
  &\text{ we have } \lh(\nu) = \kappa \text{ and } f_{\alpha}(\nu) = \eta\big\}.
\end{align*}

\mn
This defines also $f^j_\alpha:I(j,\alpha) \to I^*_j$, which
is forced to be a $(\pm\varphi)$-embedding and also just an embedding.

So now we shall define, for every $I$ such that
${}^{\kappa>}\!\lambda\subseteq I \subseteq {}^{\kappa\ge}\lambda$, 
a $\tau$-model $M_I$: clearly $I$
belongs to some $\bfV^{\bbP_j}$. Let $j=j(I)$ be the first such $j$, and let
$\alpha=\alpha(I)$ be such that $I=I(j,\alpha)$. Let
$M_{I(j,\alpha)} = N_{I^*_j}$ (and $a_\rho=f^j_\alpha(\rho)$ for 
$\rho\in I(j,\alpha))$.

We leave the details to the reader. 
\end{PROOF}

\sn
On the other hand, consistently we may easily have a better result.

\begin{lemma}\label{4.7}
Suppose that, in the universe $\bfV$,
\[
\lambda=\cf(\lambda)=\lambda^\kappa=\lambda^{<\lambda},\quad \lambda<\mu=
\mu^\lambda.
\]

\sn
For some forcing notion $\bbP$:
\mn
\begin{enumerate}
    \item[$(a)$]  $\bbP$ is as in \ref{4.6}.
\sn
    \item[$(b)$]  In $\bfV^{\bbP}$, assume that 
    \begin{itemize}
        \item $\varphi$ and the function $I \mapsto \big(M_I,\LL{\bar a}^I_\eta : \eta \in {}^{\kappa >}\!\lambda \RR \big)$ are as required in clauses (a),(b),(c) of $(*)$ of \ref{4.2},

        \item $\zeta(*)<\mu$,

        \item each $N_\zeta$ (for $\zeta<\zeta(*)$) is a model in the relevant vocabulary,
        
        \item $\sum\limits_{\zeta<\zeta(*)}\|N_\zeta\|^\kappa<\mu$ (If the vocabulary is of cardinality $<\lambda$ and each predicate or relation symbol has finite arity, then requiring just\\ 
        $\sum\big\{\|N_\zeta\| : \zeta < \zeta(*) \big\} < \mu$ will suffice.)
    \end{itemize}
    \underline{Then} for some $I$, the model $M_I$ cannot be 
    $(\pm\varphi)$-embedded into any $N_\zeta$.
\sn
    \item[$(c)$]  Assume $\mu_1=\cf(\mu_1)$, $\lambda<\mu_1 \le \mu$ 
    and $\bfV \models (\forall\chi<\mu_1)[\chi^\lambda<\mu_1]$.
    \underline{Then} in $\bfV^{\bbP}$, if $\LL M_{I_i} : i < \mu_1\RR$ 
    are pairwise non-isomorphic, 
    ${}^{\kappa>}\!\lambda\subseteq I_i\subseteq {}^{\kappa\ge} \lambda$, 
    and $M_{I_i},{\bar a}_\eta^i$ (for $\eta\in I_i$) are as in \ref{4.2}$(*)$, \underline{then} $M_{I_i}$ is not embeddable into $M_{I_j}$ for some $i \ne j$.
\sn
    \item[$(d)$]  In $\bfV^\bbP$ we can find a sequence 
    $\LL I_\zeta:\zeta<\mu \RR$ ({with} 
    ${}^{\kappa>}\!\lambda\subseteq I_\zeta\subseteq {}^{\kappa\ge} \lambda$) such that no $M_{I_\zeta}$ is $(\pm\varphi)$-embeddable into another.
\end{enumerate}
\end{lemma}

\begin{PROOF}{\ref{4.7}}
$\bbP$ is $\bbQ_0$ from the proof of \ref{4.6}. Let $\bfF$ be the
generic function that is\\ $\bigcup\{f : f\in \Name{G}_{\bbQ_0}\}$: 
clearly it is a function from $\mu$ to $\{0,1\}$.  Now clause (a) is trivial. 

Next, concerning clause (b), we are given $\big\LL
N_\zeta : \zeta < \zeta(*) \big\RR$. Clearly for some $A\in \bfV$ of 
size smaller than $\mu$, we have $A \subseteq \mu$. To compute the isomorphism
types of $N_\zeta$ for $\zeta<\zeta(*)$, it is enough to know $\bfF \rest A$.  We can force by $\{f \in \bbQ_0:\Dom(f)\subseteq A\}$, then
$\bff \rest B$ for any $B \subseteq \lambda \setminus A$
of cardinality $\lambda$ (from $\bfV$) gives us an $I$ as required.

To prove clause (c) use a $\Delta$-system argument for the names of various
$M_I$-s, and similarly for (d).
\end{PROOF}

\newpage

\section {An Application for many models in $\lambda$}

\begin{discussion}\label{5.0}
Next we consider the following:

Assume $\lambda $ is regular, $(\forall\mu<\lambda)[\mu^{<\chi}<\lambda]$.
Let $\cU_{\alpha}$ (for $\alpha<\lambda$) be pairwise disjoint stationary subsets of $\{\delta < \lambda : \cf(\delta) = \chi\}$. 

For $A \subseteq \lambda$, let
\[
\cU_{A} = \bigcup\limits_{i\in A}\cU_i.
\]

\mn
We want to define $I_A$ such that 
${}^{\chi>}\!\lambda\subseteq I_A \subseteq {}^{\chi\ge}\lambda$ and
\[
A \nsubseteq B \Rightarrow M_{I_A} \not\cong M_{I_B}.
\]

\mn
We choose $\big\LL\LL M_{I_A}^i : i < \lambda\RR : A\subseteq
\lambda \big\RR$ with $M_{I_A} = \bigcup\limits_{i<\lambda} M^i_{I_A}$,
$\|M^i_{I_A}\|<\lambda$, $M^i_{I_A}$ increasing continuous.

\noindent
Of course, we have to strengthen the restrictions on $M_{I}$. For
$\eta\in I_A\cap {}^\chi\lambda$, let $\delta(\eta)\defeq\bigcup\{\eta(i)+1:
i<\chi\}$. We are specially interested in $\eta$ {which are}
strictly increasing converging to some $\delta(\eta) \in \cU_A$; we shall put
only such $\eta$-s in $I_A$. The decision whether $\eta\in I_A$ will be
done by induction on $\delta(\eta)$ for all sets $A$. Arriving to $\eta$,
we assume we know quite a lot on the isomorphism $f: M_{I_A}\rightarrow
M_{I_B}$: specifically, we know
\[
f \rest \bigcup\limits_{\alpha<\chi}{\bar a}_{\eta\rest\alpha},
\]

\mn
which we are trying to ``kill.'' We can assume $\delta(\eta)\notin 
\cU_B$ and $\delta$ belongs to a thin enough club of $\lambda$, and using all
this information we can ``compute" what to do.

\sn
(Note: though this is the
typical case, we do not always follow it.)
\end{discussion}

\begin{notation}\label{5.0A}
1) For an ordinal $\alpha$ and a regular $\theta \ge \aleph_0$, let
$\clH_{<\theta}(\alpha)$ be the smallest set $Y$ such that:
\mn
\begin{enumerate}
    \item[(i)]   $i\in Y$ for $i<\alpha$,
\sn
    \item[(ii)]  $x \in Y$ for $x\subseteq Y$ of cardinality $<\theta$.
\end{enumerate}
\mn
2)  We can agree that $\cM_{\lambda,\theta}(\alpha)$ from
\cite[2.1=Lf2]{Sh:E59} is interpretable in $(\clH_{<\theta}(\alpha),{\in})$ 
when $\alpha\ge \lambda$, and in particular its universe is a
definable subset of $\clH_{<\theta}(\alpha)$, and also $R$ is \mgrimes{defined to be}:

\begin{align*}
R = \Big\{ \big(\sigma^*,\LL t_i:i<\gamma_x\RR,x \big) : &\ x \in 
\cM_{\lambda,\theta}({}^{\theta>}\!\alpha),\ \sigma^* \text{ is a 
$\tau_{\lambda,\kappa}$-term,} \\
  &\ \theta \le \lambda\le \alpha, \text{ and } x=\sigma^*(\LL t_i : i < \gamma_x\RR) \Big\}.
\end{align*}

\mn
Similarly {for} $\cM_{\lambda,\theta}(I)$, where $I\subseteq {}^{\kappa>}
\lambda $ is interpretable in $(\clH_{<\chi}(\lambda^*),\in)$ if
$\lambda \le \lambda^*$, $\theta \le \chi$, and $\kappa \le \chi$.
\end{notation}

\noindent
The main theorem of this section (see \cite[1.4(1)=La11]{Sh:E59}) is:

\begin{theorem}\label{5.1}
$\dot I \dot E_{\pm\varphi}(\lambda,K)=2^\lambda$, provided that:
\mn
\begin{enumerate}
    \item[$(a)$]   $\lambda=\lambda^\chi$
\sn
    \item[$(b)$]  $\varphi = \varphi(\ldots, {\bar x}_\alpha, \ldots)_{\alpha<\chi}$ is a formula in the vocabulary $\tau_K$.
\sn
    \item[$(c)$]  For every $I$ such that ${}^{\chi>}\!\lambda\subseteq I\subseteq {}^{\chi\ge}\lambda$, we have a model $M_I\in K_\lambda$, 
    a function $f_I$, and ${\bar a}_\eta\in {}^{\chi\ge}|M_I|$ for 
    $\eta\in {}^{\chi>}\!\lambda$ with 
    $\lh({\bar a}_\eta) = \lh({\bar x}_{\lh(\eta)})$ such that:
    \begin{enumerate}
        \item[$(\alpha)$] For $\eta\in {}^\chi\lambda$ we have 
        $M_I \models \varphi(\ldots,{\bar a}_{\eta\rest \alpha},\ldots)$ \underline{if and only if}  $\eta \in I$.
\sn
        \item[$(\beta)$]  $f_I : M_I \to \cM_{\mu,\kappa}(I)$, where 
        $\mu \leq \lambda$ and $\kappa = \chi^+$.
\sn
        \item[$(\gamma)$] 
        If ${\bar b}_\alpha\in M_I$  is such that 
        $\lh({\bar x}_\alpha) = \lh({\bar b}_\alpha)$ for $\alpha<\chi$ and $f_I({\bar b}_\alpha)={\bar\sigma}_\alpha({\bar t}_\alpha)$ \underline{then}: 
        \begin{itemize}
            \item The truth value of $M_I\models\varphi[\ldots,{\bar b}_\alpha,\ldots]_{\alpha<\chi}$ can be computed from\\
            $\LL{\bar\sigma}_\alpha : \alpha < \chi\RR$ and 
            $\LL{\bar t}_\alpha : \alpha<\chi\RR$ (not just its quantifier-free~type in $I$) and from the truth values of statements of the form
            \[
                \big( \exists\nu \in I\cap {}^\chi\lambda \big) \Big[\textstyle\bigwedge\limits_{i<\chi} \nu \rest \epsilon_i = {\bar t}_{\beta_i}(\gamma_i) \rest \epsilon_i \Big]
            \]
            for $\alpha_i,\beta_i,\gamma_i,\epsilon_i<\chi$
            (i.e., in a way not depending on $I$ or $f_I$).
        
            \sn
            [We can weaken this.]
        \end{itemize}
    \end{enumerate}
\end{enumerate}
\end{theorem}

\noindent
We shall first prove \ref{5.1} under stronger assumptions.

\begin{fact}\label{5.2}
Suppose
\begin{enumerate}
    \item[$(*)$]  $\lambda=\lambda^{2^\chi }$ (so $\cf(\lambda)>\chi$) and $\chi\ge \kappa$.
\end{enumerate}

\sn
\underline{Then} there are $\big\{ (M^\alpha, \eta^\alpha) : \alpha < \alpha(*) \big\}$ such that:
\mn
\begin{enumerate}[(i)]
    \item   For every model $M$ with universe $\clH_{<\chi ^+}(\lambda)$
    such that $|\tau(M)| \le \chi$ (and, e.g., 
    $\tau\subseteq \clH_{<\chi^+}(\lambda)$), 
    for some $\alpha$, we have $M^\alpha\prec M$.
\sn
    \item  $\eta^\alpha\in {}^\chi\lambda$, $(\forall i < \chi)[\eta^\alpha \rest i \in M^\alpha]$, $\eta^\alpha\notin M^\alpha$, and 
    $\alpha \ne \beta \Rightarrow \eta^\alpha \ne\eta^\beta$.
\sn
    \item  For every $\beta < \alpha < \alpha(*)$ we have 
    $\{\eta^\alpha \rest i : i < \chi\} \not\subseteq M^\beta$.
\sn
    \item  For $\beta<\alpha$, if 
    $\{\eta^\beta \rest i : i < \chi\} \subseteq M^\alpha$ then $|M^\beta|\subseteq |M^\alpha|$.
\sn
    \item  $\|M^\alpha\|=\chi$.
\end{enumerate}
\end{fact}

\begin{PROOF}{\ref{5.2}}
By \ref{6.8} + \ref{6.10} below, with $\lambda,2^\chi,
\chi$ here standing for $\lambda,\chi(*),\theta$ there.

\mn
\textbf{Proof of \ref{5.1} from the Conclusion of \ref{5.2}.}

Without loss of generality, the universe of $M_I$ is $\lambda$ in \ref{5.1}.

We shall define, for every $A\subseteq\lambda$, a set $I[A]$ 
satisfying ${}^{\chi>}\!\lambda \subseteq I[A] \subseteq {}^{\chi\ge}
\lambda$; moreover,
\[
I[A] \setminus {}^{\chi>}\!\lambda \subseteq \big\{ \eta^\alpha : \alpha < \alpha(*) \big\}.
\]

\mn
For $\alpha < \alpha(*)$, let $\cU_\alpha = \big\{\eta \in {}^\chi\lambda :
\{\eta \rest i : i < \chi\} \subseteq M^\alpha \big\}$. 
We shall define, for every $A\subseteq\lambda$, the set 
$I[A] \cap \cU_\alpha$ by induction on $\alpha$ so that on the one hand,
those restrictions are compatible (i.e. in the end we can still 
define I[A] for each $A\subseteq\lambda$), and on the other hand they guarantee the non-$(\pm\varphi)$-embeddability.

For each $\alpha$, {we argue as follows.} 
Essentially, we decide whether $\eta^\alpha\in I[A]$,
assuming that $M^\alpha$ correctly ``guesses'' both a function 
$g : M_{I_1}\to M_{I_2}$ (where $I_\ell=I[A_\ell]$) and the set
$A_\ell\cap M^\alpha$ for $\ell=1,2$, and
we make our decision to prevent this.
\medskip

\noindent
\textbf{Case I}:  there are distinct subsets $A_1,A_2$ of
$\lambda$ and $I_1,I_2$ satisfying ${}^{\chi >}\!\lambda\subseteq I_\ell
\subseteq {}^{\chi\ge}\lambda$, a $(\pm\varphi)$-embedding $g$ of
$M_{I_1}$ into $M_{I_2}$, and
\[
M^\alpha \prec \big(\clH_{<\chi^+}(\lambda),\in,R,A_1,A_2,I_1,I_2,
M_{I_1},M_{I_2},f_{I_1},f_{I_2},g \big),
\]

\sn
where
\[
R = \Big\{ \big\{(0,\sigma_x,x),(1+i,t^x_i,x) \big\} : i < i_x \text{ and } 
x \text{ has the form } \sigma_x(\LL t^x_i:i<i_x\RR) \Big\}
\]

\mn
(we choose for each $x$ a unique such term $\sigma$), $I_2\cap 
\cU_\alpha\subseteq I_2\cap (\bigcup\limits_{\beta<\alpha}\cU_\beta)$,
and $I_1,I_2$ satisfy the restrictions we already have imposed on $I[A_1]$,
$I[A_2]$ respectively, for each $\beta<\alpha$. Computing the truth value of
$M_{I_2}\models\varphi [\ldots,f({\bar a}_{\eta^\alpha \rest i}),\ldots]_{i<\chi}$ according to clause \ref{5.1}(d) (assuming 
$I_2 \cap \cU^\alpha \subseteq \bigcup\limits_{\beta<\alpha}\cU^\beta$), 
we get $\bft^\alpha$.

\mn 
\underline{Then} we restrict:
\begin{enumerate}[(i)]
    \item  If $B\subseteq\lambda$ and $B\cap |M^\alpha| = A_2 \cap |M^\alpha|$
    then $I[B]\cap \big(\cU^\alpha \setminus \bigcup\limits_{\beta<\alpha} \cU^\beta \big) = \varnothing$.
\sn
    \item  If $B\subseteq\lambda$, $B \cap |M^\alpha| = A_1 \cap |M^\alpha|$, and $\bft^\alpha$ is true, then
    \[
        I[B] \cap \big(\cU^\alpha\setminus \textstyle\bigcup\limits_{\beta<\alpha} \cU^\beta\big) = \varnothing
    \]
    or just $\eta^\alpha \notin I[B]$.
\sn
    \item  If $B\subseteq\lambda$, $B\cap |M^\alpha| = A_1 \cap |M^\alpha|$, and $\bft^\alpha$ is false, then
    \[
        I[B] \cap \big( \cU^\alpha \setminus \textstyle\bigcup\limits_{\beta<\alpha} \cU^\beta \big) = 
        \{\eta^\alpha\}
    \]
    or just $\eta^\alpha \in I[B]$.
\end{enumerate}
\medskip

\noindent
\textbf{Case II:} Not Case I.

No restriction is imposed.
\end{PROOF}

\mn
The point of this is the two facts below, which should be clear.

\begin{fact}\label{5.3A}
The choice of $A_1,A_2,I_1,I_2,g$ is immaterial (any two candidates lead to
the same decision).
\end{fact}

\begin{PROOF}{\ref{5.3A}}
Use clause (d) of \ref{5.1}.
\end{PROOF}

\begin{fact}\label{5.3B}
The $M_{I[A]}$ (for $A\subseteq\lambda$) are pairwise non-isomorphic.
Moreover, for $A \ne B \subseteq \lambda$ there is no $(\pm\varphi)$-embedding
of $M_{I[A]}$ into $M_{I[B]}$.
\end{fact}

\begin{PROOF}{\ref{5.3B}}
By the choice of the $I[A]$-s and \ref{5.2}(i). 
\end{PROOF}

\medskip
\centerline {$* \qquad * \qquad *$}
\bigskip

Still, the assumption of \ref{5.2} is too strong: it does not cover all the
desirable cases, though it covers many of them. 
However, a statement weaker
than the conclusion of \ref{5.2} holds under weaker cardinality restrictions
and the proof of \ref{5.1} above works using it, thus we will finish the
proof of \ref{5.1}.

\begin{fact}\label{5.4}
Suppose $\lambda=\lambda^\chi $.

\underline{Then} there are $\big\{ (M^\alpha,A^\alpha_1,A^\alpha_2,\eta^\alpha) : \alpha < \alpha(*) \big\}$ such that:
\mn
\begin{enumerate}
    \item[$(*)$] 
    \begin{enumerate}[(i)]
        \item For every model $M$ with universe $\clH_{<\chi^+}(\lambda)$ such that
        $|\tau(M)| \le \chi$ and $\tau(M)\subseteq \clH_{<\chi^*}(\lambda)$
        (with arity of relations and functions finite) and sets 
        $A_1 \neq A_2\subseteq\lambda$, for some $\alpha<\alpha(*)$, we have $$(M^\alpha,A^\alpha_1,A^\alpha_2)\prec (M,A_1,A_2).$$

        \item $\eta^\alpha \in {}^\chi\lambda$, 
        $\{\eta^\alpha \rest i : i < \chi\} \subseteq |M^\alpha|$, 
        $\eta^\alpha\notin M^\alpha$, and $\alpha \ne \beta\ \Rightarrow\ \eta^\alpha\neq\eta^\beta$.
\sn
        \item For every $\beta<\alpha(*)$,  if $\{\eta^\alpha \rest i : i < \chi\} \subseteq M^\beta$, then $\alpha<\beta+2^\chi$. Furthermore, 
        $\alpha + 2^ \chi = \beta + 2^\chi $ implies 
        $A^\alpha_1\cap |M^\alpha| \ne A^\beta_2\cap |M^\alpha|$.

\sn
        \item For every $\beta<\alpha$, if $\{\eta^\beta \rest i : i < \chi\} \subseteq M^\alpha$, then $|M^\beta| \subseteq |M^\alpha|$.
\sn
        \item $\|M^\alpha\|=\chi$.
    \end{enumerate}
\end{enumerate}
\end{fact}

\begin{PROOF}{\ref{5.4}}
See \ref{6.20}.

\mn
\textbf{Proof of \ref{5.1}}: Should be clear,
We act as in the proof of \ref{5.1} from the conclusion of \ref{5.2}
but now we have to use the ``or just" version in (ii),(iii) there.
\end{PROOF}

\begin{conclusion}\label{5.6}
1)  If $T \subseteq T_1$ are complete first order theories, $T$ is in the
vocabulary $\tau$, $\kappa = \cf(\kappa) < \kappa(T)$ (hence 
$T$ is un-superstable), and $\lambda = \lambda^{\aleph_0} \ge |T_1|$, then $\dot{\bbI}_\tau(\lambda,T_1) = 2^\lambda$. (For more on $\dot{\bbI}_\tau$, see \cite{Sh:E59}.)

\sn
2)  Assume $\kappa=\cf(\kappa)$, $\Phi$ is proper and almost nice for
$K^\kappa_{\tr}$ (see \cite[1.7]{Sh:E59}), ${\bar\sigma}^i$ 
($i \le \kappa$) is a finite sequence of terms,
$\tau\subseteq\tau_\Phi$, $\varphi_i({\bar x},{\bar y})$ is first order in
$\cL[\tau]$, and for $\nu\in {}^i\lambda$, $\eta\in {}^\kappa\!\lambda$,
$\nu \lhd \eta$ we have that
$$\EM({}^\kappa\lambda,\Phi)\models \varphi_i \big({\bar\sigma}^\kappa_i(x_\eta),
{\bar\sigma}^{i+1}(x_{\eta \caret \LL\alpha\RR}) \big)$$ holds
\underline{if and only if} $\alpha=\eta(i)$. \underline{Then}
\[
\big|\big\{EM_\tau(S,\Phi)/{\cong} : {}^{\kappa >}\!\lambda 
\subseteq S \subseteq {}^{\kappa\ge}\lambda \big\} \big| = 2^\lambda.
\]
\end{conclusion}

\begin{PROOF}{\ref{5.6}}
1)  By \cite[1.10]{Sh:E59} there is a template $\Phi$ which is proper for
$K^\kappa_{\tr}$, as required in part (2).

\noindent
2) By \ref{5.1}.
\end{PROOF}

\begin{discussion}\label{5.7}
What about Theorem \ref{5.1} in the case we assume
only $\lambda=\lambda^{<\chi}$? There is some
information in \cite[Ch.VIII,\S2]{Sh:a}.

Of course, concerning un-superstable $T$, that is \ref{5.6},
more is done there: the assumption is just $\lambda > |T|$.
\end{discussion}

\begin{claim}\label{2.12}
In \ref{5.1}, we can restrict ourselves to $I$ such that $I^0_{\lambda,\chi}
\subseteq I\subseteq {}^{\chi\ge}\lambda$, where

\[
I^0_{\lambda,\chi}={}^{\chi>}\!\lambda \cup \{\eta\in {}^\chi\lambda:
\eta(i)=0 \text{ for every } i<\chi \text{ large enough}\}.
\]
\end{claim}

\begin{PROOF}{\ref{2.12}}
By renaming. 
\end{PROOF}

\newpage

\section {Black Boxes}

We try to give comprehensive treatment of black boxes:
quite a few few of them are useful in some contexts
and some parts are redone here, as explained in \S0,\S1.

Note that ``omitting countable types'' is a very useful device for building
models of cardinality $\aleph_0$ and $\aleph_1$. The generalization to
models of higher cardinality, $\lambda$ or $\lambda^+$, usually requires us
to increase the cardinality of the types to $\lambda$, and even so we may
encounter problems (see \cite{Sh:E60} and background there). 
Note that we do not look mainly at the omitting type 
theorem \emph{per se}, but at its applications.

Jensen defined square and proved existence in $\bfL$: in Facts \ref{6.1}--\ref{6.4}, 
we deal with slightly weaker related principles which can be proved in ZFC.  E.g. for $\lambda$ regular $>\aleph_1$, $\{\delta < \lambda^+ : \cf(\delta) < \lambda\}$ 
is the union of $\lambda$ sets, each has square (as defined there). You can skip them in first reading --- particularly \ref{6.1} (and later take the references on faith).

Then we deal with black boxes. In \ref{6.5} we give the simplest case:
$\lambda$ regular $>\aleph_0$, $\lambda=\lambda^{<\chi(*)}$. (Really,
$\lambda^{<\theta}=\lambda^{<\chi(*)}$ is almost the same.) In \ref{6.5}
we also assume ``$S\subseteq\{\delta<\lambda:\cf(\delta)=\theta\}$ is a good
stationary set.'' In \ref{6.7} we weaken this demand such that enough sets
$S$ as required exist (provably in ZFC!). The strength of the cardinality
hypothesis ($\lambda=\lambda^{<\chi(*)}$, $\lambda^{<\theta}=\lambda^{<
\chi(*)}$, $\lambda^\theta=\lambda^{<\chi(*)}$) vary the conclusion.
In \ref{6.6}--\ref{6.7A} we prepare the ground for replacing ``$\lambda$
regular'' by ``$\cf(\lambda)\geq\chi(*)$,'' which is done in \ref{6.9}.

As we noted in \S2, it is much nicer to deal with $({\olsi M}^\beta,
\eta^\beta)$, this is the first time we deal with $\eta^\beta$, i.e., for no
$\alpha<\beta$,
\[
\{\eta^\beta \rest i : i < \theta\} \subseteq \bigcup_{i<\theta} M^\alpha_i.
\]

\mn
In \ref{6.8}, \ref{6.10} (parallel to \ref{6.5}, \ref{6.9}, respectively) we
guarantee this, at the price of strengthening $\lambda^{<\theta}=
\lambda^{<\chi(*)}$ to 

\[
\lambda^{<\theta} =
\lambda^{\chi(1)},\ \chi(1) = \chi(*) + ({<}\chi(*))^\theta.
\]

\mn
Later, in \ref{6.20}, we draw the conclusion necessary for section 2 (in its
proof the function $h$, which may look redundant, plays the major
role). This (as well as \ref{6.8}, \ref{6.10}) exemplifies how those
principles are self propagating --- better ones follow from the old variant
(possibly with other parameters).

In \ref{6.11}--\ref{6.13} we deal with the black boxes when $\theta$ (the
length of the game) is $\aleph_0$. We use a generalization of the
$\Delta$-system lemma for trees and partition theorems on trees.\footnote{See 
    Rubin-Shelah \cite[\S4]{Sh:117}, \cite[Ch.XI]{Sh:b} = \cite[Ch.XI]{Sh:f}, 
    \cite[1.10=L1.7]{Sh:E62}, \cite[1.16=L1.15]{Sh:E62} and the proof of \ref{6.11A} here; see history there, and \ref{6.3}.} 
We get several versions of the black box  --- as the cardinality restriction
becomes more severe, we get a stronger principle.

It would be better if we can use, for a strong limit $\kappa > \aleph_0 = \cf(\kappa)$,
\begin{align*}
\kappa^{\aleph_0} = \sup\big\{\lambda : & \text{ for some } \kappa_n < \kappa 
\text{ and uniform ultrafilter}\\
  &\ D \text{ on } \omega,\ \cf\big(\textstyle\prod\limits_{n<\omega}\kappa_n/D\big) = \lambda\big\}.
\end{align*}

\sn
We know this for the uncountable cofinality case (see \cite{Sh:111} or
\cite{Sh:g}), but then there are other obstacles. Now \cite{Sh:355} gives a
partial remedy, but lately by \cite{Sh:400} there are many such cardinals.

In \ref{6.15}, \ref{6.16} we deal with the case $\cf(\lambda) \le \theta$.
Note that $\cf\big(\lambda^{<\chi(*)}\big) \ge \chi(*)$ is always true, so you
may wonder: why wouldn't we replace $\lambda$ by $\lambda^{<\chi(*)}$?
This is true in many applications: but is not true, for example, when
we want to construct structures with density character $\lambda$.

Several times, we use results quoted from \cite[\S2]{Sh:331}, but there are no dependency loops.
The pcf results quoted here are gathered in \cite[\S3]{Sh:E62}, so we will refer to it throughout in addition to quoting the original place.

We end with various remarks and exercises.

\subsection{On stationary sets}

\begin{fact}\label{6.1}
1)  If $\mu^\chi = \mu < \lambda \le 2^\mu$, $\chi$ and $\lambda$ are regular
uncountable cardinals, and $S \subseteq \{\delta < \lambda : \cf(\delta)=\chi\}$ is
a stationary set, \underline{then} there are a stationary set $W\subseteq\chi$ and
functions $h_a,h_b : \lambda \to \mu$ and $\LL S_\zeta:0<\zeta<
\lambda\RR$ such that:
\mn
\begin{enumerate}[(a)]
    \item  $S_\zeta\subseteq S$ is stationary.
\sn
    \item  $\xi \ne \zeta \Rightarrow S_\xi\cap S_\zeta=\varnothing$
\sn
    \item  If $\delta\in S_\xi$, then for some increasing continuous sequence $\LL\alpha_i : i < \chi\RR$ we have $\delta = \bigcup\limits_{i<\chi} \alpha_i$, $h_b(\alpha_i)=i$, $h_a(\alpha_i)\in\{\xi,0\}$, and the set 
$\{i<\chi : h_a(\alpha_i)=\xi\}$ is stationary (in fact, it is $W$).
\end{enumerate}
\mn
2) If in (1), a sequence $\LL C_\delta : \delta < \lambda,\ \cf(\delta) \leq \chi\RR$ satisfying
\[
(\forall\alpha\in C_\delta) \big[\alpha \text{ limit } \Rightarrow\ \alpha=\sup
(\alpha\cap C_\delta) \big]
\]

\sn
is given, where $C_\delta$ is a closed unbounded subset of $\delta$ of order type
$\cf(\delta)$, \underline{then} in the conclusion we can get also $S^*$ and 
$\LL C^*_\delta : \delta \in S^*\RR$ such that in addition to (a)--(c), we have:
\mn
\begin{enumerate}
    \item[(c)$'$]  In (c), we add $C_\delta=\{\alpha_i:i<\chi\}$.
\sn
    \item[(d)]  $\bigcup\limits_{0<\xi<\lambda} S_\xi \subseteq S^* \subseteq \bigcup\limits_{0<\xi<\lambda} S_\xi \cup \{\delta<\lambda:\cf(\delta)<\chi\}$
\sn
    \item[(e)]  $W \subseteq \chi$ is $(>\aleph_0)$-closed and stationary in cofinality $\aleph_0$, which means:
    \begin{enumerate}
        \item[$(i)$]  If $i<\chi$ is a limit ordinal such that $i = \sup(i\cap W)$ has cofinality $>\aleph_0$ then $i\in W$.
\sn
        \item[$(ii)$]  $\{i\in W : \cf(i)=\aleph_0\}$ is a 
        stationary\footnote{We can add `$\notin I$' if $I$ is any normal ideal on 
        $\{i<\chi:\cf(i)=\aleph_0\}$.} subset of $\chi$.
    \end{enumerate}
\sn
    \item[(f)]  for $\delta\in\bigcup\limits_{0<\xi<\lambda}S_\xi $ we have
    \[
        C^*_\delta = \{\alpha\in C_\delta:\otp(\alpha\cap C_\delta)=\sup(W\cap
        \otp(\alpha\cap C_\delta))\}
    \]

    \item[(g)]  $C^*_\delta$ is a club of $\delta$ included in $C_\delta$ for $\delta\in S^*$, and if $\delta(1)\in C^*_\delta$, $\delta \in S^*$, 
    $\delta \in \bigcup\limits_{0<\zeta<\lambda} S_\zeta$, 
    $\delta(1) = \sup(\delta(1)\cap C^*_\delta)$, and $\cf(\delta(1)) > \aleph_0$ \underline{then} $C^*_{\delta(1)} \subseteq C^*_\delta$,
\sn
    \item[(h)]  If $C$ is a closed unbounded subset of $\lambda$ and 
    $0 < \xi < \lambda$ \underline{then} the set\\ 
    $\{\delta \in S_\xi : C^*_\delta \subseteq C\}$ is stationary.
\end{enumerate}
\end{fact}

\begin{PROOF}{\ref{6.1}}
1) We can find $\{\LL h^1_\xi ,h^2_\xi\RR:\xi< \mu\}$ such that:
\begin{enumerate}
    \item  For every $\xi$ we have $h^1_\xi:\lambda\to\mu$ and $h^2_\xi:\lambda\to\mu$.
\sn
    \item  If $A\subseteq\lambda$, $|A|\leq\chi$, and $h^1,h^2 : A \to \mu$ then for some $\xi$, $h^1_\xi \rest A=h^1$ and $h^2_\xi \rest A=h^2$.
\end{enumerate}
This holds by Engelking-Karlowicz \cite{EK}.\footnote{See for example \cite[AP]{Sh:c}; on history see e.g. \cite[\S5]{Sh:430}}

\mn
2) For $\alpha < \lambda$, let $C_\alpha^\bullet$ be a closed unbounded subset of $\alpha$ of order type $\cf(\alpha)$. Now for each $\xi < \mu$ and for  
$a \subseteq \chi$ stationary, we ask whether for every $i < \lambda$, for some $j<\lambda$, we have
\begin{enumerate}
    \item[$(*)^{\xi,a}_{i,j}$]  The following subset of $\lambda$ is stationary:
    \begin{align*}
        S^{\xi ,a}_{i,j} = \big\{\delta \in S : &\text{ (i) if }\alpha \in C_\delta,\ \otp(\alpha\cap C_\delta) \notin a \text{ then } h^1_\xi(\alpha)=0, \\
        &\text{ (ii) if }\alpha\in C_\delta,\ \otp(\alpha\cap C_\delta) \in a
        \text{ then the } h^1_\xi(\alpha)\text{-th} \\
        &\text{ member of } C_\alpha \text{ belongs to } [i,j), \\
        &\text{ (iii) if } \alpha \in C_\delta \text{ then } h^2_\xi(\alpha)=\otp(\alpha \cap C_\delta)\}
    \end{align*}
\end{enumerate}

\begin{subfact}\label{6.1A}
For some $\xi < \mu$ and a stationary set $a\subseteq\chi$, for every $i < \lambda$, for some $j\in (i,\lambda)$, the statement $(*)^{\xi,a}_{i,j}$ holds.
\end{subfact}

\begin{PROOF}{\ref{6.1A}}
If not, then for every $\xi<\mu$ and a stationary
$a \subseteq \chi$, for some $i = i(\xi,a) < \lambda$, for every $j$ such that
$i(\xi,a) < j < \lambda$, there is a closed unbounded subset $C(\xi,a,i,j)$ 
of $\lambda$ disjoint from $S^{\xi,a}_{i,j}$.

Let
\[
i(*) = \bigcup \big\{i(\xi,a)+\omega:\xi<\mu\mbox{ and } a\subseteq \chi
\text{ is stationary}\big\}.
\]

\mn
Clearly $i(*) < \lambda$.

For $i(*) \le j < \lambda$, let
$$C(j) = \bigcap \big\{C(\xi,a,i(\xi,a),j): a\subseteq\chi \text{ is stationary and }
\xi<\mu \big\} \cap \big(i(*)+\omega,\lambda\big).$$ Clearly it is a closed 
unbounded subset of $\lambda$. 

Let
\[
C^* = \big\{\delta < \lambda : \delta > i(*) \text{ and } (\forall j<\delta) \big[\delta \in C(j)\big]\big\}.
\]

\mn
So $C^*$ is a closed unbounded subset of $\lambda$ as well. Let $C^+$ be the
set of accumulation points of $C^*$. Choose $\delta(*)\in C^+\cap S$,
and we shall define
\[
h^1 : C_{\delta(*)} \to \mu, \quad h^2 : C_{\delta(*)} \to \mu.
\]

\mn
For $\alpha\in C_{\delta(*)}$, let $h^0(\alpha)$ be:
\[
\min\big\{\gamma \in (0,\chi) : \text{the $\gamma^\tthh$ member of }
C_{\alpha}^\bullet \text{ is} >i(*)\big\}
\]

\mn
if $\alpha=\sup(C_{\delta(*)}\cap\alpha)>i(*)$, and zero otherwise.
Clearly the set
\[
\{\alpha\in C_{\delta(*)}:\ h^0(\alpha )=0\}
\]

\sn
is not stationary.  Now we can define $g : C_{\delta(*)} \to \delta(*)$ by:
\[
g(\alpha) \text{ is the $h^0(\alpha)^\tthh$ member of }C_{\alpha}.
\]

\mn
Note that $g$ is pressing down and $\{\alpha \in C_{\delta(*)} : g(\alpha) \le i(*)\}$
is not stationary. So (by the variant of Fodor's Lemma speaking on an ordinal of uncountable cofinality) for some $j < \sup(C_{\delta(*)}) = \delta(*)$, the set
\[
a \defeq \{\alpha\in C_{\delta(*)}\cap C^*:i(*)<g(\alpha)<j\}
\]

\mn
is a stationary subset of $\delta(*)$. Let $h^1:C_{\delta(*)}
\to\mu$ be
\[
h^1(\alpha) = \begin{cases}
0 & \text{ if } \otp (\alpha\cap C^\bullet_\delta)\notin a\\
h^0(\alpha) & \text{ if } \otp (\alpha\cap C^\bullet_\delta)\in a.
\end{cases}
\]

\mn
Let $h^2:C_{\delta(*)} \rightarrow\mu$ be $h^2(\alpha)=\otp(\alpha \cap
C_{\delta(*)})$.  By the choice of $\LL (h^1_\xi,h^2_\xi):\xi<\mu
\RR$, for some $\xi$ we have $h^1_\xi \rest C_{\delta(*)}= h^1$ and
$h^2_\xi \rest C_{\delta(*)} = h^2$. Easily, $\delta(*) \in 
S^{\xi,a}_{i,j}$ which is disjoint to $C(\xi,a,i(*),j)$, contradicting
$\delta(*) \in C^*$ by the definition of $C(j)$ and $C^*$.

So we have proved Subfact \ref{6.1A}. 
\end{PROOF}

\textbf{Continuing the proof of \ref{6.1}}:

Having chosen $\xi$ and $a$, we define an ordinal $i(\zeta) < \lambda$ by 
induction on $\zeta < \lambda$ such that $\LL i(\zeta) : \zeta < \lambda\RR$
is increasing continuous, $i(0)=0$, and $(*)^{\xi,a}_{i(\zeta),i(\zeta+1)}$ holds.

Now, for $\alpha<\lambda$ we define $h_a(\alpha)$ as follows: it is $\zeta$
if $h^1_\xi(\alpha) > 0$ and the $h^1_\xi(\alpha)^\tthh$ member of $C_\alpha^\bullet$
belongs to $\big[i(1+\zeta),i(1+\zeta+1)\big)$, and it is zero otherwise. Lastly,
let $h_b(\alpha)\defeq h^2_\xi(\alpha)$ and $W=a$ and

\[
\begin{array}{ll}
S_\zeta\defeq\big\{\delta\in S:&\text{(i)\ \ \ for }\alpha\in C_\delta,\ \otp(
\alpha \cap C_\delta)=h_b(\alpha),\\
  &\text{(ii)\ \ for } \alpha\in C_\delta,\ h_b(i)\in a\ \Rightarrow\
h_a(\alpha) = \zeta,\\
  &\text{(iii)\ for } \alpha\in C_\delta,\ h_b(i)\notin a\ \Rightarrow h_a(i)=0\big\} .
\end{array}
\]

\mn
Now, it is easy to check that $a$, $h_a$, $h_b$, and 
$\LL S_\zeta : 0 < \zeta < \lambda\RR$ are as required.

\mn
2) In the proof of \ref{6.1}(1) we shall now consider only
sets $a\subseteq\chi$ which satisfy the demand on $W$ from
\ref{6.1}(2)(e). (This makes a difference in the definition 
of $C(j)$ during the proof of Subfact \ref{6.1A}.) Also, in $(*)^{\xi,a}_{i,j}$
in the definition of $S^{\xi,a}_{i,j}$, we change (iii) to:
\mn
\begin{enumerate}
    \item[(iii)$'$]  If $\alpha\in C_\delta$ then $h^2_\xi(\alpha)$ codes the isomorphism type of (for example)
    \[
        \Big( C_\delta^\bullet \cup \textstyle\bigcup\limits_{\beta\in C_\delta} C_\beta,<,\alpha,C_\delta^\bullet,\big\{\LL i,\beta\RR : i \in C_\beta \big\} \Big).
    \]
\end{enumerate}
\mn
In the end, having chosen $\xi$ and $a$ we can define $C^*_\delta$ and
$S^*$ in the natural way. 
\end{PROOF}

\begin{fact}\label{6.2}
1) If $\lambda$ is regular $>2^\kappa$, $\kappa$ regular, 
$S \subseteq \{\delta < \lambda : \cf(\delta) = \kappa\}$ is stationary, 
and  $C^0_\delta$ is a club of $\delta$ of order type $\kappa$
($=\cf(\delta)$) for $\delta\in S$, \underline{then} we can find a club 
$c^*$ of $\kappa$ (see \ref{6.2A}(1) below) such that for $\delta \in S$, 
$$C_\delta = C_\delta^0[c^*] \defeq \{\alpha\in C^0_\delta :
\otp(C^0_\delta \cap \alpha) \in c^*\}.$$ 
It is a club of $\delta$, and:
\mn
\begin{enumerate}
    \item[$(*)$]  For every club $C \subseteq \lambda$, we have:
    \begin{enumerate}
        \item  If $\kappa>\aleph_0$ then $\{\delta\in S:C_\delta\subseteq C\}$ is stationary.
\sn
        \item  If $\kappa = \aleph_0$, then the set
        \[
            \big\{\delta\in S : (\forall\alpha,\beta \in C_\delta)[\alpha < \beta  \Rightarrow (\alpha,\beta) \cap C \ne \varnothing] \big\}
        \]
        is stationary.
    \end{enumerate}
\end{enumerate}
\mn
2) If $\lambda$ is a regular cardinal $>2^\kappa$, \underline{then} we can find 
$\big\LL\LL C^\zeta_\delta : \delta \in S_\zeta\RR : \zeta < 2^\kappa \big\RR$ 
such that:
\mn
\begin{enumerate}
    \item  $\bigcup\{S_\zeta : \zeta < 2^\kappa\} = \{\delta < \lambda : 
    \aleph_0 < \cf(\delta) \le \kappa\}$
\sn
    \item  $C^\zeta_\delta$ is a club of $\delta$ of order type $\cf(\delta)$.
\sn
    \item  If $\alpha\in S_\zeta$, $\cf(\alpha)>\theta>\aleph_0$, then
    \[
        \big\{\beta \in C^\zeta_\alpha : \cf(\beta) = \theta,\ \beta \in S_\zeta 
        \text{ and } C^\zeta_\beta \subseteq C^\zeta_\alpha \big\}
    \]
    is a stationary subset of $\alpha$.
\end{enumerate}
\mn
3) If $\lambda$ is regular and $2^\mu \ge \lambda > \mu^\kappa$
\underline{then} we can find $\big\LL\LL C^\zeta_\delta : \delta \in 
S_\zeta\RR : \zeta < \mu\big\RR$ such that:
\mn
\begin{enumerate}
    \item  $\bigcup\{S_\zeta : \zeta < 2^\kappa\} = \{\delta < \lambda : 
    \aleph_0 < \cf(\delta) \le \kappa\}$
\sn
    \item  $C^\zeta_\delta$ is a club of $\delta$ of order type $\cf(\delta)$.
\sn
    \item  If $\alpha \in S_\zeta$, $\beta \in C^\zeta_\alpha$, 
    $\cf(\beta) > \aleph_0$, then $\beta\in S_\zeta$ and 
    $C^\zeta_\beta \subseteq C^\zeta_\alpha$.
\sn
    \item  Moreover, if $\alpha,\beta\in S_\zeta$ and 
    $\beta \in C^\zeta_\alpha$ then
    \[
        \big\{\big(\otp(\gamma\cap C^\zeta_\beta),
        \otp(\gamma\cap C^\zeta_\alpha)\big) : \gamma \in C_\beta\big\}
    \]
    depends only on $\big(\otp(\beta\cap C_\alpha),\otp(C_\alpha)\big)$.
\end{enumerate}
\mn
4) We can, in clauses (1)$(*)$(a)-(b), replace ``stationary" by
``$\notin I$" for any normal ideal I on $\lambda$.
\end{fact}

\begin{remark}\label{6.2A}
1) Here a club $C$ of $\delta$, where $\cf(\delta) = \aleph_0$, just means an
unbounded subset of $\delta$.

\sn
2) In \ref{6.2}(1) instead of $2^\kappa$, the cardinal
\[
\min\big\{|\cF| : \cF \subseteq {}^\kappa\kappa\ \wedge\ 
(\forall g\in {}^\kappa\!\kappa)
(\exists f\in {\cF})(\forall \alpha < \kappa) \big[g(\alpha) < f(\alpha)\big]\big\}
\]
suffices.

\sn
3) In \ref{6.2}(1)$(*)$(b) above, it is equivalent to ask that
\[
\big\{\delta\in S : (\forall\alpha,\beta \in C_\delta)[\alpha<\beta\ \Rightarrow\ \otp((\alpha,\beta)\cap C)>\alpha]\big\}
\]
is stationary.
\end{remark}

\begin{PROOF}{\ref{6.2}}
1) If \ref{6.2}(1) fails, then for each club $c^*$ of $\kappa$
there is a club $C[c^*]$ of $\lambda$ exemplifying its failure. 
So $C^+ \defeq \bigcap\{C[c^*]:c^* \subseteq\kappa$ a club$\}$ is a 
club of $\lambda$. Choose a $\delta\in S$ which is an accumulation 
point of $C^+$, and get a contradiction easily.

\sn
2) Let $\lambda = \cf(\lambda)>2^\kappa$, and let $C_\alpha$ be a club
of $\alpha$ of order type $\cf(\alpha)$ for each limit $\alpha < \lambda$.
Without loss of generality
\[
\beta\in C_\alpha\ \wedge\ \beta>\sup(\beta\cap C_\alpha)\ \Rightarrow\
\beta\mbox{ is a successor ordinal.}
\]

For any sequence ${\bar c}=\LL c_\theta:\aleph_0<\theta=\cf(\theta)\le
\kappa\RR$ such that each $c_\theta$ is a club of $\theta$, for $\delta
\in S^* = \{\alpha<\lambda:\aleph_0<\cf(\alpha) \le \kappa\}$ we let:
\[
C^{\bar c}_{\delta} = \{\alpha\in C_\delta:\otp(C_\delta\cap
\alpha)\in c_{\cf(\delta)}\}.
\]

\mn
Now to define $S_{\bar c}$, we define the set $S_{\bar c} \cap \delta$ 
by induction on $\delta < \lambda$: the only problem is to define whether
$\alpha\in S_{\bar c} $ knowing $S_{\bar c}\cap\delta$. We stipulate
\[
\begin{array}{lcl}
\alpha\in S_{\bar c} &\text{\underline{if and only if}} &\text{(i)}\quad 
\aleph_0 < \cf(\alpha) \le \kappa\\
  &&\mbox{(ii)} \quad \text{If } \aleph_0 < \theta = \cf(\theta) < \cf(\alpha)\\
  &&\text{\ \ \ \ \quad then the set }\{\beta \in C^{\bar c}_{\alpha} :
\cf(\beta) = \theta,\ \beta\in S_{\bar c} \cap \alpha \}\\
  &&\mbox{\ \ \ \ \quad is stationary in }\alpha.
\end{array}
\]

\sn
Let $\LL{\bar c}^\zeta:\zeta<2^\kappa \RR$ list the possible
sequences ${\bar c}$, and let $S_\zeta=S_{{\bar c}^\zeta}$ and
$C^\zeta_\delta=C^{{\bar c}^\zeta}_\delta$. To finish, note that for each
$\delta<\lambda$ satisfying $\aleph_0<\cf(\delta)\leq\kappa$, we have 
$\delta\in S_\zeta$ for some $\zeta$.

\sn
3)  Combine the proof of (2) and of \ref{6.1}.

\sn
4) Similarly. 
\end{PROOF}

\noindent
We may remark

\begin{fact}\label{6.2B}
Suppose that $\lambda$ is a regular cardinal $>2^\kappa$, 
$\kappa = \cf(\kappa) > \aleph_0$, a set 
\[
S\subseteq \{\delta < \lambda : \cf(\delta) = \kappa\}
\]
is stationary, and $I$ is a normal ideal on $\lambda$ with $S\notin I$. 
If $I$ is $\lambda^+$-saturated (i.e. in the Boolean algebra 
$\clP(\lambda)/I$, there is no family of $\lambda^+$ pairwise 
disjoint elements), \underline{then} we can find $\LL C_\delta:
\delta\in S\RR$, $C_\delta$ a club of $\delta$ of order type 
$\cf(\delta)$, such that:
\mn
\begin{enumerate}
    \item[$(*)$]  For every club $C$ of $\lambda$ we have 
    $\{\delta\in S : C_\delta\setminus C \text{ is unbounded in } \delta\}\in I$.
\end{enumerate}
\end{fact}

\begin{PROOF}{\ref{6.2B}}
For $\delta \in S$, let $C'_\delta$ be a club of $\delta$ of
order type $\cf(\delta)$. Call ${\olsi C}=\LL C_\delta:\delta\in
S^* \RR$ (where $S^* \subseteq S \subseteq \lambda$ stationary, 
$S^* \notin I$, $C_\delta$ a club of $\delta$) $I$-\emph{large} 
if for every club $C$ of $\lambda$, the set
\[
\{\delta < \lambda : \delta \in S^* 
\text{ and $C_{\delta} \setminus C$ is bounded in } \delta\}
\]
does not belong to $I$.

\sn
We call $\olsi C$ $I$-\emph{full} if above $\{\delta\in S^* : C_\delta \setminus C$
unbounded in $\delta\}\in I$.

By \ref{6.2}(4), for every stationary $S'\subseteq S$ with $S'\notin I$, there is
a club $c^*$ of $\kappa$ such that $\LL C'_\delta[c^*]:\delta\in
S'\RR$ is $I$-large. 

\sn
Now note:
\begin{enumerate}
    \item[$(*)$]  If $\LL C_\delta:\delta\in S'\RR$ is $I$-large and
$S'\subseteq S$, then for some $S''\subseteq S'$ such that $S''\notin I$, 
$\LL C_\delta : \delta \in S''\RR$ is $I$-full (hence $S'' \notin I$).
\end{enumerate}

\mn
\textbf{Proof of $(*)$}: 

\noindent
Choose, by induction on $\alpha<\lambda^+$, a club
$C^\alpha$ of $\lambda$ such that:

\begin{enumerate}
    \item For $\beta<\alpha$, $C^\alpha\setminus C^\beta$ is bounded in $\lambda$.
\sn
    \item  If $\beta=\alpha+1$ then $A_\beta\setminus A_\alpha\in I^+$, where
    \[
        A_\gamma \defeq \{\delta \in S' : C_\delta \setminus C^\gamma 
        \text{ is unbounded in }\delta\}.
    \]
\end{enumerate}
\mn
As clearly
\[
\beta < \alpha\ \Rightarrow\ A_\beta \setminus A_\alpha \text{ is
bounded in }\lambda
\]
(by (a) and the definition of $A_\alpha ,A_\beta$) and as $I$ is
$\lambda^+$-saturated, clearly for some $\alpha $ we cannot define
$C^\alpha$. This cannot be true for $\alpha=0$ or a limit $\alpha$, so
necessarily $\alpha=\beta+1$. Now $S'\setminus A_\beta$ is not in $I$ as
${\olsi C}$ was assumed to be $I$-large. Check that $S''\defeq S'\setminus
A_\beta$ is as required.

\mn
Repeatedly using \ref{6.2}(4) and $(*)$, we get the conclusion.
\end{PROOF}

\begin{claim}\label{6.3}
Suppose $\lambda = \mu^+$, $\mu = \mu^\chi$, $\chi$ is a regular cardinal and
$$S \subseteq \{\delta < \lambda : \cf(\delta) = \chi\}$$ 
is stationary. \underline{Then} we can find $S^*$, 
$\LL C_\delta : \delta \in S^*\RR$, and
$\LL S_\xi : \xi < \lambda\RR$ such that:
\mn
\begin{enumerate}
    \item[$(a)$]  $\bigcup\limits_{\zeta<\mu} S_\zeta\subseteq S^* \subseteq 
    S \cup \{\delta < \lambda : \cf(\delta) < \chi\}$
\sn
    \item[$(b)$]  $S_\zeta \cap S$ is a stationary subset of $\lambda$ 
    for each $\zeta < \mu$.
\sn
    \item[$(c)$]  For $\alpha\in S^*$, $C_\alpha$ is a closed subset of 
    $\alpha$ of order type $\le \chi$. If $\alpha\in S^*$ is a limit then
    $C_\alpha$ is unbounded in $\alpha$ (so it is a club of $\alpha$).
\sn
    \item[$(d)$]  $\LL C_\alpha : \alpha \in S_\zeta \RR$ is a square on
$S_\zeta$; i.e. $S_\zeta$ is stationary in $\sup(S_\zeta)$ and:
    \begin{enumerate}
        \item[$(i)$]  $C_\alpha$ is a closed subset of $\alpha$, unbounded if $\alpha$ is limit.
\sn
        \item[$(ii)$]  If $\alpha\in S_\zeta$ and $\alpha(1) \in C_\alpha$ 
        then $\alpha(1) \in S_\zeta$ and $C_{\alpha(1)} = C_\alpha \cap \alpha(1)$.
    \end{enumerate}
\sn
    \item[$(e)$]  For each club $C$ of $\lambda$ and $\zeta<\mu$, we have 
    $C_\delta \subseteq C$ for some $\delta \in S_\zeta$.
\end{enumerate}
\end{claim}

\begin{PROOF}{\ref{6.3}}
Similar to the proof of \ref{6.1} (or see \cite{Sh:237e}). Alternatively, see \ref{6.4} below (using \ref{6.4B}(1) for clause (e)).
\end{PROOF}

\noindent
We shall use the following in \ref{6.13}.

\begin{claim}\label{6.3A}
Suppose $\lambda=\mu^+$, $\gamma$ a limit ordinal of cofinality $\chi$,
\[
h:\gamma\to \{\theta:\theta =1 \text{ or } \theta =
\cf(\theta) \le \mu\},
\]
$\mu = \mu^{|\gamma|}$, and 
$S \subseteq \{\delta < \lambda : \cf(\delta) = \chi\}$ is stationary. 
\underline{Then} we can find $S^*$, $\LL C_\delta : \delta \in S^* \RR$ and 
$\LL S_\zeta : \zeta < \lambda\RR$ such that:
\mn
\begin{enumerate}
    \item[$(a)$]  $\bigcup\limits_{\zeta<\lambda} S_\zeta\subseteq S^* 
    \subseteq \{\delta < \lambda : \cf(\delta) \leq \chi\}$
\sn
    \item[$(b)$]  $S_\zeta\cap S$ is stationary for each $\zeta < \lambda$.
\sn
    \item[$(c)$]  For $\delta\in S^*$,
    \begin{enumerate}
        \item[$(i)$]  $C_\delta$ is a club of $\delta$ of order type $\le \gamma$ and
\sn
        \item[$(ii)$]  $\otp(C_\delta) = \gamma$ iff $\delta \in S\cap S^*$,
\sn
        \item[$(iii)$]  $\alpha \in C_\delta \wedge \sup(C_\delta \cap \alpha) < \alpha\ \Rightarrow \ \alpha$ has cofinality $h[\otp(C_\delta\cap\alpha)]$.
    \end{enumerate}
\sn
    \item[$(d)$]  If $\delta\in S_\zeta$ and $\delta(1)$ is a limit ordinal 
    $\in C_\delta$ then $\delta(1)\in S_\zeta$ and\\ 
    $C_{\delta(1)} = C_\delta\cap\delta(1)$.
\sn
    \item[$(e)$] For each club $C$ of $\lambda$ and $\zeta<\lambda$, 
    for some $\delta \in S_\zeta$, $C_\delta\subseteq C$.
\end{enumerate}
\end{claim}

\begin{PROOF}{\ref{6.3A}}
Like \ref{6.3}.
\end{PROOF}

\begin{claim}\label{6.4}
1) Suppose $\lambda$ is regular $>\aleph_1$. \underline{Then} 
$\{\delta < \lambda^+ : \cf(\delta) < \lambda\}$ is a good stationary 
subset of $\lambda^+$.
(I.e. it is in $\check I[\lambda^+]$: see \cite[3.4=Lcd1.1]{Sh:E62} 
or \cite[0.6,0.7]{Sh:88r} or \ref{6.4A}(2) below.)

\noindent
2) Suppose $\lambda$ is regular $>\aleph_1$. \underline{Then} we can find
$\LL S_\zeta:\zeta<\lambda\RR$ such that:
\mn
\begin{enumerate}
    \item[$(a)$]  $\bigcup\limits_{\zeta<\lambda} S_\zeta=\{\alpha<\lambda^+:
\cf(\alpha)<\lambda\}$
\sn
    \item[$(b)$]  On each $S_\zeta$ there is a square (see clause \ref{6.3}$(d)$). Say it is $\LL C^\zeta_\alpha : \alpha \in S_\zeta\RR$ with 
    $|C^\zeta_\delta| < \lambda$.
\sn
    \item[$(c)$]  If $\delta(*) < \lambda$ and $\kappa = \cf(\kappa) < \lambda$
    \underline{then} for some $\zeta < \lambda$, for every club $C$ of $\lambda^+$, for some accumulation point $\delta$ of $C$, $\cf(\delta) = \kappa$ and $\otp(C^\zeta_\delta \cap C)$ is divisible by $\delta(*)$.
\sn
    \item[$(d)$]  If $\cf(\delta(*)) = \kappa$ as well, then we can add in the conclusion of $(c)$:
    \[
        C^\zeta_\delta \subseteq C \text{ and } \otp(C^\zeta_\delta) = \delta(*).
    \]
\end{enumerate}
\end{claim}

\begin{remark}\label{6.4A}
1) For $\lambda=\aleph_1$ the conclusion of \ref{6.4}(1), (2)(a),(b)
becomes totally trivial. But for $\delta < \omega_1$, it means something if we
add `$\{\alpha \in S_\zeta : \otp(C^\zeta_\alpha) = \delta\}$ is stationary, and
for every club $C$ of $\lambda$ the set
$\{\alpha \in S_\delta : \otp(C^\zeta_\alpha) = \delta,\ 
C^\zeta_\alpha \subseteq C\}$ is stationary.'
So \ref{6.4}(2)(c),(d) are not so trivial, but still true. Their proofs are
similar so we leave them to the reader (they are used only in \cite[2.7]{Sh:331}).

\noindent
2)  Recall that for a regular uncountable cardinal $\mu$, the family
$\check I[\mu]$ of good subsets of $\mu$ is the family of $S\subseteq \mu$ such
that there are a sequence $\bar{a}=\LL a_\alpha:\alpha<\lambda\RR$
and a club $C\subseteq\mu$ satisfying: 
\begin{itemize}
    \item $a_\alpha\subseteq \alpha$ is of order type $< \alpha$ when $\lambda$ is a successor cardinal.

    \item $\beta \in a_\alpha\ \Rightarrow\ a_\beta = a_\alpha \cap \beta$

    \item $(\forall \delta \in S \cap C)\big[\sup(a_\delta) = \delta\ \wedge\ \otp(a_\delta) = \cf(\delta)\big].$
\end{itemize}

\mn
We may say that the sequence $\bar{a}$ as above exemplifies that $S$ is
good; if $C = \mu$ we say ``explicitly exemplifies".
\end{remark}

\begin{PROOF}{\ref{6.4}}
Appears also in detail in \cite{Sh:351} (originally proved for this work but
as its appearance was delayed we put it there, too). Of course,

\noindent
1) follows from (2).

\noindent
2) Let $S = \{\alpha < \lambda^+ : \cf(\alpha) < \lambda\}$. For each 
$\alpha \in S$, choose $\bar{A}^\alpha$ such that:

\begin{enumerate}
    \item[$(\alpha)$]  $\bar{A}^\alpha = \LL A^\alpha_i : i < \lambda\RR$ 
    is an increasing continuous sequence of subsets of $\alpha$ of 
    cardinality $< \lambda$ such that 
    $\bigcup\limits_{i<\lambda} A^\alpha_i=\alpha\cap S$.
\sn
    \item[$(\beta)$]  If $\beta\in A^\alpha_i \cup \{\alpha\}$, $\beta$ 
    is a limit ordinal and $\cf(\beta)<\lambda$ (this actually follows 
    from the first two conditions), \underline{then} 
    $\beta = \sup(A^\alpha_i \cap \beta)$.
\sn
    \item[$(\gamma$)]  If $\beta\in A^\alpha_i \cup \{\alpha\}$ is limit and
    $\aleph_0 < \cf(\beta) < \lambda$ \underline{then} $A^\alpha_i$ contains 
    a club of $\beta$.
\sn
    \item[$(\delta)$]  $0\in A^\alpha_i$ and $\big(\beta \in S\ \wedge\ \beta + 1 \in A^\alpha_i \cup \{\alpha\} \big)\ \Rightarrow\ \beta \in A^\alpha_i$.
\sn
    \item[$(\eps)$]  The closure of $A^\alpha_i$ in $\alpha$ 
    (in the order topology) is included in $A^\alpha_{i+1}$.
\end{enumerate}
\mn
There are no problems with choosing $\bar{A}^\alpha$ as required.

We define $B^\alpha_i$ (for $i<\lambda$, $\alpha \in S$) by induction on
$\alpha$ as follows:
\[
B^\alpha_i = \begin{cases}
\closure(A^\alpha_i)\cap \alpha & \mbox{ if }
\cf(\alpha)\neq\aleph_1\\
\bigcap\big\{\bigcup\limits_{\beta\in C} B^\beta_i : C\mbox{ a club of }\alpha\big\}
&\mbox{ if } \cf(\alpha)=\aleph_1.
 \end{cases} 
\]

\mn
For $\zeta<\lambda$ we let:
\begin{align*}
S_\zeta=\big\{ \alpha\in S : &\mbox{ (i)\ \ } B^\alpha_{\zeta}\mbox{ is a closed subset of }\alpha,\\
&\mbox{ (ii)\ \ if }\beta\in B^\alpha_\zeta, \mbox{ then }B^\beta_\zeta =
B^\alpha_\zeta \cap\beta \mbox{ and}\\
&\mbox{ (iii) if }\alpha\mbox{ is limit, then }\alpha=\sup(B^\alpha_\zeta)
 \big\}
\end{align*}

\mn
and for $\alpha\in S_\zeta$ let $C^\zeta_\alpha=B^\alpha_\zeta$.

Now, demand (b) holds by the choice of $S_\zeta$. To prove clause (a) we
shall show that for any $\alpha\in S$, for some $\zeta<\lambda$, $\alpha\in
S_\zeta$; moreover we shall prove
\mn
\begin{enumerate}
    \item[$(*)^0_\alpha$]  $E_\alpha \defeq \{\zeta < \lambda : 
    \text{if $\cf(\zeta) = \aleph_1$ then }\alpha \in S_\zeta\}$ 
    contains a club of $\lambda$.
\end{enumerate}
\mn
For $\alpha \in S$ define $E^0_\alpha = \{\zeta < \lambda : \text{if } \cf(\zeta) 
= {\aleph_1} \text{ then }B^\alpha_\zeta =
\closure(A^\alpha_\zeta)\cap\alpha\}$. 
We shall prove by induction on $\alpha \in S$ that 
$E_\alpha\cap E^0_\alpha$ contains a club of $\lambda$,
and then we will choose such a club $E^1_\alpha$. 

Arriving to $\alpha$, let
\[
E^*_\alpha = \{\zeta < \lambda : \text{if } \beta\in A^\alpha_\zeta 
\text{ then } \zeta \in E^1_\beta \text{ and } A^\beta_\zeta = 
A^\alpha_\zeta \cap \beta\}.
\]
Clearly $E^*_\alpha$ is a club of $\lambda$. Let $\zeta\in E^*_\alpha$ 
and $\cf(\zeta) = \aleph_1$, and we shall prove that 
$\alpha \in S_\zeta \cap E_\alpha \cap E^0_\alpha$: clearly this will suffice. 
By the choice of $\zeta$ 
(and the definition of $E$) we have: if $\beta$ belongs 
to $A^\alpha_\zeta$ then $A^\beta_\zeta=
A^\alpha_\zeta\cap A$ and $B^\beta_\zeta = \closure(A^\beta_\zeta)\cap
\beta$, so 
\begin{enumerate}
    \item[$(*)_1$]  $\beta \in A^\alpha_\zeta\ \Rightarrow\ B^\beta_\zeta = \closure(A^\alpha_\zeta) \cap \beta$. 
\end{enumerate}
Let us check the three conditions for ``$\alpha\in S_\zeta$;"
this will suffice for clause (a) of the claim.

\mn
\textbf{Clause (i)}:  $B^\alpha_\zeta$ is a closed subset of $\alpha$.

If $\cf(\alpha) \ne \aleph_1$ then $B^\alpha_\zeta =
\closure(A^\alpha_\zeta)\cap \alpha$, hence necessarily it is a 
closed subset of $\alpha$.

If $\cf(\alpha)=\aleph_1$ then $B^\alpha_\zeta = \bigcap \big\{ \bigcup\limits_{\beta\in C} B^\beta_\zeta : C \text{ is a club of } \alpha \big\}$. 
Now, for any club $C$ of $\alpha$, $C\cap A^\alpha_\zeta$ is an unbounded 
subset of $\alpha$ (see clause $(\gamma)$ above). By $(*)_1$ above,
\[
\bigcup_{\beta\in C} B^\beta_\zeta \supseteq 
\bigcup_{\beta\in C\cap A^\alpha_\zeta}
B^\beta_\zeta = \closure(A^\alpha_\zeta)\cap\beta.
\]
To finish proving clause (i), it suffices to note that we have gotten
\begin{enumerate}
    \item[$(*)_2$]   $\alpha\in E^0_\zeta$.
\end{enumerate}

\sn
[Why?  If $\cf(\alpha)=\aleph_1$ see above, if $\cf(\alpha)\ne\aleph_1$
this is trivial.]

\mn
\textbf{Clause (ii)}:  If $\beta\in B^\alpha_\zeta$ then
$B^\beta_\zeta=B^\alpha_\zeta\cap\beta$.

We know that $B^\alpha_\zeta = \closure(A^\alpha_\zeta) \cap \alpha$ by
$(*)_2$ above. If $\beta\in A^\alpha_\zeta$ then (by $(*)_1$) we have
$B^\beta_\zeta =\closure(A^\alpha_\zeta)\cap \beta$, so we are done. So
assume $\beta\notin A^\alpha_\zeta$. Then by clause $(\eps)$,
necessarily:
\begin{enumerate}
    \item[$\odot$] If $\eps < \zeta$ then  $\beta > \sup(A^\alpha_\eps \cap \beta)$ and $\sup(A^\alpha_\eps \cap \beta) \in A^\alpha_{\eps+1} \subseteq A^\alpha_\zeta$.
\end{enumerate}

\mn
But $\beta\in B^\alpha_\zeta = \closure(A^\alpha_\zeta)$ by $(*)_2$,
hence together $A^\alpha_\zeta$ contains a club of $\beta$ and $\cf(\beta)=
\cf(\zeta)$, but $\cf(\zeta)=\aleph_1$, so $\cf(\beta)=\aleph_1$. Now, as in
the proof of clause (i), we get $B^\beta_\zeta=\bigcup\{B^\gamma_\zeta:
\gamma\in A^\alpha_\zeta\cap\beta\}$, so by the induction hypothesis we are
done.

\mn
\textbf{Clause (iii)}: If $\alpha$ is limit then $\alpha=\sup(A^\alpha_i)$.

By clause $(\beta)$ we know $A^\alpha_\zeta$ is unbounded in $\alpha$, but
$A^\alpha_\zeta\subseteq B^\alpha_\zeta$ (by $(*)_2$) and we are done.
\medskip

So we have finished proving $(*)^0_\alpha $ by
 induction on $\alpha$ hence clause (a) of the claim.

For proving (c) of \ref{6.4}(2), note that above, if $\alpha$ is limit, $C$
is a club of $\alpha$, $C\subseteq S$, and $|C| < \lambda$, then for every
$i$ large enough, $C \subseteq A^\alpha_i$ and even $C \subseteq B^\alpha_i$.

Now assume that the conclusion of (c) fails (for fixed $\delta(*)$ and
$\kappa$). Then for each $\zeta<\lambda$ we have a club $E^0_\zeta$
exemplifying it. Now, $E^0 \defeq\bigcap\limits_{\zeta<\lambda} E^0_\zeta$ is a
club of $\lambda^+$, hence for some $\delta\in E^0$, $\otp (E^0\cap\delta)$
is divisible by $\delta(*)$ and $\cf(\delta)=\kappa$. Choose an unbounded
in $\delta$ set $e \subseteq E^0\cap\delta$ of cardinality $<\lambda$ and 
order type divisible by $\delta(*)$. Then, for a final segment of 
$\zeta < \lambda$ we have $e \cap \delta \subseteq C^\zeta_\delta$.

Note that for any set $C_1$ of ordinals, $\otp(C_1)$ is divisible by $\delta
(*)$ if $C_1$ has an unbounded subset of order type divisible by
$\delta(*)$, so we get a contradiction because
by $(*)^0_{\delta(*)}$ for some $\zeta \in E_{\delta(*)}$ 
(so $\delta(*) \in S_\zeta$) by 
$E^0_\zeta \cap C^\zeta_\delta \supseteq E^0 \cap \delta \supseteq e$, 
$\sup(e)=\delta$ and $e$ has order type divisible by $\delta(*)$.

We are left with clause (d) of \ref{6.4}(2).  Fix $\kappa$, $\delta(*)$, 
and $\zeta$ as above, we may add $\le \lambda$ new sequences of the form
$\LL C_\alpha:\alpha \in S_\zeta\RR$ as long as each is a square. 
First assume that for every $\gamma$, $\beta<\lambda$, such that 
$\cf(\beta)=\kappa=\cf(\gamma)$, $\gamma$
divisible by $\delta(*)$ we have
\sn
\begin{enumerate}
\item[$(*)^3_{\beta ,\gamma}$]  There is a club $E_{\beta ,\gamma}$ of
$\lambda^+$ such that for no $\delta\in S_\zeta$ do we have
$\otp(C^\zeta_\delta)=\beta$ and $\otp (C^\zeta_{\delta} \cap 
E_{\beta,\gamma})=\gamma$.
\end{enumerate}
\sn
Then let
\[
E \defeq \bigcap \big\{E_{\beta,\gamma}:\gamma<\lambda,\ \beta<\lambda,\ \cf(\beta)
= \kappa = \cf(\gamma),\ \gamma \mbox{ divisible by } \delta(*) \big\}.
\]
Applying part (c) we get a contradiction.

So for some $\gamma$, $\beta < \lambda$, $\cf(\beta) = \kappa = \cf(\gamma)$,
$\gamma$ divisible by $\delta(*)$ and $(*)^3_{\beta,\gamma}$
fails. Also there is a club $E^*$ of $\lambda^+$ such that for every club
$E\subseteq E^*$ for some $\delta \in S_\zeta$, $\otp(C^\zeta_\delta) =
\beta$, $\otp(C^\zeta_\delta \cap E) = \gamma$ and $C^\zeta_\delta \cap E =
C^\zeta_\gamma \cap E^*$ (by \ref{6.4B} below). Let $e \subseteq \gamma = \sup(e)$
be closed and such that $\otp(e)=\delta(*)$ and
\[
\epsilon \in e \text{ is limit } \Rightarrow \epsilon=\sup
(e\cap\epsilon).
\]
We define ${}^*\!C^\zeta_\delta$ (for $\delta \in S_{\zeta}$) as follows:
if $\delta \notin E^*$ then
\[
{}^*\! C^\zeta_\delta \defeq C^\zeta_\delta\setminus (\max(\delta\cap
E^*)+1),
\]
if $\delta\in E^*$ and $\otp(C^\zeta_\delta\cap E^*) \in e \cup \{\gamma\}$ then
\[
{}^*\!C^\zeta_\delta=\{\alpha\in C^\zeta_\delta\cap E^*:\otp(\alpha\cap
C^\zeta_\delta\cap E^*)\in e\},
\]
and if $\delta\in E^*$, $\otp(C^\zeta_\delta\cap E^*) \notin e \cup \{\gamma\}$ let
\[
{}^*\!C^\zeta_\delta = C^\zeta_\delta \setminus \big(\max\! \big\{\alpha :
\otp(C^\zeta_\delta \cap E^* \cap\alpha) \in e \cup \{\gamma\} \big\}+1\big).
\]
One easily checks that (d) and square hold for $\LL {}^*\!C^\zeta_\delta :
\delta \in S_\zeta\RR$. So, we just have to add $\LL
{}^*\!C^\zeta_\delta : \delta \in S_\zeta\RR$ to $\{\LL
C^\zeta_\delta : \delta\in S_\zeta\RR : \zeta < \lambda\}$ for any
$\zeta,\delta(*),\kappa$ (for which we choose $\zeta$ and $E^*$).
\end{PROOF}

\begin{claim}\label{6.4B}
1) Assume that $\aleph_0 < \kappa = \cf(\kappa)$, $\kappa^+ < \lambda
= \cf(\lambda)$, $S \subseteq \{\delta < \lambda : \cf(\delta) = \kappa\}$ is
stationary, $C_\delta$ is a club of $\delta$ (for $\delta\in S$), and
$(\forall \delta \in S)\big[|C_\delta| = \kappa\big]$ (or at least
$\sup\limits_{\delta\in S} |C_\delta|^+<\lambda$).  \underline{Then} for some
club $E^* \subseteq\lambda$, for every club $E \subseteq E^*$, the set
$\{\delta\in S^*:C_\delta\cap E^* \subseteq E\}$ is stationary, where
\[
S^* \defeq \{\delta\in S:\delta \in \acc(E^*)\}.
\]

\mn
2) Assume that $\kappa = \cf(\kappa)$, $\kappa^+ < \lambda = \cf(\lambda)$, 
$S \subseteq \{\delta < \lambda : \cf(\delta) = \kappa\}$ is stationary,
$C_\delta$ is a club of $\delta$ (for $\delta\in S$), $\sup\limits_{\delta
\in S} |C_{\delta}|^+<\lambda$, $I_\delta$ is an ideal on $C_\delta$
{which includes} the bounded subsets, and for every club $E$ of 
$\lambda$, for stationarily many $\delta\in S$, we have 
$C_\delta\cap E \notin I_\delta$ (or $C_\delta \setminus E \in I_\delta$).

\underline{Then} for some club $E^*$ of $\lambda$, for every club $E\subseteq
E^*$ of $\lambda$ the set $\{\delta\in S^*:C_\delta\cap E^*
\subseteq E\}$ is stationary, where
\begin{align*}
S^* \defeq \big\{\delta \in S: &\ \delta\in \acc(E^*),\ \delta =
\sup(C_\delta \cap E^*) \text{ and} \\
  &\ C_\delta \cap E^* \notin I_\delta \text{ (or } C_\delta \setminus E^*
\in I_\delta) \big\}.
\end{align*}
\end{claim}

\begin{remark}\label{6.4C} 
This also was written in \cite{Sh:365}.
\end{remark}

\begin{PROOF}{\ref{6.4B}}
1) If not, choose by induction on $i<\mu\defeq\sup\limits_{\delta\in S}
(|C_\delta|^+)$ a club $E^*_i \subseteq\lambda $, decreasing with
$i$, $E^*_{i+1}$ exemplifies that $E^*_i$ is not as required, i.e.,
\[
\{\delta\in S^*(E^*_i):C_\delta\cap E^*_i\subseteq E^*_{i+1}\}=\varnothing.
\]
Now, $\acc(\bigcap\limits_{i<\mu} E^*_i\big)$ is a club of $\lambda$, so there is 
$\delta\in S \cap \acc\big(\bigcap\limits_{i<\mu} E^*_i\big)$. The sequence 
$\LL C_\delta \cap E^*_i : i < \mu\RR$ is necessarily strictly decreasing, 
and we get an easy contradiction.

\sn
2) Similarly.
\end{PROOF}

\subsection{Black Boxes: First round}

Now we turn to the main issue: black boxes.

\begin{lemma}
\label{6.5}
Suppose that $\lambda,\theta$ and $\chi(*)$ are regular cardinals and
$\lambda^\theta = \lambda^{<\chi(*)}$, $\theta < \chi(*) \le \lambda$, and 
a set $S\subseteq\{\delta < \lambda : \cf(\lambda) = \theta\}$ is stationary 
and in $\check I[\lambda]$.\footnote{If
$\theta = \aleph_0$ this holds trivially; see \cite[3.4=Lcd1.1]{Sh:E62}, 
\cite[0.6,0.7]{Sh:88r}, or just \ref{6.4A}(2).} 

\underline{Then} we can find
\[
\bfW = \big\{\big(\olsi M^\alpha,\eta^\alpha \big) : \alpha < \alpha(*) \big\}
\]
(pedantically, $\bfW$ is a sequence) and functions $\dot\zeta : \alpha(*) \to S$
and $h : \alpha(*) \to \lambda$ such that (so $\alpha,\beta < \alpha(*)$):
\mn
\begin{enumerate}
    \item[$(a0)$]  $h(\alpha)$ depends only on $\dot\zeta(\alpha)$, and $\dot\zeta$
is non-decreasing (but not necessarily strictly increasing).
\sn
    \item[$(a1)$]  We have:
    \begin{enumerate}
        \item[$(\alpha)$]   $\olsi M^\alpha = \LL M^\alpha_i : i \le \theta\RR$ 
        is an increasing continuous chain. ($\tau(M^\alpha_i)$, the vocabulary,
        may be increasing.) 
\sn
        \item[$(\beta)$]  Each $M^\alpha_i$ is an expansion of a submodel of $(\clH_{<\chi(*)}(\lambda),\in,<)$ belonging to $\clH_{<\chi(*)}(\lambda)$ and $M^\alpha_i$ is transitive (i.e. considering the ordinals as atoms, $x \in M^\alpha_i \Rightarrow x \subseteq M^\alpha_i$),
        so $M^\alpha_i$ necessarily has cardinality $<\chi(*)$. (Of course the 
        order means the order on the ordinals, and for transparency the vocabulary
        belongs to $\clH_{<\chi(*)}(\chi(*))$.)
\sn
        \item[$(\gamma)$]   $M^\alpha_i \cap \chi(*)$ is an ordinal, 
        $\chi(*) = \chi^+ \Rightarrow \chi+1 \subseteq M^\alpha_i$, and 
        $M^\alpha_i \in \clH_{<\chi(*)}(\eta^\alpha(i))$.
\sn
        \item[$(\delta)$]   $M^\alpha_i \cap \lambda\subseteq \eta^\alpha(i)$
\sn
        \item[($\eps)$]   $\LL M^\alpha_j:j\le i \RR\in M^\alpha_{i+1}$
\sn
        \item[$(\zeta)$]  $\eta^\alpha\in {}^\theta\!\lambda$ is increasing with limit $\dot\zeta(\alpha) \in S$ such that for $i < \theta$,\\ 
        $\eta^\alpha \rest(i+1) \in M^\alpha_{i+1}$.
    \end{enumerate}
\sn
    \item[$(a2)$]  In the following game, $\Game(\theta,\lambda,\chi(*),\bfW,h)$, Player I has no winning strategy. A play lasts $\theta$ moves. In the 
    $i^\tthh$ move Player I chooses a model $M_i \in \clH_{<\chi(*)}(\lambda)$, 
    and then Player II chooses $\gamma_i < \lambda$.  In the first move, Player I also chooses $\beta < \lambda$. In the end Player II wins the play if 
    $(\alpha) \Rightarrow (\beta)$, where: 
    \begin{enumerate}
        \item[$(\alpha)$]  The pair $\big(\LL M_i : i < \theta \RR,\LL \gamma_i : i < \theta \RR\big)$ satisfies the relevant demands on the 
            pair\footnote{So $\LL M_j : j \le i\RR$ is an increasing continuous chain, $M_i \cap \chi(*)$ an ordinal, $\chi(*) = \chi^+ \Rightarrow \chi+1\subseteq M_i$, $\LL M_\epsilon : \epsilon \le j \RR \in M_{j+1}$ and $\LL \gamma_\epsilon : \epsilon \le j \RR \in M_{j+1}$ for $j<i$, $M_i \in \clH_{<\chi(*)}(\gamma_i)$, and 
            $\LL\gamma_i : j \le i\RR \in M_{i+1}$.}
        $(\olsi M^i \!\rest \theta,\eta^\alpha)$ in clause $(a1)$.
\sn
        \item[$(\beta)$]  For some $\alpha<\alpha(*)$, 
        $\eta^\alpha = \LL\gamma_i : i < \theta\RR$, $M_i = M^\alpha_i$ for $i<\theta$, and $h(\alpha) = \beta$.
    \end{enumerate}
\sn
    \item[$(b0)$]  $\eta^\alpha\ne\eta^\beta$ for $\alpha\ne\beta$.
\sn
    \item[$(b1)$]  If $\{\eta^\alpha \rest i : i < \theta\} \subseteq 
    M^\beta_\theta$ \underline{then} 
    \begin{enumerate}[$\bullet_1$]
        \item $\dot\zeta(\alpha) \leq \dot\zeta(\beta)$

        \item $x \in M^\alpha_\theta \Rightarrow x \in M^\beta_\theta$

        \item $\alpha + (<\chi(*))^\theta = \beta + (<\chi(*))^\theta$ (see \ref{6.5A}(2) below).
    \end{enumerate}
    
\sn
    \item[$(b2)$]  If in addition $\lambda^{<\theta} = \lambda^{<\chi(*)}$, 
    \underline{then} for every $\alpha < \alpha(*)$ and $i < \theta$, there is 
    $j < \theta$ such that $\eta^\alpha \rest j \in M^\beta_\theta$ implies $M^\alpha_i \in M^\beta_\theta$ (hence $M^\alpha_i \subseteq M^\beta_\theta$).
\sn
    \item[$(b3)$]  If $\lambda = \lambda^{<\chi(*)}$ and 
    $\eta^\alpha \rest (i+1) \in M^\beta_j$ \underline{then} $M^\alpha_i \in M^\beta_j$ (and hence $M^\alpha_i \subseteq M^\beta_j$, so $x \in M^\alpha_i \Rightarrow x \in M^\beta_j$) and
    \[
        \eta^\alpha \rest i \ne \eta^\beta \rest i\ \Rightarrow\ 
        \eta^\alpha(i)\ne \eta^\beta(i).
    \]
\end{enumerate}
\end{lemma}

\begin{remark}\label{6.5A}
1)  If $\bfW$ (with $\dot\zeta,h,\lambda,\theta,\chi(*)$) 
satisfies (a0), (a1), (a2), (b0), (b1) we call it a \emph{barrier}.

\sn
2) Remember, $({<}\chi)^\theta\defeq\sum\limits_{\mu<\chi} \mu^\theta$.

\sn
3) The existence of a good stationary set 
$S \subseteq \{\delta < \lambda : \cf(\delta) = \theta\}$ follows, 
for example, from $\lambda = \lambda^{<\theta}$ (see \cite[3.4=Lcd1.1]{Sh:E62} or \cite[0.6,0.7]{Sh:88r}) and from ``$\lambda$ is the successor of a regular cardinal and $\lambda>{\theta^+}$."  But see \ref{6.7}(1),(2),(3).

\sn
4) Compare the proof below with \cite[Lemma 1.13,pg.49]{Sh:227} and
\cite{Sh:140}.
\end{remark}

\begin{PROOF}{\ref{6.5}}
First assume $\lambda = \lambda^{<\chi(*)}$.

Let $\LL S_\gamma : \gamma < \lambda\RR$ be a sequence of pairwise
disjoint stationary subsets of $S$ such that $S = \bigcup\limits_{\gamma<\lambda}
S_\gamma$, and without loss of generality $\gamma<\min(S_\gamma)$.
We define $h^* : S \to \lambda$ by $h^*(\alpha) =$ ``the
unique $\gamma$ such that $\alpha \in S_\gamma$", and below we shall
let $h(\alpha) \defeq h^*(\dot\zeta(\alpha))$.

Let $\cd = \cd_{\lambda,\chi(*)}$ be a one-to-one function from
$\clH_{<\chi(*)}(\lambda)$ onto $\lambda$ such that
$\cd(\LL\alpha,\beta\RR)$ is an ordinal 
$$\max(\alpha,\beta) < \cd(\LL\alpha,\beta\RR) < \max(|\alpha+\beta|^+,\omega),$$
and $x \in \clH_{<\chi(*)}(\cd(x))$ for every relevant $x$. For $\xi \in S$ let:
\begin{enumerate}
    \item[$(*)_1$] 
    \begin{enumerate}
        \item $\bfW_\xi^0 \defeq$
        \begin{align*}
            \Big\{ \big(\olsi M,\eta\big) : &\text{ the pair } (\olsi M,\eta)
            \text{ satisfies (a1) of \ref{6.5}, } \\
            &\ \sup\{\eta(i) : i < \theta\} = \xi, \text{ and for every } i < \theta,\text{ for}\\
            & \text{ some } y \in \clH_{<\chi(*)}(\lambda),\ \eta(i) = \cd\big(\LL\olsi M \rest (i+1),\eta \rest i,y\RR\big)\Big\}.
        \end{align*}

        \item $\bfW = \bigcup\{\bfW_\xi^0 : \xi \in S\}$
    \end{enumerate}
\end{enumerate}
Below, we shall choose $\big\LL (\olsi M^\alpha\!,\eta^\alpha) : \alpha < \alpha(*) \big\RR$ listing $\bfW$.

So (a0), (a1), (b0), (b3) (hence (b2)) should be clear.

We can choose $\big\LL(\olsi M^\alpha\!,\eta^\alpha) : \alpha < \alpha(*) \big\RR$
an enumeration of $\bigcup\limits_{\xi\in S} \bfW^0_\xi$ to 
satisfy (b1) (and $\dot\zeta(\alpha) = \sup\rang(\eta^\alpha)$, of
course) because:
\begin{enumerate}
    \item[$(*)_2$] If $(\olsi M^*,\eta^*) \in \bigcup\limits_\xi \bfW^0_\xi$ then
    \[
        \big|\big\{\eta\in {}^\theta\!\lambda : \{\eta \rest i : i < \theta\} \subseteq
        M^*_\theta \big\} \big| \le \|M^*_\theta\|^\theta \le \big({<}\chi(*)\big)^\theta.
    \]
\end{enumerate}
Clearly $(*)_2$ holds, but why does it suffice for choosing our 
$\big\LL (\olsi M^\alpha\!,\eta^\alpha) : \alpha < \alpha(*) \big\RR$? 
\begin{enumerate}
    \item[$(*)_{2.1}$] We define the partial order $\leq_\bfW$ on $\bfW$ by
    $$(\olsi M,\eta) \leq_\bfW (\olsi M',\eta') \text{ \underline{iff} } M_\theta \subseteq M_\theta'.$$
\end{enumerate}
For each $\xi \in S$, try to choose a sequence $\bfx_{\xi,\gamma} = \big\LL 
(\olsi M^{\xi,\gamma},\eta^{\xi,\gamma}) : \alpha < \alpha_\gamma\big\RR$ 
by induction on the order $\gamma < \|\bfW_\xi^0\|^+$, so it will be 
$\lhd$-increasing with $\gamma$ such that:
\begin{enumerate}
    \item[$(*)_{2.2}$] 
    \begin{enumerate}
        \item $\big(\olsi M^{\xi,\alpha},\eta^{\xi,\alpha}\big) \in \bfW_\xi^0$ 
        for $\alpha < \alpha_\gamma$.

        \item If $(\olsi M,\eta) \in \bfW_\xi^0$ and 
        $(\olsi M,\eta) \leq_\xi \big(\olsi M^{\xi,\alpha},\eta^{\xi,\alpha}\big)$
        for some $\alpha < \gamma_\alpha$ \underline{then} 
        $(\olsi M,\eta) = \big(\olsi M^{\xi,\beta},\eta^{\xi,\beta}\big)$ 
        for some $\beta < \alpha_\gamma$.
    \end{enumerate}
\end{enumerate}
How do we carry the induction? For $\gamma = 0$ let $\alpha_\gamma = 0$; 
also, for $\gamma$ a limit ordinal the choice of $\bfx_{\xi,\gamma}$ is 
obvious. For $\gamma = \gamma_1 + 1$ ($= \gamma(1) + 1$), if 
$$\big\{ \big(\olsi M^{\xi,\alpha},\eta^{\xi,\alpha}\big) : \alpha < \alpha_{\gamma_1} \big\} = \bfW_\xi^0$$ 
then we stop. Otherwise, choose 
$$\big(\olsi N^{\gamma_1},\eta^{\gamma_1}\big) \in \bfW_\xi^0 \setminus \big\{ \big(\olsi M^{\xi,\alpha},\eta^{\xi,\alpha}\big) : \alpha < \alpha_{\gamma_1} \big\}$$
and let
\begin{align*}
    \bfW_{\xi,\gamma_1} = \Big\{(\olsi M,\eta) \in \bfW_\xi^0 : &\ (\olsi M,\eta) \leq_\xi (\olsi N^{\gamma_1},\eta^{\gamma_1}), \text{ but}\\  
&\ (\olsi M,\eta) \notin \big\{ \big(\olsi M^{\xi,\alpha}, \eta^{\xi,\alpha}\big) : \alpha < \alpha_{\gamma_1}\big\}\Big\},
\end{align*}
so $\bfx_\gamma$ is defined by letting 
$\alpha_\gamma \defeq \alpha_{\gamma(1)} + \|\bfW_{\xi,\gamma(1)}\|$ and 
$$\big\LL \big(\olsi M^{\xi,\alpha},\eta^{\xi,\alpha}\big) : \alpha \in[\alpha_{\gamma(1)}, \alpha_\gamma) \big\RR$$
list the elements of $\bfW_{\xi,\gamma(1)}$.

So for some $\gamma[\xi] < |\bfW_\xi^0|^+$, $\bfx_{\xi,\gamma[\xi]}$ lists $\bfW_\xi^0$. Lastly, we choose $\alpha(*) = \sum\limits_{\xi \in S} \gamma[\xi]$, and $(\olsi M^\alpha\!,\eta^\alpha) =(\olsi M^{\xi,\beta}\!,\eta^{\xi,\beta})$ when $\alpha = \sum\{\gamma[\xi'] : \xi' < \xi\} + \beta$ and $\beta < \gamma[\xi]$. 

This, in fact, defines the function $\dot\zeta$ as follows: we have
$\dot\zeta(\alpha)=\xi$ \underline{if and only if} 
$(\olsi M^\alpha\!,\eta^\alpha)\in \bfW^0_\xi$.

We are left with proving (a2). Let \textsf{G} be a strategy for Player I.

Let $\LL C_\delta:\delta<\lambda\RR$ exemplify ``$S$ is a good
stationary subset of $\lambda$" (see \ref{6.4A}(2)) and let 
$R = \{(i,\alpha) : i \in C_\alpha,\ \alpha < \lambda\}$.

Let $\LL \cA_i : i < \lambda\RR$ be a representation of the model
\[
\cA = \big(\clH_{<\chi(*)}(\lambda),\in,G,R,\cd \big),
\]

\mn
i.e. it is increasing continuous, $\|\cA_i\| < \lambda$, and
$\bigcup\limits_i \cA_i = \cA$. Without loss of generality
$\cA_i \prec \cA$ and $|\cA_i|\cap\lambda$ is an
ordinal for $i < \lambda$.

Let \textsf{G} ``tell" Player I to choose $\beta^* < \lambda$ in
his first move. So there is a $\delta \in S_{\beta^*}$ (hence
$\delta > \beta^*$: see the beginning of the proof) such that $|\cA_\delta|\cap\lambda=\delta$.
Now, necessarily $C_\delta\cap\alpha\in \cA_\delta$ for $\alpha <
\delta$.  Let $\{\alpha_i:i<\cf(\delta)\}$ list $C_\delta$ in
increasing order.

Lastly, by induction on $i$, we choose $M_i,\eta(i)$ as follows:
\[
\eta(i) = \cd\Big(\Big\LL \LL M_j : j \le i\RR,\LL\eta(j):
j < i\RR, \LL\alpha_j : j < i\RR \Big\RR\Big),
\]

and $M_i$ is what the strategy \textsf{G} ``tells" Player I to choose in his
$i^\tthh$ move if Player II has chosen $\LL\eta(j) : j < i\RR$ so far.

Now for each $i < \theta$, the sequences $\LL M_j : j \le i\RR$,
$\LL\eta(j) : j < i\RR$ are definable in $\cA_\delta$ with 
$\LL\alpha_j : j \le i\RR$ as the only parameter, hence they belong to 
$\cA_\delta$. So $$\sup\{\eta(j) : j < \theta\} \le \delta.$$ However, by the choice
of $\eta(i)$ (and $\cd$), $\eta(i) \ge \sup\{\alpha_j : j < i\}$ and hence\\
$\sup\{\eta(j) : j < \theta\}$ is necessarily $\delta$. Now check.
\medskip

We {have} finish{ed} the proof, but only by including the assumption 
$\lambda = \lambda^{<\chi(*)}$.  The case 
$\lambda < \lambda^{<\theta} = \lambda^{<\chi(*)}$ is similar. For a
set $A\subseteq\theta$ of cardinality $\theta$ we let $\cd^A=\cd^A_{\lambda,
\chi(*)}$ be a one-to-one function from $\clH_{<\chi(*)}(\lambda)$
onto $A_\lambda$, where
\[
A_\lambda = \{h:h \text{ is a function from } A \text{ to } \lambda\}.
\]
We strengthen (b2) to 
\mn
\begin{enumerate}
    \item[$(b2)'$]  Let $A_i \defeq \{\cd(i,j):j<\theta\}$ for 
    $i\in [1,\theta)$ and $A_0 \defeq \theta \setminus 
    \bigcup\{A_{1+i} : i < \theta\}$, so $\LL A_i : i < \theta \RR$ is a 
    sequence of pairwise disjoint subsets of $\theta$, each of cardinality 
    $\theta$, with $\min(A_i) \ge i$, and we have 
    \begin{enumerate}
        \item[$(*)$]  $\eta^\alpha \rest A_{i} = \cd^{A_i}\big({\olsi M^\alpha\!
  \rest i}, \eta^\alpha \rest i\big)$. \qedhere
    \end{enumerate}
\end{enumerate}
\end{PROOF}

\medskip
\centerline {$* \qquad * \qquad *$}
\bigskip

What can we do when $S$ is not good?  As we say in \ref{6.5A}(3), in many
cases a good $S$ exists (note that for singular $\lambda$ we will not have
one).

The following rectifies the situation in the other cases (but is interesting
mainly for $\lambda$ singular). We shall, for a regular cardinal $\lambda$,
remove this assumption in \ref{6.7}(1)--(3), while \ref{6.7A} helps for
singular $\lambda $. (This is carried in \ref{6.9}).

\begin{definition}
\label{6.6}
Let $\partial$ be an ordinal greater than 0, and for $\alpha < \partial$ 
let $\kappa_\alpha$ be a regular uncountable cardinal and 
$S_\alpha \subseteq \{\delta < \kappa_\alpha : \cf(\delta) = \theta )\}$ 
be a stationary set. Assume $\theta$, $\chi$ are regular cardinals such 
that for every $\alpha < \partial$ we have $\theta < \chi \leq \kappa_\alpha$. 
Let $\bar{S} = \LL S_\alpha : \alpha < \partial\RR$,
$\bar{\kappa} = \LL\kappa_\alpha : \alpha < \partial\RR$.
If $\partial = 1$ we may write $S_0,\kappa_0$.  

\sn
1) We say that $\bar{S}$ is \emph{good for} $(\bar{\kappa},\theta,\chi)$ 
when\footnote{Note that we can compute $\partial$ from $\bar \kappa$.} for every large enough $\mu$ and model $\cA$ expanding
$(\clH_{<\chi}(\mu),\in)$ with $|\tau(\cA)| \le \aleph_0$, there are 
$M_i$ (for $i < \theta$) such that:
\mn
\begin{itemize}
    \item   $M_i \prec \cA$ and $\bar S \in M_i$.  
\sn
    \item   $\LL M_j : j \le i\RR \in M_{i+1}$, $\|M_i\| < \chi$, 
    $M_i \cap \chi \in \chi$, and
    $\chi = \chi^+_1 \Rightarrow \chi_1 + 1 \subseteq M_i$. 
\sn
    \item   $\alpha < \partial$, $\alpha \in \bigcup\limits_{j<\theta} M_j$ implies that $\sup\!\big(\kappa_\alpha \cap (\bigcup\limits_{j<\theta} M_j)\big)$ belongs to $S_\alpha$.
\end{itemize}

\noindent
2) If $\bar \kappa$ is constant (i.e. $i < \partial \Rightarrow \kappa_i = \kappa$) then we may say $\bar S$ is good for $(\kappa,\partial,\theta,\chi)$. We may omit $\partial$ if $\partial = \kappa$.

\sn
3) If $\partial=1$, we  may write $S_0,\kappa_0$ instead of $\bar{S},\bar{\kappa}$.
If $\partial < \chi$ then we can  demand $\partial \subseteq M_0$.
\end{definition}

\begin{definition}\label{6.6C}
For regular uncountable cardinal $\lambda$ and regular $\theta < \lambda$,
let $\check J_\theta[\lambda]$ be the family of subsets $S$ of $\lambda$
such that $\{\delta \in S : \cf(\delta) = \theta\}$ is not good for $(\lambda,\theta,\lambda)$ (i.e. for $(\lambda,\lambda,\theta,\lambda)$).
\end{definition}

\begin{claim}\label{6.7}
Assume $\theta = \cf(\theta) < \chi = \cf(\chi) \le \kappa = \cf(\kappa)$.

\sn
1) \underline{Then} $\{\delta < \kappa : \cf(\delta) = \theta\}$ is good for 
$(\kappa,\theta,\chi)$, i.e. is not in $\check J_\theta[\lambda]$.

\sn
2) Any $S\subseteq\kappa$ good for $(\kappa,\theta,\chi)$ is the
union of $\kappa$ pairwise disjoint such sets.

\sn
3) In \ref{6.5} it suffices to assume that $S$ is good for
$(\lambda,\theta,\chi)$.

\sn
4) $\check J_\theta[\lambda]$ is a normal ideal on $\lambda$ and
there is no stationary $S \subseteq \{\delta < \lambda:\cf(\delta) 
= \theta\}$ which belongs to $\check J_\theta[\lambda] \cap 
\check I[\lambda]$.

\sn
5) In Definition \ref{6.6}, any $\mu>\lambda^{<\chi}$ is OK;
we can pre-assign $x\in \clH_{< \chi}(\mu)$ and demand $x\in M_i$.

\sn
6)  In \ref{6.5} we can replace the assumption 
``$S \subseteq \{\delta < \lambda : \cf(\delta) = \theta\}$ 
is stationary and in $\check I[\lambda]$" by 
``$S \subseteq \{\delta < \lambda : \cf(\delta) = \theta\}$ 
is stationary not in $\check J_\theta[\lambda]$" (which holds for
$S=\{\delta<\kappa:\cf(\delta)=\theta\}$).
\end{claim}

\begin{PROOF}{\ref{6.7}}
1)  Straightforward (play the game).

\sn
2) Similar to the proof of \ref{6.1}.

\sn
3)  Obvious.

\sn
4) Easy.

\sn
5) Easy.

\sn
6) Follows.  
\end{PROOF}

\begin{claim}\label{6.7A}
Assume that $\bar{\kappa},\theta,\chi$ are as in \ref{6.6} with
$|\partial |\leq\chi$. 

\sn
1) \underline{Then} the sequence $\LL\{\delta<\kappa_i:\cf(\delta)=\theta\}:
i<\partial\RR$ is good for $(\bar{\kappa},\theta,\chi)$.

\sn
2) If $\partial_1<\partial$ and $\LL S_i:i<\partial_1\RR$ is good
for $(\bar{\kappa} \rest \partial_1,\theta,\chi)$ \underline{then}
\[
\big\LL S_i : i < \partial_1 \big\RR \caret \big\LL\{ \delta < \kappa_i : \cf(\delta) = \theta\} : \partial_1 \le i < \partial\big\RR
\]
is good for $(\bar{\kappa},\theta,\chi)$.

\sn
3) If $\LL S_i:i<\partial_1\RR$ is good for $(\bar{\kappa},\theta,\chi)$ and $i(*)<\partial$, \underline{then} we can partition $S_{i(*)}$ to pairwise disjoint sets $\LL S_{i(*),\epsilon} : \epsilon<\kappa_i\RR$ such that 
for each $\epsilon<\kappa_i$, the sequence
\[
\big\LL S_i:i<i(*) \big\RR \caret \big\LL S_{i(*),\epsilon}\big\RR
\caret
\big\LL\{\delta<\kappa_i:\cf(\delta)=\theta\} : i(*) < i < \partial \big\RR
\]
is good for $(\bar{\kappa},\theta,\chi)$.

\sn
4) $\bar{S}$ good for $(\bar{\kappa},\theta,\chi)$ implies that $S_i$
is a stationary subset of $\kappa_i$ for each $i<\lh(\bar{\kappa})$.
\end{claim}

\begin{PROOF}{\ref{6.7A}}
Like \ref{6.7}. [In \ref{6.7A}(3) we choose, for $\delta\in S_{i(*)}$, a
club $C_\delta$ of $\delta$ of order type $\cf(\delta)$; for $j<\theta$,
$\epsilon<\kappa_{i(\alpha)}$, let $S^j_{i(*),\epsilon} =
\{\delta\in S_{i(*)}:\epsilon$ is the $j^\tthh$ member of $C_\delta\}$;
for some $j$ and unbounded $A \subseteq\kappa_{i(*)}$, $\LL
S^j_{i(*),\epsilon }:\epsilon \in A\RR$ are as required.]
\end{PROOF}

\noindent
Now we remove from \ref{6.5} (and subsequently \ref{6.8}) the hypothesis
``$\lambda$ is regular" when $\cf(\lambda) \ge \chi(*)$.

\begin{lemma}\label{6.9}
Suppose $\lambda^\theta = \lambda^{<\chi(*)}$, $\lambda$ is singular,
$\theta$ and $\chi(*)$ are regular, $\theta < \chi(*)$ and
$\cf(\lambda) \ge \chi(*)$.  Suppose further that 
$\lambda = \sum\limits_{i<\cf(\lambda)} \mu_i$ and each $\mu_i$ is 
regular $>\chi(*)+\theta^+$. \underline{Then}
we can find $\bfW = \{(\olsi M^\alpha\!,\eta^\alpha) : \alpha < \alpha(*)\}$
and functions $\dot\zeta : \alpha(*) \to \cf(\lambda)$, 
$\dot\xi : \alpha(*) \to \lambda$, and $h : \alpha(*) \to \lambda$, 
and $\{\mu'_i:i<\cf(\lambda)\}$ such that 
($\{\mu'_i:i<\cf(\lambda)\}= \{\mu_i:i<\cf(\lambda)\}$ and):
\mn
\begin{enumerate}
    \item[$(a0)$]  $h(\alpha)$ depends only on 
    $\LL \dot\zeta(\alpha),\dot\xi(\alpha)\RR$, $\alpha < \beta\ \Rightarrow\ \dot\zeta(\alpha)\le \dot\zeta(\beta)$,
    \[
        \alpha < \beta\wedge \dot\zeta(\alpha) = \dot\zeta(\beta)\ \Rightarrow\ 
        \dot\xi(\alpha) \le \dot\xi(\beta),
    \]
    and $\dot\xi(\alpha) < \mu'_{\dot\zeta(\alpha)}$.
\sn
    \item[$(a1)$]  As in \ref{6.5}, except that: 
    $\LL\eta^\alpha(3i) : i < \theta\RR$ is strictly increasing with limit 
    $\dot\zeta(\alpha)$ and $\LL \eta^\alpha(3i+1) : i < \theta\RR$ is 
    strictly increasing with limit $\dot\xi(\alpha)$ for $i<\theta$,
    \[
        \sup\big(|M^\alpha_i| \cap \mu'_{\zeta(\alpha)}\big) < \dot\xi(\alpha) = 
        \sup\big(|M^\alpha_\theta| \cap \mu'_{\dot\zeta(\alpha)}\big),
    \]
    and for every $i<\theta$,
    \[
        \sup\big(|M^\alpha_i|\cap\cf(\lambda)\big) < \dot\zeta(\alpha) = \sup\big(|M^\alpha_\theta| \cap \cf(\lambda)\big).
    \]
\sn
    \item[$(a2)$]   As in \ref{6.5}.
\sn
    \item[$(b0),(b1),(b2)$]  As in \ref{6.5}, but in clause (b3) we demand 
    $i = 2\mod 3$.
\end{enumerate}
\end{lemma}

\begin{remark}\label{6.9A} 
To make it similar to \ref{6.5}, we can fix $S^a$,
$S^a_i$, $S^b_i$, $S^b_{i,a}$, $\mu'_i$ as in the first paragraph of the
proof below.
\end{remark}

\begin{PROOF}{\ref{6.9}}
First, by \ref{6.7} [(1)+(2)], we can find pairwise disjoint
$S^a_i \subseteq \cf(\lambda)$ for $i < \cf(\lambda)$, each good for
$(\cf(\lambda),\theta,\chi(*))$ (and $\alpha\in S^a_i \Rightarrow
\alpha > i \wedge \cf(\alpha) = \theta$), and let $S^a=\bigcup\limits_{i<
\cf(\lambda)} S^a_i$. We define $\mu'_i\in\{\mu_j:j<i\}$ such that
$$\big(\forall i < \cf(\lambda)\big)[j\in S^a_i \Rightarrow \mu'_j = \mu_i].$$

Then for each $i$, by \ref{6.7A}(2),(3) (with $1,2,S_0,\kappa_0,
\kappa_1$ standing for $\sigma_1$, $\sigma$, $S^a_i$, $\cf(\lambda)$,
$\mu'_i)$, we can find pairwise disjoint subsets $\LL S^b_{i,\alpha}:
\alpha<\mu'_i\RR$ of $\{\delta<\mu'_i:\cf(\delta)=\theta\}$ such that
for each $\alpha<\mu'_\alpha$, $(S^a_i,S^b_{i,\alpha})$ is good for
$(\LL \cf(\lambda), \mu'_i\RR,\theta,\chi)$. Let $S^b_i=\bigcup
\{S^b_{i,\alpha }: \alpha<\mu'_i\}$.

Let $\cd$ be as in \ref{6.5}'s proof coding only for ordinals $i = 2 \mod 3$,
 and for $\zeta\in S^a_i$, $\xi\in S^a_{i,j}$ let
\begin{align*}
\bfW^0_{\zeta,\xi} = \Big\{ \big(\olsi M,\eta\big) : &\ \olsi M \text{ satisfies (a1), }
\zeta=\sup\{\eta(3i):i<\theta\}, \\
  &\ \xi=\sup\{\eta(3i+1):i<\theta\} \text{ and} \\
  &\ \text{for each } i<\theta, \text{ for some } y \in \clH_{<\chi(*)}(\lambda),\\
  &\ \eta(3i+2) = \cd\big(\LL M_j:j \le 3i+1\RR,\eta \rest (3i+1),y\big)\Big\}.
\end{align*}
\mn
The rest is as in \ref{6.5}'s proof.
\end{PROOF}

\medskip

\centerline {$* \qquad * \qquad *$}

\bn
The following Lemma improves \ref{6.5} when $\lambda$ satisfies a stronger
requirement, making the distinct $(\olsi M^\alpha\!,\eta^\alpha)$ interact
less. Lemmas \ref{6.8} + \ref{6.9} were used in the proof of \ref{5.2} (and
\ref{5.1}).

\begin{lemma}\label{6.8}
1)In \ref{6.5}, if $\lambda=\lambda^{\chi(*)}$ and $\chi(*)^\theta=
\chi(*)$ \underline{then} we can strengthen clause (b1) to
\sn
\begin{enumerate}
    \item[$(b1)^+$]  If $\alpha \ne \beta$ and 
    $\{\eta^\alpha \rest i : i < \theta\} \subseteq M^\beta$ then 
    $\alpha < \beta$ and $x\in M^\alpha_\theta \Rightarrow x \in M^\beta_\theta$.
\end{enumerate}

\sn
2) To clause \ref{6.5}(b1), we can add 
\begin{itemize}
    \item Moreover, if $\alpha < \chi(*) \Rightarrow |\alpha|^{\aleph_0}< \chi(*)$ then $\alpha < \beta + \big({<}\chi(*)\big)^\theta$.
\end{itemize}
\end{lemma}

\begin{PROOF}{\ref{6.8}}
1) Apply \ref{6.5} (actually, its proof) but using $\lambda$, $\chi(*)^+$,
$\theta$, instead of $\lambda,\chi(*),\theta$; and get $\bfW =
\{(\olsi M^\alpha\!,\eta^\alpha) : \alpha < \alpha(*)\}$, and the 
functions $\dot\zeta,h$.

Let $\cd$ be as in the proof of \ref{6.5}. Let $<^*$ be some well
ordering of $\clH_{<\chi(*)}(\lambda)$, and let $\cU$ be the
set of ordinals $\alpha<\alpha(*)$ such that for $i<\theta$,
$M^\alpha_i$ has the form $(N^\alpha_i,\in^\alpha_i,<^\alpha)$ and
$(|N^\alpha_i|,\in^\alpha_i,<^\alpha) \prec (\clH_{<\chi(*)}
(\lambda),\in,<^*)$.

Let $\alpha\in \cU$, by induction on $\epsilon<\chi(*)$ we define
$M^{\epsilon,\alpha}_i,\eta^{\epsilon,\alpha}$ as follows:
\mn
\begin{enumerate}
    \item[$(A)$]  $\eta^{\epsilon,\alpha}(i)$ is 
    $\cd(\LL\eta^\alpha(i),\epsilon \RR)$, (which is an ordinal $<\lambda$ but $>\eta^\alpha(i)$ and $>\epsilon)$
\sn
    \item[$(B)$]  $M^{\epsilon,\alpha}_i\prec N^\alpha_i$ is the Skolem Hull of 
    $\{\eta^{\epsilon,\alpha} \rest (j+1) : j < i\}$ inside $N^\alpha_i$, 
    using as Skolem functions the choice of the $<^*$-first element and making $M^{\epsilon,\alpha}_i\cap \chi(*)$ an ordinal. [If we want, we can use
    $\eta^{\epsilon,\alpha}$ such that it fits the definition in the proof of \ref{6.5}].
\end{enumerate}
\mn
Note that $\chi(*) = \chi^+\ \Rightarrow\ \chi+1 \subseteq M^\alpha_i$ and
$M^{\epsilon,\alpha}_i$ is definable in $M^{\epsilon,\alpha}_{i+1}$ as
$M^{\epsilon,\alpha }_i \in M^{\epsilon,\alpha}_{i+1}$ (by the definition of
$\bfW^0_\xi$ in the proof of \ref{6.5}). Similarly, $\LL M^{\epsilon,
\alpha}_j:j \le i\RR$ is definable in $M^\alpha_{i+1}$. It is easy to
check that the pair $(\olsi M^{\epsilon,\alpha},\eta^{\epsilon,\alpha})$
satisfies condition (a1) of \ref{6.5}.

Next we choose $\epsilon(\alpha) < \chi(*)$ by induction on 
$\alpha \in \cU$ as follows:
\mn
\begin{enumerate}
    \item[$(C)$]  $\epsilon(\alpha)$ is the first $\epsilon < \chi(*)$ such that if $\beta < \alpha$ but $\beta +\chi(*) > \alpha$ then
    \begin{enumerate}
        \item[$(*)$]  $\{\eta^{\alpha,\epsilon} \rest j : j < \theta\} 
        \nsubseteq M^{\beta,\epsilon(\beta)}_\theta$.
    \end{enumerate}
\end{enumerate}
\mn
This is possible and easy, as for $(*)$ it suffices to have for each
suitable $\beta$, $\epsilon\notin M^{\beta,\epsilon(\beta)}_\theta$, so
each $\beta$ ``disqualifies" $<\chi(*)$ ordinals as candidates for
$\epsilon(\alpha)$, and there are $<\chi(*)$ such $\beta$-s, and
$\chi(*)$ is by the assumptions (see \ref{6.5}) regular.

Now
\[
\bfW' = \big\{(\olsi N^{\alpha,\epsilon(\alpha)},\eta^{\alpha,\epsilon(\alpha)}) :
\alpha \in \cU \big\},
\]

\mn
$\dot\zeta \rest \cU,h \rest \cU$ are as required except that we should replace $\cU$ by an ordinal (and
adjust $\zeta,h$ accordingly). In the end replace $N^\alpha_i$ by
$N^\alpha_i\cap \clH_{<\chi(*)}(\lambda)$.

\sn
2) We have to prove the version of (b1) with the ``Moreover.''

Let $\clS \subseteq [\clH_{<\chi(*)}(\lambda)]^{\aleph_0}$ be MAD (that is, 
$u \neq v \in \clS \Rightarrow |u \cap v| < \aleph_0$ and $\clS$ 
is maximal under $\subseteq$) such that $\clS \cap [\clH_{<\chi(*)}(\zeta)]^{\aleph_0}$ is MAD for every $\zeta < \lambda$, {and demand} $|M_\theta^\alpha| \cap \clS \subseteq \big[|M_\theta^\alpha|\big]^{\aleph_0}$ is MAD. So it is well-known that the order $(\bfW,\leq_\bfW)$ is well founded.\footnote{So we use $\zeta < \chi(*) \Rightarrow |\zeta|^{\aleph_0} < \chi(*)$ to ensure that we can demand that $M_\theta^\alpha$ is as required. However, $\lambda \not\to (\omega)_2^{<\omega}$ will suffice.}
\end{PROOF}

\begin{claim}\label{6.10}
If in \ref{6.9} we add ``$\lambda=\lambda^{\chi(*)^\theta}$" (or the
condition from \ref{6.8}) then we can replace (b1) by
\mn
\begin{enumerate}
    \item[$(b1)^+$]  If $\{\eta^\alpha \rest i : i < \theta\} 
    \subseteq M^\beta_\theta$ \underline{then} $\alpha \le \beta$.
\end{enumerate}
\end{claim}

\begin{PROOF}{\ref{6.10}}
The same as the proof of \ref{6.8} combined with the proof of \ref{6.9}.
\end{PROOF}

\bigskip

\subsection{Black Boxes: For $\theta$ countable}

Next we turn to the case (of black boxes with) $\theta=\aleph_0$. We shall
deal with several cases.

\begin{lemma}\label{6.11}
Suppose that
\mn
\begin{enumerate} 
    \item[$(*)$]   $\lambda$ is a regular cardinal, $\theta = \aleph_0$, 
    $\mu = \mu^{<\chi(*)} < \lambda \le 2^\mu$,\\ 
    $S \subseteq \{\delta < \lambda : \cf(\delta) = \aleph_0\}$ is stationary, 
    and $\aleph_0 < \chi(*) = \cf(\chi(*))$.
\end{enumerate}
\mn
\underline{Then} we can find
\[
\bfW = \{(\olsi M^\alpha ,\eta^\alpha):\alpha<\alpha(*)\}
\]
and functions
\[
\dot\zeta : \alpha(*) \to S \text{ and } h:\alpha(*) \to \lambda
\]
such that:
\mn
\begin{enumerate}
    \item[$(a0)$-$(a2)$] As in \ref{6.5}.
\sn
    \item[$(b0)$-$(b2)$] As in \ref{6.5}, and even
    \begin{enumerate}
        \item[$(b1)^*$]  $\alpha \ne \beta$, $\{\eta^\alpha \rest n : n < \omega\} \subseteq M^\beta_\omega$ implies $\alpha < \beta$ and even $\dot\zeta(\alpha) < \dot\zeta(\beta)$.
    \end{enumerate}
\sn
    \item[$(c1)$]  If $\dot\zeta(\alpha) = \dot\zeta(\beta)$ \underline{then} $|M^\alpha_\omega| \cap \mu = |M^\beta_\omega| \cap \mu$, there is an 
    isomorphism $h_{\alpha,\beta}$ from $M^\alpha_\omega$ onto $M^\beta_\omega$,
    mapping $\eta^\alpha(n)$ to $\eta^\beta(n)$ and $M^\alpha_n$ to $M^\beta_n$ 
    for $n<\omega$, and $h_{\alpha,\beta} \rest (|M^\alpha_\omega| \cap |M^\beta_\omega|)$ is the identity.
\sn
    \item[$(c2)$]  There is $\olsi C = \LL C_\delta : \delta \in S\RR$,
    $C_\delta$ an $\omega$-sequence converging to $\delta$, $0 \notin C_\delta$,
    and letting $\LL\gamma^\delta_n : n < \omega\RR$ enumerate 
    $\{0\} \cup C_\delta$ we have, when $\dot\zeta(\alpha) = \delta$:
    \begin{enumerate}
        \item[$(i)$]  $\lambda\cap |M^\alpha_n|\subseteq\gamma^\delta_{n+1}$ but $\lambda \cap |M^\alpha_n|$ is not a subset of $\gamma^\delta_n$, (hence $M^\alpha_n \cap [\gamma^\delta_n,\gamma^\delta_{n+1}) \ne \varnothing$).
\sn
        \item[$(ii)$]  $C_\delta\cap |M^\alpha_\omega|=\varnothing$
\sn
        \item[$(iii)$] If in addition $\dot\zeta(\beta)=\delta$ then for each $n$, $h_{\alpha,\beta}$ maps $|M^\alpha_\omega| \cap [\gamma^\delta_n,\gamma^\delta_{n+1})$ onto $|M^\beta_\omega| \cap [\gamma^\delta_n,\gamma^\delta_{n+1}]$.
\sn
        \item[$(iv)$]  If $\dot\zeta(\beta) = \delta = \dot\zeta(\alpha)$ and $\lambda = \lambda^{<\chi(*)}$ \underline{then} $|M^\alpha_\omega| \cap \gamma^\delta_1 = |M^\beta_\omega| \cap \gamma^\delta_1$.
    \end{enumerate}
\end{enumerate}
\end{lemma}

\begin{remark}\label{6.11u}
1) We only use $\lambda \leq 2^\mu$ in order to get
``$h_{\alpha,\beta} \rest (|M^\alpha_\omega| \cap |M^\beta_\omega|) = \id$"
in condition {(c1)}.

\sn
2)  Below we quote ``guessing of clubs" --- that is clause (ii) in the 
proof; without this we just get a somewhat weaker conclusion.
\end{remark}

\begin{PROOF}{\ref{6.11}}
Let $S$ be the disjoint union of stationary $S_{\alpha,\beta,\gamma}$, for $\alpha < \mu$, $\beta < \lambda$, $\gamma < \lambda$.

For each $\alpha$, $\beta$, $\gamma$, let 
$\LL C_\delta : \delta \in S_{\alpha,\beta,\gamma}\RR$ satisfy:
\mn
\begin{enumerate}
    \item[$\boxtimes$] 
    \begin{enumerate}[(i)]
        \item $C_\delta$ is an unbounded subset of $\delta$ of order type $\omega$.

        \item For every club $C$ of $\lambda$, for stationarily many 
        $\delta \in S_{\alpha,\beta,\gamma}$, we have $C_\delta \subseteq C$.

        \item $0 \notin C_\delta$
    \end{enumerate}
\end{enumerate}
\mn
(exists by \cite[2.2]{Sh:331} or \cite{Sh:365}=\cite[Ch.III]{Sh:g}).  

Let $\bfW^*$ be the family of quadruples $(\delta,\olsi M,\eta,C)$
 such that:
\mn
\begin{enumerate}
    \item[$\circledast$]    
    \begin{enumerate}
        \item[$(\alpha)$]  $(\olsi M,\eta)$ satisfies the requirement (a1) 
        (so $\olsi M = \LL M_n : n < \omega\RR$).
\sn
        \item[$(\beta)$]  $0 \notin C$, and letting $\{\gamma_n : n < \omega\}$ enumerate $C \cup \{0\}$ in increasing order, we have $\lambda \cap M_n$ is a subset of $\gamma_{n+1}$ but not of $\gamma_n$, and $\bigcup\limits_{n<\omega}\gamma_n = \delta$ and $C \cap (\bigcup\limits_n M_n) = \varnothing$.
\sn
        \item[$(\gamma)$]  $\bigcup\limits_n |M_n|\subseteq \clH_{<\chi(*)}(\mu+\mu)$
\sn
        \item[$(\delta)$]  In $\tau(M_n)$ there is a two-place relation $R$ 
        and a one-place function $\cd$. (We do not necessarily require 
        $\cd\rest M_n = \cd^{M_n}$; similarly for $R$ --- see below. 
        Recall that as usual, $\tau(M_n) \in \clH_{<\chi(*)}(\chi(*))$ for transparency.)
    \end{enumerate} 
\end{enumerate}
\mn
As $\mu^{<\chi(*)} = \mu$, clearly $|\bfW^*| = \mu$, so let
\[
\bfW^* = \big\{(\delta^\alpha,\LL M_{\alpha,n} : n < \omega\RR,\eta_\alpha,C^\alpha) : \alpha < \mu \big\}.
\]

\mn
If $\lambda = \lambda^{<\chi(*)}$ let $\{N_\beta : \beta < \lambda\}$ 
list the models $N \in \clH_{<\chi(*)}(\lambda)$ with 
$\tau(N) \in \clH_{<\chi(*)}(\chi(*))$.  

Also, let $\LL A_\alpha : \alpha < \lambda\RR$ be a sequence of pairwise
distinct subsets of $\mu$, and define the two place relation $R$ on $\lambda$ by
\[
\gamma_1\; R\; \gamma_2 \Leftrightarrow  [\gamma_1 < \mu\ \wedge\ \gamma_1\in
A_{\gamma_2}].
\]

\mn
Lastly, for $\delta\in S_{\alpha,\beta,\gamma}$ let $\bfW^0_\delta$
be the set of pairs $(\olsi M,\eta)$ such that:
\begin{enumerate}
    \item[$\oplus$]
    \begin{enumerate}
        \item $\olsi M = \LL M_n : n < \omega \RR$, $\eta\in {}^\omega\!\lambda$

        \item $(\olsi M,\eta)$ satisfies \ref{6.5}(a1).  In particular:
        \begin{enumerate}
            \item[$(\alpha)$] $\eta$ is increasing with limit $\delta$.

            \item[$(\beta)$] there is an isomorphism $h$ from $\bigcup\limits_{n<\omega} M_n$ onto 
            $\bigcup\limits_{n<\omega} M_{\alpha,n}$.
            
            \item [$(\gamma)$] $h$ maps $\eta(n)$ to $\eta^\alpha(n)$ and $M_n$ onto $M_{\alpha,n}$.
 
            \item[$(\delta)$] $h$ preserves $\in$, $R$, $\cd(x) = y$ and their negations. (For $R$ and $\cd$: in $\bigcup\limits_{n<\omega} M_n$ 
            we mean the standard $\cd$ restricted to 
            $\bigcup\limits_{n<\omega} M_{\alpha,n}$ 
            {as in clause $\circledast(\delta$) above}.)
        \end{enumerate}

        \item $(\forall\epsilon<\lambda)\big[\epsilon\in \textstyle\bigcup\limits_n M_n\ \Rightarrow\ \otp(C_\delta \cap \epsilon) = \otp\big(C^\alpha\cap h(\epsilon) \big)\big].$

        \item If $\lambda = \lambda^{<\chi(*)}$ then 
        $N_\beta = \big(\bigcup\limits_n M_n\big) \rest \big\{x \in \bigcup\limits_n M_n : \cd(x) < \min(C_\delta) \big\}$.        
    \end{enumerate}
\end{enumerate}

We proceed as in the proof of \ref{6.5} after $\bfW^0_\delta$ was 
defined (only $\dot\zeta(\alpha)=\delta\in
S_{\alpha_1,\beta_1,\gamma_1} \Rightarrow h(\alpha)=\gamma_1)$.

Suppose \textsf{G} is a winning strategy for Player I. So suppose that if Player II
has chosen $\eta(0),\eta(1),\ldots,\eta(n-1)$, Player I will choose
$M_\eta$. So $|M_\eta|$ is a subset of $\clH_{<\chi(*)}(\lambda)$ of
cardinality $<\chi(*)$ and $\Rang(\eta) \subseteq M_\eta$. For $\eta \in
{}^\omega\!\lambda$ we define $M_\eta=\bigcup\limits_{\ell<\omega}
M_{\eta \rest \ell}$.

Let $\clT_n$ be the set of $\eta \in {}^n\!\lambda$ such that $M_\eta$
is well defined, so $\bigcup \{\clT_n : n < \omega\}$
is a subtree of $({}^{\omega >}\!\lambda,\lhd)$
with each node having $\lambda$ immediate successors.

We can find a function $\bfc_n$ from $\clT_n$ into $\mu$ 
such that $\bfc_n(\eta) = \bfc_n(\nu)$ iff there is an
isomorphism $h$ from $M_\eta$ onto $M_\nu$ mapping $M_{\eta \rest k}$ onto
$M_{\nu \rest k}$ for every $k < n$.
By \cite[1.10=L1.7]{Sh:E62}, or \cite{Sh:117}, or the proof of \ref{6.11A} below,
there is $\clT$ such that:
\begin{enumerate}
    \item[$(*)$]
    \begin{enumerate}
        \item $\clT \subseteq {}^{\omega>}\!\lambda$
        
        \item $\clT$ is closed under initial segments.

    \item $\LL\ \RR\in \clT$ 
    
    \item $\eta\in \clT \Rightarrow (\exists^\lambda\alpha) \big[\eta \caret \LL\alpha\RR \in \clT \big]$

    \item $\bfc_n \rest (\clT \cap \clT_n)$ is constant. 
    \end{enumerate} 
\end{enumerate}

It follows that for any $\nu_* \in \lim(\clT)$ we can find 
$\LL h_\eta : \eta \in \clT \RR$ such that $h_\eta$ is an isomorphism
from $M_{\nu_* \rest \lh(\eta)}$ onto $M_\eta$ increasing with $\eta$.

Note that above, all those isomorphisms are unique
as the interpretation of $\in$ satisfies comprehension.
Also, clause (c1) follows from the use of $R$.

\sn 
The rest should be clear.
\end{PROOF}

\begin{lemma}\label{6.11A}
Let $S$, $\lambda$, $\mu$, $\theta$, $\chi(*)$ be as in \ref{6.11}$(*)$, and in addition:
\[
\aleph_0 \le \kappa = \cf(\kappa) < \chi(*) = \cf(\chi(*)),\quad 
\big(\forall\chi < \chi(*)\big)\big[\chi^{<\kappa} < \chi(*)\big],\quad 
(\forall\alpha < \lambda)\big[|\alpha|^{<\kappa} < \lambda\big].
\]

\sn
\underline{Then} we can find $\bfW = \{(\olsi M^\alpha\!,\eta^\alpha):\alpha<
\alpha(*)\}$ and functions $\dot\zeta:\alpha(*) \to S$ and
$h:\alpha(*) \to \lambda$ such that:
\mn
\begin{enumerate}
    \item[$(a0),(b0),(b2)$]  As in \ref{6.11} (i.e. as in \ref{6.5}).
\sn
    \item[$(b1)^*,(c1),(c2)$]   As in \ref{6.11}.
\sn
    \item[$(a1)^*$]  As in  in \ref{6.5}$(a1)$, except that we omit 
    ``$\LL M_j : j \le i\RR \in M_{i+1}$" and add: 
    $\left[a \subseteq |M_i| \wedge |a| < \kappa\right] \Rightarrow a \in M_i$, 
    and for $i<j$, $M_i \cap \lambda$ is an initial segment of $M_j \cap \lambda$.
\sn
    \item[$(a2)^*$]  For every expansion $\cA$ of $(\clH_{<\chi(*)}(\lambda),\in,<)$ by $\chi < \chi(*)$ relations (with $\tau(\cA) \subseteq \clH_{<\chi(*)}(\chi(*))$), for some $\alpha < \alpha(*)$, for every $n$, 
    $M^\alpha_n \prec \cA$. In fact, for stationarily many $\zeta \in S$, 
    there is such $\alpha$ satisfying $\dot\zeta(\alpha) = \zeta$.
\end{enumerate}
\end{lemma}

\begin{remark}\label{6.11D}
We can retain $(a1)^*$ and add 
$a \subseteq M_i \wedge |a| < \kappa \Rightarrow a \in M_i$.
\end{remark}

\begin{PROOF}{\ref{6.11A}}
Similar to \ref{6.11}, using the proof of \cite{Sh:247},
but for completeness we give details.

We choose $\LL S_{ \alpha,\beta,\gamma}: \alpha < \mu,\ 
\beta < \lambda,\ \gamma < \lambda \RR$ as there.
The main point is that defining $\bfW^*$ we have one additional demand:
\mn
\begin{enumerate}
    \item[$(\eps)$]  If $n < \omega$ and $u \subseteq M_n$ has cardinality 
    $< \kappa$, then $u \in M_n$.
\end{enumerate}
\mn
We then define $\bfW^0_\delta$ and $\LL N_\alpha : \alpha < \lambda\RR$ as there.

This gives the changed demand in (a1)$^*$, but it creates extra work in verifying
the demand (a2)$^*$.

So let a model $\cA$ and cardinal $\chi = \chi^{< \kappa} < \chi(*)$ be 
given as there; as usual, $\tau(\cA) \in \clH_{ < \chi(*)}(\chi(*))$
and $\cA$ expands $(\clH_{< \chi(*)}(\lambda),\in,<)$.
For every 
$$\bfx = (\delta_\bfx,\olsi M_\bfx,\eta_\bfx,C_\bfx) \in \bfW^*$$
we define a family $\cF_\bfx$, a
function $n : \cF_\bfx \rightarrow \omega$ and a function $\rank_\bfx$
from $\cF_\bfx$ into $\Ord \cup \{\infty\}$ as follows:
\mn
\begin{enumerate}
    \item[$(\alpha)$]  $\cF_\bfx = \bigcup \{\cF_{\bfx,n}: n < \omega\}$
\sn
    \item[$(\beta)$]   $\cF_{\bfx,n} = \{f : f$ is an elementary embedding of $M_{\bfx,n}$ into $\cA\}$
\sn
    \item[$(\gamma)$]  $n(f) = k$ if and only if $f \in \cF _{\bfx,k}$.
\sn
    \item[$(\delta)$]  $\rank(f) = \bigcup \big\{ \epsilon + 1 : \text{ for every } \alpha < \lambda$ there is $g \in \cF_ {\bfx,n(f)}$ extending $f$ such that $\beta = \rank_\bfx(g)$ and $\Rang(g) \cap \alpha = \Rang(f) \cap \lambda\big\}$.
\end{enumerate}
\mn
Now

\mn
\textbf{Case 1}: For no $\bfx \in \bfW^*$ and
$f \in \cF_{\bfx,0}$ do we have $\rank_\bfx(f) =  \infty$.

For every $\bfx \in \bfW^*$ and $f \in \cF_\bfx$ let 
$\beta(f,\bfx)$ be the first ordinal $\alpha < \lambda $ such that
if $\rank_\bfx(f) = \epsilon$ \underline{then} there is no $g \in 
\cF_{\bfx,n(f)+1}$ extending $f$ with $\rank_\bfx(g) = \epsilon$
and $\Rang(g) \cap \alpha = \Rang(f) \cap \lambda$.

Next, let $\bar \cA = \LL \cA_i:i < \lambda \RR$ be an increasing
continuous sequence of elementary submodels of $\cA$,
each of cardinality $< \lambda$ such that $\LL \cA_j:j \le i \RR
\in \cA_{i+1}$.

Easily the set $E = \{i < \lambda : \cA_i \cap \lambda = i > \mu\}$
is a club of $\lambda$.

Choose, by induction on $n < \omega$, an ordinal $i_n$ increasing
with $n$ such that $i_n \in E$ is of cofinality $\kappa$ (this is possible as 
$\kappa = \cf(\kappa)<  \lambda$) hence $\cA_{i_n}$ is 
an elementary submodel of $\cA$ of cardinality $< \lambda$.

Choose $M \prec \cA$ of cardinality $\chi$, including 
$\{i_n : n < \omega\} \cup \{\bar \cA, \bfW^*\} \cup (\chi+1)$ such that every 
$u \subseteq M$ of cardinality $< \kappa$ belongs to $M$.

Note that if $u \subseteq \cA_{i_n}$ has cardinality $< \kappa$ then
$u \in \cA_{i_n}$ because $i_n \in E$ and $\cf(i_n) = \kappa$, hence this holds for every $\cA_{i_n} \cap M$.

Let $M^*_n$ be $\cA \rest (\cA_{i_n} \cap M)$; easily
$M^*_n \in \cA_{i_n}$, so $[u \subseteq M^*_n \wedge |u| < \kappa] 
\Rightarrow u \in M^*_n$.  We can find $\bfx \in \bfW$,
and isomorphism $f_n$ from $M_{\bfx,n}$ onto $M^*_n$ increasing
with $n$.  Now clearly $\bfx \in \cA_{i_n}$. 

\sn
[Why? As $\mu = \mu^{<\chi(*)}$ and $\mu + 1 \subseteq \cA_{i_n}$.
Also, $f_n \in \cF_{\bfx,n}$
and {these $f_n$ are unique} 
as those models
expand a submodel of $(\clH_{< \chi(*)}(\lambda),\in,<)$ and are
necessarily transitive over the ordinals.]

Similarly by the choice of $\bfx$, we have $f_n \subseteq f_{n+1}$.  So
$\LL \rank_ \bfx(f_n):n < \omega \RR$ is constantly $\infty$ as
otherwise we get an infinite decreasing sequence of ordinals.

But this contradicts our case assumption.
\medskip

\noindent
\textbf{Case 2}:  Not case 1.

So we choose $\bfx \in \bfW^*$ and $f \in \cF_{\bfx,0}$ such that
$\rank_\bfx(f) = \infty$.

We easily get the desired contradiction and even a $\Delta$-system
tree of models.  How?  Let $\LL \eta_\alpha:\alpha < \lambda \RR$
list ${}^{\omega >}\!\lambda $ such that $\eta_\alpha \lhd \eta_\beta$
implies $\alpha < \beta$.

Now we choose a pair $(f_{\eta_\alpha},\gamma_\alpha)$
by induction on $\alpha < \lambda$ such that
\mn
\begin{enumerate}
    \item[$(i)$]  $f_{\eta_\alpha} \in \cF_{\bfx,\lh(\eta_\alpha)}$
\sn
    \item[$(ii)$]  $\gamma_\alpha = \sup\!\big( \bigcup\{\lambda \cap \Rang(f_{\eta_\beta}):\beta < \alpha\}\big)$
\sn
    \item[$(iii)$]   if $\eta_\beta \lhd \eta_\alpha$ and $\lh(\eta_\alpha) = (\lh(\eta_\beta)+1$ then $\gamma_\alpha \cap \Rang(f_{\eta_\alpha}) = \lambda \cap \Rang(f_{\eta_\beta})$.
\end{enumerate}
\mn
There is no problem to carry the induction. This finishes the proof.
\end{PROOF}

\begin{lemma}\label{6.12}
1)  In \ref{6.11A}, if in addition $\lambda=\mu^+$ \underline{then} we can add:
\mn
\begin{enumerate}
    \item[$(c3)$]  If $\dot\zeta(\alpha) = \dot\zeta(\beta)$, 
    \underline{then} $|M^\alpha_\omega|\cap |M^\beta_\omega| \cap\lambda$ 
    is an initial segment of $|M^\alpha_\omega| \cap \lambda$ and of
    $|M^\beta_\omega|\cap\lambda$, so when $\alpha \ne \beta$ it is a 
    bounded subset of $\dot\zeta(\alpha)$.
\end{enumerate}
\mn
2) In \ref{6.11A} (and \ref{6.12}), when $\kappa > {\aleph_0}$
\underline{then} it follows that:
\mn
\begin{enumerate}
    \item[$(c4)^*$]  If $\alpha \ne \beta$ and 
    $\{\eta^\alpha \rest n : n < \omega\} \subseteq M^\beta_\omega$ 
    then $\olsi M^\alpha,\bar{\eta}^\alpha\in M^\beta_\omega$.
\end{enumerate}
\mn
3)  Assume $\lambda = \mu^+$, $\mu = \mu^\kappa$, 
$S \subseteq \{\delta < \lambda : \cf(\delta) = {\aleph_0}\}$ 
is a stationary subset of $\lambda$, and $\LL C_\delta : \delta \in S \RR$
guesses clubs (and $C_\delta$ is an unbounded subset of $\delta$
of order type $\omega$, of course). 

\underline{Then} we can find $\LL \olsi N_\eta:\eta \in \Gamma \RR$
such that:
\begin{enumerate}
    \item[$(a)$]  $\Gamma = \bigcup \{\Gamma_\delta : \delta \in S\}$, where $\Gamma_\delta \subseteq \{\eta : \eta$ an increasing $\omega$-sequence of ordinals $< \delta$ with limit $\delta\}$ and $\delta(\eta) = \delta$ when $\eta \in \Gamma_\delta$ and $\delta \in S$.
\sn
    \item[$(b)$]  $\olsi{N}_\eta$ is $\LL N_{\eta,n}: n \le \omega \RR$, 
    which {is} $\prec$-increasing continuous, and we let $N_\eta = N_{\eta,\omega}$.
\sn
    \item[$(c)$]  Each $N_\eta$ is a model of cardinality $\kappa$ 
    (with vocabulary $\subseteq \clH(\kappa^+)$ for notational simplicity),
    universe $\subseteq \delta \defeq \delta(\eta)$, 
    $N_{\eta,n} = N_\eta \rest \gamma^\delta_n$ (where $\gamma^\delta_n$ 
    is the $n^\tthh$ member of $C_\delta$), and $N_\eta \cap (\gamma^\delta_n,\gamma^\delta_{n+1}) \ne \varnothing$.
\sn
    \item[$(d)$]  For every distinct $\eta,\nu \in \Gamma_\delta$ 
    with $\delta \in S$, for some $n < \omega$, we have 
    $N_\eta \cap N_\nu = N_{\eta,n} = N_{\nu,n}$.
\sn
    \item[$(e)$]  For every $\eta,\nu \in \Gamma_\delta$ the models 
    $N_\eta,N_\nu $ are isomorphic; moreover, there is such an isomorphism 
    $f$ which preserves the order of the ordinals and maps $N_{\eta,n}$ 
    onto $N_{\nu,n}$.
\sn
    \item[$(f)$]  If $ \cA$ is a model with universe $\lambda$ and vocabulary $\subseteq \clH(\kappa^+)$ then for stationarily many $\delta \in S$, 
    for some $\eta \in \Gamma_\delta \subseteq \Gamma$, we have 
    $N_\eta \prec \cA$.  Moreover, if $\kappa^\partial = \kappa$ and $h$ 
    is a one to one function from ${}^{\partial} \lambda$ into $\lambda$ \underline{then} we can add: if $\rho \in{}^{\partial}(N_{\eta,n})$ 
    then $h(\rho) \in N_{\eta,n}$.
\end{enumerate}
\end{lemma}

\begin{PROOF}{\ref{6.12}}
1) Let $g^0,g^1$ be two place functions from $\lambda\times\lambda$ 
to $\lambda$ such that for $\alpha \in[\mu,\lambda]$, 
$\LL g^0(\alpha,i) : i < \mu\RR$ enumerates $\{j:j<\mu\}$ without repetition 
and $g^1(\alpha,g^0(\alpha,i)) = i$ for $i < \lambda$.

Now we can restrict ourselves to $\olsi M^\alpha$ such that each
$M^\alpha_i$ (for $i \le \omega$) is closed under $g^0,g^1$. Then {(c3)}
follows immediately from
\[
\dot\zeta(\alpha) = \dot\zeta(\beta) \Rightarrow 
|M^\alpha_\omega| \cap \mu = |M^\beta_\omega| \cap \mu
\]
(required in {(c1)}).

\sn
2) Should be clear.

\sn
3) This just rephrases what we have proved above.
\end{PROOF}

\begin{lemma}\label{6.13}
Suppose that $\lambda = \mu^+$, 
$\mu = \kappa^{\aleph_0} = 2^\kappa > 2^{\aleph_0}$, $\cf(\kappa) = \aleph_0$ 
and $S \subseteq \{\delta < \lambda : \cf(\delta)=\aleph_0\}$ is stationary, 
$\theta = \aleph_0$, $\aleph_0 < \chi(*) = \cf(\chi(*)) < \kappa$. 
\underline{Then} we can find
$\bfW = \{(\olsi M^\alpha\!,\eta^\alpha): \alpha<\alpha(*)\}$ and functions
\[
\dot\zeta : \alpha(*) \to S,\quad h:\alpha(*) \to \lambda
\]

\mn
and $\LL C_\delta:\delta\in S\RR$ with $\LL
\gamma^\delta_n:n < \omega \RR$ listing $C_\delta$ in increasing order
such that:
\mn
\begin{enumerate}
    \item[$(a0)$-$(a1)$]  As in \ref{6.5}.
\sn
    \item[$(a2)^*$]  As in \ref{6.11A}.
\sn
    \item[$(b0)$-$(b2)$]  As in \ref{6.5}, and even
    \begin{enumerate}
        \item[$(b1)^*$]  $\alpha\ne \beta$, 
        $\{\eta^\alpha \rest n : n < \omega\} \subseteq M^\beta_\omega$ 
        implies $\alpha < \beta$ and even $\dot\zeta(\alpha) < \dot\zeta(\beta)$. 
    \end{enumerate}
\sn
    \item[$(c1)$-$(c3)$]  As in \ref{6.11} + \ref{6.12}(1).
\sn
    \item[$(c4)$]  If $\dot\zeta(\alpha) = \dot\zeta(\beta) = \delta$ but 
    $\alpha \ne \beta$ \underline{then} for some $n_0 \ge 1$, there are no 
    $n > n_0$ and $\alpha_1 \le \beta_2 \le \alpha_3$ satisfying:
    \[
    \begin{array}{l}
        \alpha_1\in |M^\alpha_\omega|\cap [\gamma^\delta_n,\gamma^\delta_{n+1}), \\
        \beta_2\in |M^\beta_\omega|\cap [\gamma^\delta_n,\gamma^\delta_{n+1}),\\
        \alpha_3\in |M^\alpha_\omega|\cap [\gamma^\delta_n,\gamma^\delta_{n+1}),
    \end{array}
    \]
    i.e., either 
    $$\sup\big([\gamma^\delta_n,\gamma^\delta_{n+1}) \cap |M^\alpha_\omega|\big) < \min\big([\gamma^\delta_n,\gamma^\delta_{n+1}) \cap |M^\beta_\omega|\big)$$
    or $\sup\big([\gamma^\delta_n,\gamma^\delta_{n+1}) \cap |M^\beta_\omega|\big) < \min\big([\gamma^\delta_n,\gamma^\delta_{n+1})\cap |M^\alpha_\omega|\big)$.
\sn
    \item[$(c5)$]  If $\Upsilon <\kappa$ and there is $B\subseteq {}^\omega\kappa$, $|B| = \kappa^{\aleph_0}$ which contains no perfect 
    set with density $\Upsilon$ (this holds trivially if $\kappa$ is strong 
    limit), \underline{then} also $\{\eta^\alpha:\alpha<\alpha(*)\}$ does 
    not contain such a set. (See \ref{6.13A}.)
\end{enumerate}
\end{lemma}

\begin{PROOF}{\ref{6.13}}
We repeat the proof of \ref{6.11} with some changes.

Let $\LL S_{\alpha,\beta,\gamma} : \alpha < \mu ,\ \beta < \lambda,\ \gamma <
\lambda\RR$ be pairwise disjoint stationary subsets of $S$. Let
$g^0,g^1$ be as in the proof of \ref{6.12}. By \ref{6.3A} there is a
sequence $\LL C_\delta: \delta\in S\RR$ such that:
\begin{enumerate}
    \item[$(i)$]  $C_\delta$ is a club of $\delta$ of order type $\kappa$ (not $\omega$!), $0\notin C_\delta$.
\sn
    \item[$(ii)$]  For $\alpha<\mu$, $\beta<\lambda$, $\gamma<\lambda$, for every club $C$ of $\lambda$, the set
    \[
        \{\delta \in S_{\alpha,\beta,\gamma} : C_\delta \subseteq C\}
    \]
    is stationary.
\end{enumerate}
\mn
We then define $\bfW^*,(\delta^j,\LL M_{j,n} : n < \omega \RR, 
\eta_j,C^j)$ for $j<\mu$, $A_\alpha$ for $\alpha < \lambda$, and $R$ as 
in the proof of \ref{6.11}.

Now, for $\delta\in S_{\alpha,\beta,\gamma}$ let $\bfW^1_\delta$ be the
collection of all systems $\LL M_\rho,\eta_\rho : \rho \in {}^{\omega>} \!\kappa\RR$ such that:
\begin{enumerate}
    \item[$(i)$]  $\eta_\rho$ is an increasing sequence of ordinals of length
$\lh(\rho)$.
\sn
    \item[$(ii)$]  $\otp\big(C_\delta\cap\eta_\rho(\ell)\big) = 1 + \rho(\ell)$ 
    for $\ell < \lh(\rho)$.
\sn
    \item[$(iii)$]  There are isomorphisms $\LL h_\rho : \rho \in {}^{\omega>} \!\kappa\RR$ such that $h_\rho$ maps $M_\rho$ onto $M_{\alpha,\lh(\rho)}$ 
    preserving $\in,R$, $\cd(x)=y$, $g^0(x_1,x_2)=y$, $g^1(x_1,x_2)=y$ 
    (and their negations).
\sn
    \item[$(iv)$]  If $\rho\lhd\nu$ then $h_\rho \subseteq h_\nu$, 
    $M_\rho\prec M_\sigma$, and $M_\rho\in M_\nu$.
\sn
    \item[$(v)$]  $M_\rho \cap C_\delta = \varnothing$, and 
    $M_\rho \cap \lambda \subseteq \bigcup\limits_\ell [\gamma_{\rho(\ell)},\gamma_{\rho(\ell)+1})$, where $\gamma_\zeta$ 
    is the $\zeta^\tthh$ member of $C_\delta$.
\sn
    \item[$(vi)$]  If $\rho\in {}^{\omega>}\!\kappa$, $\ell < \lh(\rho)$, 
    and $\gamma$ is the $(1+\rho(\ell))^\tthh$ member of $C_\delta$ then 
    $M_\ell \cap \gamma$ depends only on $\rho \rest \ell$ and 
    $M_\rho \rest \gamma \prec M_\rho$.
\sn
    \item[$(vii)$]  $N_\beta = M_{\LL\ \RR}$.
\end{enumerate}
\mn
Now clearly $|\bfW^1_\delta| \le \mu$, so let $\bfW^1_\delta = \big\{ \big\LL (M^j_\rho,\eta^j_\rho) : \rho \in {}^{\omega >}\!\kappa \big\RR : j < \mu \big\}$.  Let $\LL\rho_j : j < \mu\RR$ be a list of
distinct members of ${}^\omega \kappa $, for {(c5)} --- choose as there.

Let
\[
M^j_\ell = \bigcup\limits_{\ell<\omega} M^j_{\rho_j \rest \ell},\quad
\eta^j = \big\LL \eta^j_{\rho_j \rest (\ell +1)}(\ell+1) : \ell \le
\omega \big\RR.
\]
Now,
\[
\big\{\LL M^j_\ell : \ell < \omega\RR : j < \mu \big\}
\]

\mn
is as required in {(c4)}. Also, {(c5)} is straightforward, as taking
union for all $\delta$-s changes little. (Of course, we are omitting
$\delta$-s where we get unreasonable pairs.)

The rest is as before. 
\end{PROOF}

\begin{remark}\label{6.13A}
The existence of $B$ as in {(c5)} is proved for some $\Upsilon$, for all
strong limit $\kappa$ of cofinality $\aleph_0$. By
\cite[Ch.II,6.9,pg.104]{Sh:g}, much stronger conclusions hold. If 
$2^\kappa$ is regular and belongs to $\{\cf(\prod \kappa_n/D) : D \text{ an 
ultrafilter on } \omega,\ \kappa_n < \kappa\}$, or $2^\kappa$ is singular 
and is the supremum of this set, then it exists for 
$\Upsilon = (2^{\aleph_0})^+$. Now, if above we
replace $D$ by the filter of co-bounded subsets of $\omega$, then we get it
even for $\Upsilon = \aleph_0$; by \cite[Part D]{Sh:E12} the requirement holds,
e.g., for $\beth_\delta$ for a club of $\delta<\omega_1$.

Moreover, under this assumption on $\kappa$ we can demand
(essentially, this is expanded in \ref{6.14A})
We strengthen clause (c4) to:
\begin{enumerate}
    \item[$(c4)^*$] If 
    $\dot\zeta(\alpha) = \dot\zeta(\beta) = \delta$ but $\alpha \ne \beta$
    \underline{then} for some $n_0 \ge 1$, either for every $n \in [n_1,\omega)$
    we have $$\sup\big([\gamma^\delta_n,\gamma^\delta_{n+1}) \cap |M^\alpha_\omega|\big) < \min\big([\gamma^\delta_n,\gamma^\delta_{n+1}) \cap |M^\beta_\omega|\big)$$ 
    or for every $n \in [n_1,\omega)$  we have
    $$\sup\big([\gamma^\delta_n,\gamma^\delta_{n+1})\cap |M^\beta_\omega|\big) < \min\big([\gamma^\delta_n,\gamma^\delta_{n+1})\cap |M^\alpha_\omega|\big).$$
\end{enumerate}
\end{remark}

\begin{lemma}\label{6.13B}
We can combine \ref{6.13} with \ref{6.11A}.
\end{lemma}

\begin{PROOF}{\ref{6.13B}}
Left to the reader.
\end{PROOF}

\begin{lemma}\label{6.13C}
Suppose $\aleph_0 = \theta < \chi(*) = \cf(\chi(*))$ and
$\lambda^{\aleph_0} = \lambda^{<\chi(*)}$, $\chi(*) \le \lambda$, 
$\lambda = \lambda^+_1$, and $(*)_{\lambda_1}$ (see below) holds.

\sn
Then
\begin{enumerate}
    \item[$(*)_\lambda$]  We can find $\bfW = \{(\olsi M^\alpha\!,\eta^\alpha) : \alpha < \alpha(*)\}$ and functions $\dot\zeta : \alpha(*) \to S$ and 
    $h : \alpha(*) \to \lambda$ such that:
    \begin{enumerate}
        \item[$(a0)$-$(a2)$]  Are as in \ref{6.5}.
\sn
        \item[$(b0)$-$(b2)$]  As in \ref{6.5}, and even
\sn
        \item[$(c3)$]  If $\dot\zeta(\alpha) = \dot\zeta(\beta)$ then 
        $|M_\alpha| \cap |M_\beta|$ is a bounded subset of $\dot\zeta(\alpha)$.
    \end{enumerate}
\end{enumerate}
\end{lemma}

\begin{PROOF}{\ref{6.13C}}
Left to the reader.
\end{PROOF}

\begin{lemma}\label{6.14}
Suppose that $\lambda$ is a strongly inaccessible uncountable cardinal,
\[
\cf(\lambda) \ge \chi(*) = \cf(\chi(*)) > \theta = \aleph_0,
\]
and let $S \subseteq\lambda$ consist of strong limit singular cardinals of
cofinality $\aleph_0$ and be stationary. \underline{Then} we can find 
$\bfW = \{(\olsi M^\alpha\!,\eta^\alpha) : \alpha < \alpha(*)\}$ and functions\\ $\dot\zeta : \alpha(*) \to S$ and $h : \alpha(*) \to\lambda$ such that:
\mn
\begin{enumerate}
    \item[$(a0)$-$(a2)$]  As in \ref{6.5} (except that $h(\alpha)$ does 
    not only depend on $\dot\zeta(\alpha)$).
\sn
    \item[$(b0),(b3)$]  As in \ref{6.5}.
\sn
    \item[$(b1)^+$]   As in \ref{6.8}.
\sn
    \item[$(c3)^-$]  If $\dot\zeta(\alpha) = \delta = \dot\zeta(\beta)$ 
    then $|M^\alpha_\omega| \cap |M^\beta_\omega| \cap \delta$ is a bounded 
    subset of $\delta$.
\end{enumerate}
\end{lemma}

\begin{remark}\label{6.14X}
1)  See \cite{Sh:45} for a use of what is essentially a weaker version.

\sn
2) We can generalize \ref{6.11A}.
\end{remark}

\begin{PROOF}{\ref{6.14}}
See the proof of \cite[1.10(3)]{Sh:331} (but there
$\sup(N_{\LL\ \RR} \cap \lambda) < \delta$).
\end{PROOF}

\begin{lemma}\label{6.14A}
1) Suppose that $\lambda=\mu^+$, $\mu=\kappa^\theta = 2^\kappa$,
$\theta < \cf(\chi(*)) = \chi(*) < \kappa$, $\kappa$ is strong limit,
$\kappa > \cf(\kappa) = \theta > \aleph_0$, and 
$S \subseteq \{\delta < \lambda : \cf(\delta) = \theta\}$ is stationary.

\underline{Then} we can find $\bfW = \{(\olsi M^\alpha\!,\eta^\alpha) : \alpha <
\alpha(*)\}$ (actually, a sequence), functions $\dot\zeta : \alpha(*) \to S$ 
and $h : \alpha(*) \to \lambda$, and $\LL C_\delta : \delta \in S\RR$ such that:
\begin{enumerate}
    \item[$(a1),(a2)$]  As in \ref{6.5}.
\sn
    \item[$(b0)$]  $\eta^\alpha \ne \eta^\beta$ for $\alpha \ne \beta$.
\sn
    \item[$(b1)$]  If $\{\eta^\alpha \rest i : i < \theta\} \subseteq  M^\beta_\theta$ and $\alpha \ne \beta$ then $\alpha < \beta$ and even $\dot\zeta(\alpha) < \dot\zeta(\beta)$.
\sn
    \item[$(b2)$]  If $\eta^\alpha \rest (j+1) \in M^\beta_\theta$ then $M^\alpha_j \in M^\beta_\theta$.
\sn
    \item[$(c2)$]  $\olsi C = \LL C_\delta : \delta \in S\RR$, $C_\delta$ is
    a club of $\delta$ of order type $\theta$, and every club of $\lambda$ 
    contains $C_\delta$ for stationarily many $\delta \in S$.
\sn
    \item[$(c3)$]  If $\delta \in S$, 
    $C_\delta = \{\gamma_{\delta,i} : i < \theta\}$ is the increasing 
    enumeration, and $\alpha < \alpha^*$ satisfies $\dot\zeta(\alpha) = \delta$,
    then there is $\big\LL\LL\gamma^-_{\alpha,i},\gamma^+_{\alpha,i} \RR : i < \theta \text{ odd} \big\RR$ such that $\gamma^-_{\alpha,i}\in M^\alpha_i$,
    $M^\alpha_i \cap \lambda \subseteq \gamma^+_{\alpha,i}$, $\gamma_{\delta,i}<\gamma^-_{\alpha,i} < \gamma^+_{\alpha,i} < \gamma_{\delta,i+1}$, and
    \sn
    \begin{enumerate}
        \item[$(*)$]  If $\dot\zeta(\alpha) = \dot\zeta(\beta)$ and 
        $\alpha < \beta$ then for every large enough odd $i < \theta$ 
        we have $\gamma^+_{\alpha,i} < \gamma^-_{\beta,i}$ (hence 
        $[\gamma^-_{\alpha,i},\gamma^+_{\alpha,i}) \cap [\gamma^-_{\beta,i}, \gamma^+_{\beta,i}) = \varnothing$) and $[\gamma^-_{\beta,i}, \gamma^+_{\beta,i}) \cap M^\alpha_\theta = \varnothing$.
    \end{enumerate}
\end{enumerate}
\mn
2) In part (1), assume $\theta=\aleph_0$ and $\pp(\kappa)=^+
2^\kappa$. Then the conclusion holds; moreover, {(c3)} (from
\ref{6.12}) does as well.
\end{lemma}

\begin{remark}\label{6.14Z}
The assumption $\pp(\kappa)=2^\kappa$ holds (for example) for 
$\kappa = \beth_\delta$ for a club of $\delta < \omega_1$ 
(and for a club of $\delta < \theta$ when $\aleph_0 < \theta = \cf(\theta) < \kappa$: see \cite[\S5]{Sh:400}).
\end{remark}

\begin{PROOF}{\ref{6.14A}}
1) By \ref{6.3} we can find $\olsi C = \LL C_\delta : \delta \in S \RR$,
$C_\delta$ a club of $\delta$ of order type $\kappa$
such that for any club $C$ of $\lambda$, for stationarily many 
$\delta\in S$, we have $C_\delta \subseteq C$.

\mn
\textbf{First Case}:  Assume $\mu$ ($=2^\kappa$) is regular.

By \cite[Ch.II,5.9]{Sh:g}, we can find an increasing sequence 
$\LL\kappa_i : i < \theta\RR$ of regular cardinals $>\chi(*)$ such that 
$\kappa = \sum\limits_{i<\theta} \kappa_i$ and 
$\prod\limits_{i<\theta}\kappa_i / J^\bd_\theta$ has true cofinality 
$\mu$, and let $\LL f_\epsilon : \epsilon < \mu\RR$ exemplify this. This means
\[
\epsilon < \zeta < \mu\ \Rightarrow\ f_\epsilon < f_\zeta\mod J^\bd_\theta
\]

\mn
and for every $f \in \prod\limits_{i<\theta}\kappa_i$, for some 
$\epsilon < \mu$, we have $f < f_\epsilon \mod J^\bd_\theta$. We may assume 
that if $\epsilon$ is limit and $\bar{f} \rest \epsilon$ has a 
$<_{J^\bd_\theta}$-l.u.b. then $f_\epsilon$ is a $<_{J^\bd_\theta}$-l.u.b., 
and we know that if $\cf(\epsilon) > 2^\theta$ then this holds, 
and that without loss of generality 
$\bigwedge\limits_{i<\theta} \cf(f_\epsilon(i))=\cf(\epsilon)$.
Without loss of generality $\kappa_i>f_\epsilon(i)>\bigcup\limits_{j <
i}\kappa_j$.

We shall define {$\bfW$} later. Let $\St$ be a strategy for Player I 
in the game from \ref{6.5}(a2). By the choice of $\olsi C$, for some 
$\delta\in S$, for every $\alpha\in C_\delta$ of cofinality $>\theta$, 
$\clH_{<\chi(*)}(\alpha)$ is closed under the strategy $\St$. 
Let $C_\delta = \{\alpha_i :  i < \kappa\}$ be increasing continuous. 
For each $\epsilon<\mu$ we choose a play of the game with Player I
using $\St$. For a play, $\LL M^\epsilon_j,\eta^\epsilon_j : j < \theta\RR$ satisfies:

\[
\LL M^\epsilon_j:j \le j_1\RR\in \clH_{<\chi(*)}
(\alpha_{f_\epsilon (j_1)+1}),
\]
\[
\eta^\epsilon_\gamma = \big\LL\cd \big(\alpha_{f_\epsilon(i)},\LL
M^\epsilon_i : i \leq j\RR \big) : j < \gamma \big\RR,
\]
\[
\text{and } \eta^\epsilon_{j+1}\in M^\epsilon_{j+1}.
\]

\mn
Then let $g_\epsilon\in\prod\limits_{i<\theta}\kappa_i$ be
$g_\epsilon(i) = \sup\big(\kappa_i \cap \textstyle\bigcup\limits_{j<\theta} M^\epsilon_j \big)$,

\mn
so for some $\beta_\epsilon\in (\epsilon,\mu)$, we have $g_\epsilon <
f_{\beta_\epsilon} \mod J^\bd_{\theta}$.

On the other hand, if $\cf(\epsilon)=(2^\theta)^+$ then without loss of
generality $\cf\big(f_\epsilon(i)\big) = \cf(\epsilon)$ for every $i<\theta$
(see \cite[Ch.II,\S1]{Sh:g}), so there is $\gamma_\epsilon < \epsilon$ such that
\[
h_\epsilon < f_{\gamma_\epsilon}\mod J^\bd_\theta,\ \text{ where }
h_\epsilon(i) = \sup(f_\epsilon(i) \cap \bigcup\limits_{j<\theta} M^\epsilon_j).
\]
So for some $\gamma(*)<\mu$ we have:
\[
S_\delta[\St] = \{\epsilon < \mu : \cf(\epsilon) = (2^\theta )^+
\text{ and } \gamma_\epsilon = \gamma(*)\} \text{ is stationary}.
\]

\medskip
Now, for each $\delta \in S$ we can consider the set $\bfC_\delta$ of all
possible such $\LL (\olsi M^\epsilon,\eta^\epsilon) : 
\epsilon < \mu\RR$,
where $\olsi M^\epsilon = \LL M^\epsilon_j : j < i\RR$ and
$\eta^\epsilon_\theta$ are as above (letting $\St$ vary on all 
strategies of Player I for which
$[\alpha\in C_\delta\ \wedge\ \cf(\alpha) > \theta]\ \Rightarrow\ 
[\clH_{<\chi(*)}(\alpha) \text{ is closed under }\St])$.

A better way to write the members of $\bfC_\delta$ is $\big\LL\LL (\olsi M^\epsilon_j, \eta^\epsilon_j) : j < \theta\RR : \epsilon < 
\mu \big\RR$, but for $j<\theta$, 
$$f_{\epsilon(1)} \rest j = f_{\epsilon(2)} \rest j\ \Rightarrow\ 
\Big[\olsi M^{\epsilon(1)}_j = M^{\epsilon(2)}_j\ \wedge\ \eta^{\epsilon(1)}_j =
\eta^{\epsilon(2)}_j \Big].$$ 
Actually, it is a function from 
$\{f_\epsilon \rest j : \epsilon < \mu,\ j < \theta\}$ to
$\clH_{<\chi(*)}(\delta)$.  But the domain has power
$\kappa$, the range has power $|\delta|\leq\mu$. So $|\bfC_\delta|
\leq\mu^\kappa=(2^\kappa)^\kappa=2^\kappa=\mu$.

So we can well order $\bfC_{\delta} $ in a sequence of length $\mu $, and
choose by induction on $\epsilon < \mu$ a representative of each for
$\bfW$ satisfying the requirements.

\mn
\textbf{Second case}:  Assume $\mu $ is singular.

So let $\mu = \sum\limits_{\xi<\cf(\mu)}\mu_\xi$ with $\mu_\xi$ regular. 
Without loss of generality 
$$\mu_\xi > \big(\textstyle\sum\{\mu_\epsilon : \epsilon < \xi\} \big)^+ + \cf(\mu)^+.$$ 
We know that $\cf(\mu) > \kappa$, and again by \cite[Ch.VIII,\S1]{Sh:g} 
there are  $\LL\kappa_{\xi,i} : i < \theta\RR$, 
$\LL\kappa_i : i < \theta\RR$ such that
\[
\tcf(\textstyle\prod\limits_{i< \theta}\kappa_{\xi,i}/J^\bd_\theta) = \mu_\xi,
\qquad \tcf (\textstyle\prod\limits_{i<\theta}\kappa_i/J^\bd_\theta) = \cf(\mu),
\]
\[
\kappa^a_i < \kappa_{\xi i} < \kappa^b_i,\quad \kappa^a_i < \kappa_i < \kappa^b_i
\quad \text{ and }\quad i < j\ \Rightarrow\ \kappa^b_i < \kappa^a_j
\]

\mn
(we can even get $\kappa^a_i > \prod\limits_{j<i}\kappa^b_j$ as we can
uniformize on $\xi$).

Let $\LL f^\xi_\epsilon : \epsilon < \mu_\xi\RR$, $\LL f_\epsilon :
\epsilon < \cf(\mu)\RR$ witness the true cofinalities. Now, for every
$f \in \prod\limits_{i<\theta} \kappa_i$ (for simplicity, every $f$ such that 
$f(i) > \sum\limits_{j<i} \kappa_j$ and $\bigwedge\limits_i\cf(f(i)) = (2^\theta)^+$)
and $\xi$ we can repeat the previous argument for $\LL f+f^\xi_\epsilon:
\epsilon<\mu_\epsilon\RR$. After ``cleaning inside," replacing by a
subset of power $\mu_\xi$, we find a common bound below
$\prod\limits_{i<\theta}\kappa_i$ and below $\prod f$, and we can uniformize
on $\xi$.

Thus we apply $\cf(\epsilon)=(2^\theta)^+$ on every $f_\epsilon$, and use
the same argument on the bound we have just gotten.

\mn
2) Should be clear.
\end{PROOF}

\noindent
Similarly to \ref{6.11}, with $\omega^2$ for $\theta$ 
(not a cardinal!) we have: 

\begin{claim}\label{6.14B}
Suppose that
\mn
\begin{enumerate}
    \item[$(*)$]   $\lambda$ is a regular cardinal, $\theta = \aleph_0$, 
    $\mu = \mu^{<\chi(*)} < \lambda \le 2^\mu$, \\
    $S \subseteq \{\delta < \lambda : \cf(\delta) = \aleph_0\}$ is 
    stationary, and $\aleph_0 < \chi(*) = \cf(\chi(*))$.
\end{enumerate}
\mn
\underline{Then} we can find
\[
\bfW = \{(\olsi M^\alpha\!,\eta^\alpha):\alpha<\alpha(*)\}
\]
and functions
\[
\dot\zeta:\alpha(*)\to S\ \text{ and }\ h:\alpha(*)
\to\lambda
\]

\mn
such that:
\mn
\begin{enumerate}
    \item[$(a0)$]  As in \ref{6.5}. 
\sn
    \item[$(a1)$]  $\olsi M^\alpha = \LL M^\alpha_i : i \leq \omega^2\RR$ is an increasing continuous elementary 
    chain,\footnote{$\tau(M^\alpha_i)$, the vocabulary, may be increasing too and belongs to $\clH_{< \chi(*)}(\chi(*))$.}
    each $M^\alpha_i$ is a model belonging to $\clH_{<\chi(*)}(\lambda)$ 
    [so necessarily has cardinality $<\chi(*)$], $M^\alpha_i\cap\chi(*)$ is an
    ordinal, $[\chi(*) = \chi^+\ \Rightarrow\ \chi+1 \subseteq M^\alpha_i]$, 
    $\eta^\alpha\in {}^{\omega^2}\lambda$ is increasing with limit 
    $\dot\zeta(\alpha)\in S$, $\eta^\alpha \rest i \in M^\alpha_{i+1}$, 
    $M^\alpha_i$ belongs to $\clH_{<\chi(*)}(\eta^\alpha(i))$, and 
    $\LL M^\alpha_i : i \leq j\RR$ belongs to $M^\alpha_{j+1}$.
\sn
    \item[$(a2)$]  Like \ref{6.5} (with $\omega^2$ instead $\theta)$.
\sn
    \item[$(b0)$-$(b2)$]  As in \ref{6.5}.
\sn
    \item[$(b1)^*$] As in \ref{6.11}.
\sn
    \item[$(c1)$]  If $\dot\zeta(\alpha) = \dot\zeta(\beta)$ then 
    $M^\alpha_{\omega^2}\cap \mu=M^\beta_{\omega^2} \cap\mu$, there is an
    isomorphism $h_{\alpha,\beta}$ from $M^\alpha_{\omega^2}$ onto $M^\beta_{\omega^2}$ mapping $\eta^\alpha(i)$ to $\eta^\beta(i)$ and
    $M^\alpha_i$ to $M^\beta_i$ for $i<\omega^2$, and 
    $h_{\alpha,\beta} \rest \big(|M^\alpha_{\omega^2}| \cap |M^\beta_{\omega^2}| \big)$ is the identity.
\sn
    \item[$(c2)$]  As in \ref{6.11}, using 
    $\LL M^\alpha_{\omega n} : n < \omega\RR$.
\sn
    \item[$(c3)$]  As in \ref{6.12}, assuming $\lambda = \mu^+$.
\sn 
    \item[$(c4)$]  $\eta^\alpha(i) > \sup\!\big(|M^\alpha_i|\cap\lambda\big)$ 
    (so $\sup\!\big(|M^\alpha_{\omega (n+1)}| \cap \lambda\big) = 
    \bigcup\limits_\ell \eta^\alpha(\omega n+\ell)$).
\end{enumerate}
\end{claim}

\begin{PROOF}{\ref{6.14B}}
We use $\LL\olsi M^{\alpha,0}:\alpha<\alpha(*)\RR$, which we got
in \ref{6.11}. Now for each $\alpha$ we look at $\bigcup\limits_{n<\omega}
M^{\alpha,0}_n$ as an elementary submodel of $(\clH_{<\chi(*)}
(\lambda),\in)$ with a function $\St$ (intended as a strategy
for Player I in the play for (a2) above).

Play in $\bigcup\limits_{n<\omega} M^{\alpha,0}_n$ and get
\[
\begin{array}{c}
\LL M^\alpha_i,\eta^\alpha(i):i<\omega n\RR\in M^{\alpha ,0}_n,\\
\sup\{\eta^\alpha(i):i<\omega n\}\in M^{\alpha ,0}_{n+1},\\
\eta^\alpha(\omega n)>\sup(M^{\alpha ,0}_n\cap\lambda).
\end{array}
\]
\end{PROOF}

\subsection{Black Boxes: third round}

\begin{lemma}\label{6.14C}
Assume that $\lambda \ge \chi(*)>\theta$ are regular cardinals, $S
\subseteq\{\delta<\lambda :\cf(\delta)=\theta\}$ is a stationary set,
$\lambda^{<\chi(*)} = \lambda$, and the conclusion of \ref{6.14A} holds for
them. \underline{Then} it holds for $\lambda^+$ as well as $\lambda$.
\end{lemma}

\begin{PROOF}{\ref{6.14C}}
By \cite[2.10(2)]{Sh:331} (or see \cite{Sh:365}) we know
\mn
\begin{enumerate}
    \item[$(*)$]  There are $\LL C_\delta : \delta < \lambda^+,\ 
    \cf(\delta) = \theta\RR$, $\LL e_\alpha : \alpha < \lambda^+\RR$ such that:
    \begin{enumerate}
        \item[$(i)$]  $C_\delta$ is a club of $\delta$ of order type $\theta$  such that
        $$\alpha\in C_\delta\ \wedge\ \alpha > \sup(C_\delta \cap \alpha) \ \Rightarrow\ \cf(\alpha) = \lambda.$$
        
        \item[$(ii)$]  $e_\alpha$ is a club of $\alpha$ of order type $\cf(\alpha)$; we let 
        $e_\alpha = \{\beta^\alpha_i : i < \cf(\alpha)\}$ (increasing continuous).
\sn
        \item[$(iii)$]  If $E$ is a club of $\lambda^+$ then for stationarily 
        many $\delta < \lambda^+$ we have $\cf(\delta) = \theta$, 
        $C_\delta \subseteq E$, and the set
        \[
            \{i < \lambda : \mbox{for every }\alpha\in C_\delta,\ \cf(\alpha) = \lambda\ \Rightarrow\ \beta^\alpha_{i+1}\in E\}
        \]
        is unbounded in $\lambda$.
    \end{enumerate}
\end{enumerate}
\mn
Now copying the black box of $\lambda$ on each $\delta < \lambda^+$ with
$\cf(\delta) = \theta$, we can finish easily.
\end{PROOF}

\begin{lemma}
\label{6.14D}
If $\lambda$, $\mu$, $\kappa$, $\theta$, $\chi(*)$, $S$ are as in
\ref{6.14A}, and
\[
\alpha < \chi(*)\ \Rightarrow\ |\alpha|^\theta < \chi(*)
\]
\underline{then} there is a stationary 
$S^* \subseteq [\lambda]^{<\chi(*)}$
and a one-to-one function $\cd$ from $S^*$ to $\lambda$ such that
\[
[A \in S^*\ \wedge\ B \in S^*\ \wedge\ A\subsetneq B]\ 
\Rightarrow\ \cd(A)\in B.
\]
\end{lemma}

\begin{remark}\label{6.14Dx}
This gives another positive instance to a problem of
Zwicker. (See \cite{Sh:247}.)
\end{remark}

\begin{PROOF}{\ref{6.14D}}
Similar to the proof of \ref{6.14A}, only choose
$\cd : [\lambda]^{<\chi(*)} \to \lambda$
one-to-one, and then define
$S^* \cap [\alpha]^{<\chi(*)}$
by induction on $\alpha$.
\end{PROOF}

\begin{problem}
\label{6.14E}
1) Can we prove in ZFC that for some regular $\lambda>\theta$:
\mn
\begin{enumerate}
    \item[$(*)_{\lambda,\theta,\chi(*)}$] We can define, for 
    $\alpha \in S^\lambda_\theta = \{\delta < \lambda : \aleph_0 \le \cf(\delta) = \theta\}$, a model $M_\alpha$ with a countable vocabulary and universe an unbounded subset of
$\alpha$ of power $<\chi(*)$, {such that} $M_\delta\cap \chi(*)$ is an ordinal
such that for every model $M$ with countable vocabulary and universe
$\lambda$, for some\footnote{Equivalently, stationarily many.} $\delta\in
S^\lambda_\kappa$, we have $M_\delta \subseteq M$.
\end{enumerate}
\mn
2) The same. dealing with relational vocabularies only. (We call it
$(*)^{\rel}_{\lambda,\theta ,\kappa}$.)
\end{problem}

\begin{remark}\label{6.14F}
Note that by \ref{6.4}, if $(*)_{\lambda,\theta,\kappa}$ and
$\mu = \cf(\mu) > \lambda$ then $(*)_{\mu^+,\theta,\kappa}$.
\end{remark}

\medskip
\centerline {$* \qquad * \qquad *$}
\bigskip

Now (in \ref{6.15}--\ref{6.19}) we return to black boxes for singular
$\lambda$: i.e. we deal with the case $\cf(\lambda) \le \theta$.

\begin{lemma}\label{6.15}
Suppose that $\lambda^\theta = \lambda^{<\chi(*)}$, $\lambda$ is a singular
cardinal, $\theta$ is regular, and $\chi(*)$ is regular $>\theta$. 

Assume further
\begin{enumerate}
    \item[$(\alpha)$]  $\cf(\lambda) \le \theta$
\sn
    \item[$(\beta)$]  $\lambda = \sum\limits_{i\in w} \mu_i$, $|w| \le \theta$,
    $w \subseteq \theta^+$ (usually $w = \cf(\lambda)$), 
    $[i < j\ \Rightarrow\ \mu_i < \mu_j]$, each $\mu_i$ is regular $<\lambda$, and
    \[
        \cf(\lambda) > \aleph_0\ \wedge\ \cf(\lambda) = \theta\ \Rightarrow\ w = \cf(\lambda).
    \]

    \item[$(\gamma)$]  $\mu>\lambda$, $\mu$ is a regular cardinal, $D$ is a uniform filter on $w$ (so $\{\alpha \in w : \alpha > \beta\}\in D$ for 
    each $\beta \in w$), $\mu$ is the true cofinality of 
    $\prod\limits_{i\in w}(\mu_i,<)/D$ (see \cite[3.7(2)=Lc18]{Sh:E62} or \cite{Sh:g}).
\sn
    \item[$(\delta)$]  $\bar{f} = \LL f_i / D : i < \mu\RR$ exemplifies ``the 
    true cofinality of $\prod\limits_i (\mu_i,<)/D$ is $\mu$:'' i.e.
    \[
    \begin{array}{c}
        \alpha < \beta < \lambda\ \Rightarrow\ f_\alpha/D < f_\beta/D,\\
        f \in \prod\limits_i\mu_i\ \Rightarrow\ \bigvee\limits_\alpha 
        f/D < f_\alpha/D.
    \end{array}
    \]

    \item[$(\eps)$]  $S\subseteq\{\delta<\mu:\cf(\delta)=\theta\}$ 
    is good for $(\mu,\theta,\chi(*))$.
\sn
    \item[$(\zeta)$] If $\theta > \cf(\lambda)$, $\delta \in S$, then 
    for some $A_\delta \in D$ and unbounded $B_\delta \subseteq \delta$ we have
    \[
        \alpha,\beta\in B_\delta\ \wedge\ \alpha<\beta\ \wedge\ i\in A_\delta\ \Rightarrow\ f_\alpha(i) < f_\beta(i)
    \]
    i.e. $\LL f_\alpha \rest A_\delta : \alpha \in B_\delta\RR$ is $<$-increasing.
\end{enumerate}
\mn
\underline{Then} we can find 
$\bfW = \{(\olsi M^\alpha\!,\eta^\alpha) : \alpha < \alpha(*)\}$ (pedantically, a sequence) and functions $\dot\zeta : \alpha(*) \to S$ and $h : \alpha(*) \to \mu$ such that:
\mn
\begin{enumerate}
    \item[$(a0)$-$(a2)$]  As in \ref{6.5}, except that we replace $(a1)(*)$ by
    \begin{enumerate}
        \item[$(*)'$] 
        \begin{enumerate}
            \item[$(i)$] $\eta^\alpha\in {}^\theta\!\lambda$
\sn
            \item[$(ii)$] If $i < \cf(\lambda)$ then $\sup(\mu_i \cap \Rang(\eta^\alpha)) = \sup(\mu_i \cap M^\alpha_\theta)$.
\sn
            \item[$(iii)$] If $\xi < \dot\zeta(\alpha)$ then 
            $$f_\xi/E < \big\LL \sup(\mu_i \cap M^\alpha_\theta) : i < \cf(\lambda) \big\RR/E \le f_{\dot\zeta(\alpha)}/E.$$
        \end{enumerate}
    \end{enumerate}
\sn
    \item[$(b0)$-$(b3)$]  As in \ref{6.5}.
\end{enumerate}
\end{lemma}

\begin{PROOF}{\ref{6.15}}
For $A\subseteq\theta$ of cardinality $\theta$, let
$\cd^A_{\lambda,\chi(*)} : \clH_{<\chi(*)}(\lambda) \to {}^A\!\lambda$ be one-to-one and $G:\lambda\to\lambda$ be such
that for $\gamma$ divisible by $|\gamma|$ and $\alpha < \gamma \le \lambda$
({and} $\mu \ge \aleph_0$), the set $\{\beta<\gamma:G(\beta)=\alpha\}$ is unbounded
in $\gamma$ and of order type $\gamma$. Let $\bar{A}=\LL
A_i:i<\theta\RR$ be a sequence of pairwise disjoint subsets of $\theta$
each of cardinality $\theta$. 

For $\delta\in S$, let
\begin{align*}
  \bfW^0_\delta = \Big\{\big(\olsi M,\eta\big) : &\ \olsi M,\eta 
  \text{ satisfy {(a1)}, and for some} \\
  &\ y \in \clH_{<\chi(*)}(\lambda),\text{ for every }i < \theta, \text{ we have}\\
  &\ \big\LL G(\eta(i)) : i \in A_j \big\RR = \cd^A_{\lambda,\chi(*)}\big(\LL
\olsi M \rest j,\eta \rest j,y\RR \big) \Big\}.
\end{align*}
The rest is as before.
\end{PROOF}

\begin{claim}
\label{6.16}
Suppose that $\lambda^\theta = \lambda^{<\chi(*)}$, $\lambda$ is singular,
$\theta$ and $\chi(*)$ are regular, and $\chi(*) > \theta$.

\sn
1) If $(\forall \alpha < \lambda) \big[|\alpha|^{<\chi(*)} < \lambda \big]$ then
by $\lambda^\theta=\lambda^{<\chi(*)}$ we know that either
$\cf(\lambda) \ge \chi(*)$ (and so lemma \ref{6.9} applies) or
$\cf(\lambda) \le \theta$.

\sn
2)  We can find regular $\mu_i$ (for $i < \cf(\lambda)$) increasing with
$i$ such that $\lambda=\sum\limits_{i<\cf(\lambda)}\mu_i$.

\sn
3)  For $\LL\mu_i : i \in w\RR$ as in \ref{6.15}($\beta$), we can
find $D,\mu,\bar{f}$ as in \ref{6.15}($\gamma$),($\delta$) with $D$ the
co-bounded filter plus one unbounded subset of $\omega$.

\sn
4) For $\LL\mu_i : i \in w\RR$, $D,\mu,\bar{f}$ as in 
($\beta$),($\gamma$),($\delta$) of \ref{6.15}, we can find $\mu$ and 
pairwise disjoint $S \subseteq \mu$ as required in \ref{6.15}($\delta$)($\eps$)
provided that $\theta > \cf(\lambda) \Rightarrow 2^\theta < \mu$ [equivalently,
$<\lambda$].

\sn
5)  If $\cf(\lambda)>\aleph_0$, $(\forall \alpha < \lambda)
\big[|\alpha|^{\cf(\lambda)} < \lambda \big]$, and 
$\lambda < \mu = \cf(\mu) \le \lambda^{\cf(\lambda)}$ \underline{then} we 
can find $\LL\mu_i : i < \cf(\lambda)\RR$ and the co-bounded filter $D$ on 
$\cf(\lambda)$ as required in \ref{6.14}$(\beta),(\gamma)$.
\end{claim}

\begin{PROOF}{\ref{6.16}}
Now (1)-(3) are trivial; for (5) see \cite[\S9]{Sh:345}. As for (4), we should
recall that \cite[\S5]{Sh:345} actually says:

\begin{fact}\label{6.17A}
If $\LL \mu_i : i \in w\RR,\bar{f},D$ are as in \ref{6.15}, then
\begin{align*}
S = \big\{\delta < \mu : &\ \cf(\delta) = \theta \text{ and there are }
A_\delta \in D \text{ and unbounded } B_\delta\subseteq\delta \\
  &\text{ such that }[\alpha,\beta\in B_\delta\wedge
\alpha<\beta \wedge i\in A_\delta \Rightarrow f_\alpha(i)<f_\beta(i)] \big\}.
\end{align*}
is good for $(\mu,\theta,\chi(*))$.
\end{fact}
\end{PROOF}

\begin{lemma}\label{6.18}
Let $\chi(1) = \chi(*) + (<\chi(*))^\theta$. 

In \ref{6.15}, if $\lambda^\theta=\lambda^{\chi(1)}$, we can strengthen
{(b1)} to {(b1)}$^+$ (of \ref{6.8}).
\end{lemma}

\begin{PROOF}{\ref{6.18}}
Combine proofs of \ref{6.15}, \ref{6.8}.
\end{PROOF}

\begin{lemma}\label{6.19}
$\frac{3.17}{3.11} \times 3.29$ and $\frac{3.19}{3.11} \times 3.37$
 hold (we need also the parallel to \ref{6.14A}).
\end{lemma}

\begin{PROOF}{\ref{6.19}}
Left to the reader.
\end{PROOF}

\bigskip

\subsection{Conclusion}
Now we draw some conclusions.

\sn
The first, \ref{6.20}, gives what we need in \ref{5.4} (so \ref{5.1}).
\begin{conclusion}
\label{6.20}
Suppose $\lambda^\theta = \lambda^{<\chi(*)}$, $\cf(\lambda) \ge \chi(*) + \theta^+$, $\theta = \cf(\theta) < \chi(*) = \cf(\chi(*))$.
\underline{Then} we can find
\[
\bfW = \{(\olsi M^\alpha\!,\eta^\alpha) : \alpha < \alpha(*)\},
\]
where
\[
M^\alpha_i = (N^\alpha_i,A^\alpha_i,B^\alpha_i),\ 
A^\alpha_i \subseteq \lambda \cap |N^\alpha_i|,\ B^\alpha_i \subseteq
\lambda \cap |N^\alpha_i|,\ A^\alpha_i \ne B^\alpha_i,
\]
and functions $\dot\zeta,h$ such that:
\mn
\begin{enumerate}
    \item[$(a0),(a1)$]  As in \ref{6.5}.
\sn
    \item[$(a2)$]  As in \ref{6.5}, except that in the game, Player I can 
    choose $M_i$ only as above.
\sn
    \item[$(b0),(b1),(b2)$] As in \ref{6.5}.
\sn
    \item[$(b1)''$]  If $\{\eta^\alpha \rest i : i < \theta\} \subseteq M^\beta$ but $\alpha < \beta$ (so $\beta < \alpha + ({<}\chi(*))^\theta$) \underline{then}:
    \[
    \begin{array}{c}
        A^\alpha_\theta \cap \big(|M^\alpha_\theta|\cap |M^\beta_\theta|\big) \neq
        B^\beta_\theta \cap \big(|M^\alpha_\theta| \cap |M^\beta_\theta|\big),\\
        B^\alpha_\theta \cap \big(|M^\alpha_\theta|\cap |M^\beta_\theta|\big) \neq
        A^\beta_\theta \cap \big(|M^\alpha_\theta|\cap|M^\beta_\theta|\big).
    \end{array}
    \]
\end{enumerate}
\end{conclusion}

\begin{PROOF}{\ref{6.20}}
First assume $\lambda$ is regular, and $\bfW=\{(\olsi M^\alpha\!,\eta^\alpha):
\alpha <\alpha(*)\},\dot\zeta,h^*$ be as in the conclusion of \ref{6.5} 
(with $h^*$ here standing in for $h$ there). Let
$w=\{\cd(\alpha,\beta):\alpha,\beta<\lambda\}$, and
$G_1,G_2:w\to\lambda$ be such that for $\alpha\in E$,
$\alpha=\cd(G_1(\alpha),G_2(\alpha))$. 

Let
\begin{align*}
Y = \big\{\alpha < \alpha(*) : &\ \olsi M^\alpha_i \text{ has the form } 
(N^\alpha_i,A^\alpha_i,B^\alpha_i),\\
  &\ A^\alpha_i,B^\alpha_i \text{ distinct subsets of }\lambda \cap
|N^\alpha_i|\\
  &\text{ (equivalently, monadic relations), and}\\ 
  &\ G_2(h(\alpha))= \min(A^\alpha_i\setminus
B^\alpha_i \cup B^\alpha_i \setminus A^\alpha_i) \big\}.
\end{align*}

\mn
Now we let
\[
\bfW^* = \big\{(\olsi M^\alpha\!,\eta^\alpha) : \alpha\in Y \big\},\ 
\dot\zeta^* = \dot\zeta \rest Y,\text{ and } h = G_1 \circ h^*.
\]

\mn
They exemplify that \ref{6.20} holds.

What if $\lambda$ is singular? Still, $\cf(\lambda) \ge \chi(*) + \theta^*$, 
and we can just use \ref{6.9} instead \ref{6.5}.
\end{PROOF}

\begin{claim}\label{6.21}
1) In \ref{6.5}, if $\lambda=\lambda^{<\chi(*)}$ we can let 
$h : S \to \clH_{<\chi(*)}(\lambda)$ be onto. Generally, we can
still make $\Rang(h)$ be $\subseteq A$ whenever $|A| = \lambda$.

\sn
2) In \ref{6.5}, by its proof, whenever $S'\subseteq S$ is stationary, and $$\textstyle\bigwedge\limits_{\zeta}\big[h^{\text{--}1}(\zeta) \cap S' 
\text{ stationary} \big]$$
then
$\big\{(\olsi M^\alpha\!,\eta^\alpha) : \alpha<\alpha(*),\ \dot\zeta(\alpha) \in
S' \big\}$ satisfies the same conclusion.

\sn
3) For any unbounded $a \subseteq\theta$, we can let Player I also choose
$\eta(i)$ for $i \in \theta \setminus a$ without changing our
conclusions.

\sn
4) Similar statements hold for the parallel claims.

\sn
5) It is natural to have $\chi(*)=\chi^+$.
\end{claim}

\begin{PROOF}{\ref{6.21}}
Straightforward.
\end{PROOF}

\begin{fact}\label{6.22}
We can make the following changes in {(a1)}, {(a2)} of \ref{6.5}
(and in all similar lemmas here) getting equivalent statements:
\mn
\begin{enumerate}
    \item[$(*)$]  $M^\alpha_i\in \clH_{<\chi(*)}(\lambda+\lambda)$: 
    in the game, for some arbitrary $\lambda^* \ge \lambda$ (but fixed 
    during the game) Player I chooses the $M^\alpha_i\in \clH(\lambda^*)$ 
    of cardinality $<\chi(*)$, and in the end instead of 
    ``$\bigwedge\limits_{i<\theta} [M_i = M^\alpha_i]$" we have
    \begin{itemize}
        \item There is an isomorphism from $M_\theta$ onto $M^\alpha_\theta$ taking $M_i$ onto $M^\alpha_i$, is the identity on 
        $M_\theta\cap \clH_{<\chi(*)}(\lambda)$, maps 
        $|M_\theta| \setminus \clH(\lambda)$ into $\clH_{<\chi(*)}(\lambda+\lambda) \setminus \clH_{<\chi(*)}(\lambda)$, and preserves $\in$, $\notin$, and `[is/is not] an ordinal'.
    \end{itemize}
\end{enumerate}
\end{fact}

\begin{exercise}\label{6.23}
If $D$ is a normal fine filter on $\clP(\mu)$, $\lambda$ is regular,
$\lambda\le \mu$,\\ $S \subseteq \{\delta < \lambda : \cf(\delta) = \theta\}$ is
stationary, and furthermore
\begin{enumerate}
    \item[$(*)_{D,S}$]   $\{ a \subseteq \mu : \sup(a\cap\lambda) \in S\} 
    \ne \varnothing\mod D$. 
\end{enumerate}
\mn
\underline{then} we can partition $S$ to $\lambda$ stationary disjoint subsets $\LL S_i : i < \lambda\RR$ such that $i < \lambda \Rightarrow (*)_{D,S_i}$.

\medskip
[Hint: like the proof of \ref{6.2}.]
\end{exercise}

\begin{notation}\label{6.24}
1)  Let $\kappa$ be an uncountable regular cardinal. We let
$\seq^\alpha_{<\kappa}(\cA)$ (where $\cA$ is an expansion of a
submodel of some $\clH_{\le \mu}(\lambda)$ with $|\tau(\cA)| \le \chi$) 
be the set of sequences $\LL M_i : i < \alpha\RR$ which
are increasing continuous with $M_i\prec \cA$, $\|M_i\|<\kappa$, $M_i\cap
\kappa\in\kappa$, $\kappa = \kappa_1^+ \Rightarrow \kappa_1+1 \subseteq
M_i$, and $\LL M_j : j \le i\RR \in M_{i+1}$. (If $\alpha=\delta$ is
limit, $M_\delta \defeq \bigcup\limits_{i<\delta} M_i)$.

\sn
2) If $\kappa=\kappa^+_1$, we may write $\le \kappa_1$ instead $<\kappa$.
\end{notation}

\noindent
We repeat the definition of filters introduced in \cite[Definition 3.2]{Sh:52}.

\begin{definition}\label{6.25}
1)  $\cE^\theta_{<\kappa}(A)$ is a filter on $[A]^{<\kappa}$ defined as follows:
$Y \in \cE^\theta_{<\kappa}(A)$ \underline{iff} for (every) $\chi$
large enough, for some $x\in \clH(\chi)$, the set
$$\big\{\big(\textstyle\bigcup\limits_{i<\theta} M_i \big) \cap A : 
\LL M_i : i < \theta\RR \in \seq^\theta_{<\kappa} \big(\clH(\chi),\in,x\big) \big\}$$ 
is included in $Y$.
\end{definition}

\begin{exercise}\label{6.26}
Let $\lambda$, $\kappa$, $\theta$, and $Y\subseteq [\lambda]^{<\kappa}$ be
given. Then
\[
\mathbf{(a)}\ \Rightarrow\ \mathbf{(b)}\ \Rightarrow\ \mathbf{(c)},
\]
where
\begin{enumerate}
    \item[\textbf{(a)}]  For some $\bfW = \{(\olsi M^\alpha\!,\eta^\alpha) : \alpha < \alpha(*)\}$, $\dot\zeta$, and $h$ satisfying \ref{6.5}, we have
    \[
        Y = \{M^\alpha_\theta \cap \lambda : \alpha < \alpha(*)\}
    \]
    and

    \sn
    \begin{enumerate}
        \item[$(*)$]  $\alpha \ne \beta \wedge \bigwedge\limits_{i<\theta} [\eta^\alpha_i\in M^\beta_\theta] \Rightarrow \alpha < \beta$.
    \end{enumerate}
\sn
    \item[\textbf{(b)}]  $\diamondsuit_{E^\theta_{<\kappa}(\lambda)}$ holds.
\sn
    \item[\textbf{(c)}]  Like \textbf{(a)}, but without $(*)$.
\end{enumerate}
\end{exercise}

\begin{exercise}\label{6.26A}
If $\lambda^{2^\kappa} = \lambda$ and $\theta\leq\kappa$ then
$\diamondsuit_{E^\theta_{<\kappa}}$. (Main case: $\kappa=\theta$.)
\end{exercise}

\begin{exercise}\label{6.26B}
If $\lambda=\mu^+$, $\lambda^\kappa=\lambda$, $\theta=\aleph_0$, $\kappa=
\kappa^\theta$, then there is a coding set with diamond (see \cite{Sh:247}).
\end{exercise}

\begin{exercise}\label{6.26D}
Suppose that $\cf(\lambda) > \aleph_0$, $2^\lambda = \lambda^{\cf(\lambda)}$,
$\chi(*) \ge \theta > \cf(\lambda)$,\\ $(\forall\alpha < \lambda)
\big[|\alpha|^{\chi(*)} < \lambda \big]$, and $\gC$ is a model expanding 
$(\clH_{<\chi(*)}(\lambda),\in),|\tau(\gC)| \le \aleph_0$. 
\underline{Then} we can find $\{\olsi M^\alpha: \alpha<\alpha(*)\}$ such that:
\mn
\begin{enumerate}
    \item[$(i)$]  $\olsi M^\alpha = \LL M^\alpha_i : i < \sigma\RR$, 
    $M^\alpha_i \in \clH_{<\chi(*)}(\lambda)$, $M^\alpha_i \cap \chi(*)$ 
    is an ordinal,\\ $M^\alpha_i \rest \tau (\gC) \prec \gC$, 
    $[i < j\ \Rightarrow\ M^\alpha_i \prec M^\alpha_j]$, and 
    $\LL M^\alpha_j : j \le i\RR \in M^\alpha_{i+1}$.
\sn
    \item[$(ii)$]  If $f_n$ is a $k_n$-place function from $\lambda$ to $\clH_{<\chi(*)}(\lambda)$ then for some $\alpha$, \\
    $M^\alpha_\sigma \prec (\gC,f_n)_{n<\omega}$.
\end{enumerate}
\end{exercise}

\begin{exercise}\label{6.26E}
Suppose $\theta = \cf(\mu) < \mu$, $(\forall\alpha < \mu) \big[|\alpha|^\theta < \mu \big]$,
$2^\mu = \mu^\theta$ and $\lambda = (2^\mu)^+$, and 
$S \subseteq \{\delta < \lambda : \cf(\delta) = \theta\}$. Let 
$\mu = \sum\limits_{i<\theta}\mu_i$, $\mu_i$
regular strictly increasing, and $\cf(\prod\mu_i/E) = 2^\mu$. Then we can find
\[
\bfW = \big\{(\olsi M^\alpha\!,\eta^\alpha) : \alpha < \alpha(*) \big\},\quad 
\dot\zeta : \alpha(*) \to S,\quad h : \alpha(*) \to \lambda
\]

\mn
such that:
\mn
\begin{enumerate}
    \item[$(*)$]  For $\delta\in S$ there is a club $C_\delta$ of $\delta$ of
    order type $\theta$ such that
    \[
        \alpha \in C_\delta \wedge \otp(\alpha \cap C_\delta) = \gamma + 1 \
        \Rightarrow \ \cf(\alpha) = \mu_\gamma.
    \]
\end{enumerate}
\end{exercise}

\begin{remark}\label{6.26F}
We do not know if the existence of a Black Box for $\lambda^+$ with $h$
one-to-one follows from ZFC (of course it is a consequence of
$\diamondsuit$). On the other hand, it is difficult to get rid of such a
Black Box (i.e., prove the consistency of non-existence).

If $\lambda = \lambda^{<\lambda}$ then we have $h : S \to \lambda$,
$S \subseteq \{\delta<\lambda^+:\cf(\delta)<\lambda\}$ such that $C_\delta$ is
a club of $\delta$, $\otp (C_\delta)=\cf(\delta)$ and
\[
{(\forall \alpha \in C_\delta)(\forall \text{clubs }C \subseteq \alpha)} \big[\cf(\alpha) > \aleph_0\ \wedge\ \min\limits_{C' \text{ club of } C_\alpha}
\sup (h \rest C') = \otp(C \cap \alpha) \big].
\]

\mn
This is hard to get rid of (i.e. it is hard to find a forcing notion making it
no longer a black box without collapsing too many cardinals); compare with
Mekler-Shelah \cite{Sh:274}.
\end{remark}

\sn
Recall

\begin{definition}\label{4.2A}
For $\lambda > \theta = \cf(\theta) > \aleph_0$ and stationary $S \subseteq [\lambda]^{<\theta}$, let $\diamondsuit_S$ be defined as follows:

If $\tau$ is a countable vocabulary, {then} there is a diamond sequence $\olsi N = \LL N_a : a \in S\RR$ witnessing it, which means
\begin{itemize}
    \item If $N$ is a $\tau$-model with universe $\lambda$ \underline{then} for stationarily many $a \in S$ we have $N_a \prec N$.
\end{itemize}
(Pedantically, we only consider $a \in S\setminus \varnothing$.)
\end{definition}

\newpage

\section {On Partitions to stationary sets}

We present some results on the club filter on
$[\kappa]^{\aleph_0}$ and $[\kappa]^{\theta}$ and some relatives, and on
$\diamondsuit$ (see Definition \cite[4.6=Ld12]{Sh:E62} or 
\ref{4.4new}(2) here). There are overlaps of the
claims, hence redundant parts, but we believe they are still of some interest.

\begin{claim}\label{4.1new}
Assume $\kappa$ is a cardinal $>\aleph_1$. \underline{Then}
$[\kappa]^{\aleph_0}$ can be partitioned to $\kappa^{\aleph_0}$ (pairwise
disjoint) stationary sets.
\end{claim}

\begin{PROOF}{\ref{4.1new}}
Follows by \ref{4.2new} below. In detail,
let $\tau$ be the vocabulary $\{c_n : n < \omega\}$ where each $c_n$ 
is an individual constant.  By \ref{4.2new} below there is a sequence
$\olsi M =  \LL M_u : u \in [\kappa]^{\aleph_0} \RR$
of $\tau$-models, with $M_u$ having universe $u$ 
such that $\olsi M$ is a diamond sequence.

For each $\eta \in {}^{\omega}\!\lambda$, let $\clS_\eta$ be the set
$u \in [\kappa]^{\aleph_0}$ such that for every $n < \omega$ we have
$c^{M_u}_n = \eta(n)$.

By the choice of $\olsi M$, each set $\clS_\eta$ is necessarily
a stationary subset of $[\kappa]^{\aleph_0}$, and trivially those
sets are pairwise disjoint.
\end{PROOF}

\begin{claim}\label{4.2new}
Let $\kappa>\aleph_1$. Then we have diamond on $[\kappa]^{\aleph_0}$ (modulo
the filter of clubs on it: see \ref{4.2A} or
\cite[4.6=Ld12]{Sh:E62}), and we can find 
$A_\alpha\subseteq [\kappa]^{\aleph_0}$ for $\alpha<\lambda \defeq 
2^{\kappa^{\aleph_0}}$ such that each is stationary but the intersection of any
two is not.
\end{claim}

\begin{PROOF}{\ref{4.2new}}
The existence of the $A_\alpha$-s for $\alpha < \lambda$
follows from the first result.
Let $\tau$ be a countable vocabulary and $\tau_1 = \tau \cup \{<\}$. 
First we prove it when $\kappa = \aleph_2$. Without loss of generality $\kappa \leq 2^{\aleph_0}$, as otherwise the claim follows by 
\ref{6.12}(3), with $(\aleph_2,\aleph_1,\aleph_0)$ here standing in for $(\lambda,\mu,\kappa)$ there.
Let $\omega \setminus \{0\}$ be the disjoint union of $s_n$ for
$n < \omega$, each $s_n$ is infinite with the first element $>n+3$ when
$n>0$. By \cite[2.2]{Sh:331} or \cite{Sh:365}=\cite[Ch.III]{Sh:g} we can 
choose a sequence $\LL C_\delta : \delta \in S^2_0\RR$ which guesses 
clubs (where $S^2_0 = \{\delta < \omega_2 : \cf(\delta) = \aleph_0\}$)
such that $C_\delta \subseteq \delta = \sup(C_\delta)$ 
has order type $\omega$.

Let $\big\LL (\gA^\zeta,\bar{\alpha}^\zeta) : \zeta < 2^{\aleph_0}\big\RR$
list the pairs $(\gA,\bar{\alpha})$ without repetitions, with $\gA$ 
a model with vocabulary $\tau_1$ and universe a limit countable 
ordinal, and $\bar{\alpha} = \LL \alpha_n : n < \omega \RR$ 
an increasing sequence of ordinals with limit
$\sup(\gA)$ and $\gA \!\rest\! \alpha_n \prec \gA$.
Let $E_n$ be the following equivalence relation
relation on $2^{\aleph_0}$: $\eps\ E_n\ \zeta$ iff
$(\gA^\eps \!\rest\! \alpha^\eps_n,\bar\alpha^\eps \!\rest\! n)$
is isomorphic to $(\gA^\zeta \!\rest\! \alpha^\zeta_n,\bar\alpha^\zeta
\!\rest\! n)$. By this we mean there is an isomorphism $f$ from
$\gA \rest \alpha^\eps_n$ onto $\gA^\zeta 
\rest \alpha^\zeta_n$ which maps $\gA^\eps 
\rest \alpha^\eps_k$ onto $\gA^ \zeta \rest 
\alpha^\zeta_k$ for $k < n$ and is an order preserving function 
(for the ordinals, alternatively we restrict ourselves to the case where $<$ 
is interpreted as a well ordering).

We can find subsets $t^\zeta$ of $\omega$ (for $\zeta < 2^{\aleph_0}$) such that:
\mn
\begin{enumerate}
    \item[$(*)$]
    \begin{enumerate}
        \item For $\zeta,\eps < 2^{\aleph_0}$ {and $n < \omega$} 
        we have $t^\zeta\cap s_n = t^\eps \cap s_n$ \underline{iff} 
        $\gA^\zeta \rest \alpha^\zeta_n = \gA^\eps \rest\alpha^\eps_n$ and $\alpha^\zeta_k=\alpha^\eps_k$ for $k \le n$. 

        \item If $\zeta < 2^{\aleph_0}$ and $n < \omega$ then 
        $t^\zeta\cap s_n$ is infinite.


        \item $t^\zeta \cap s_n$ depends only on $\zeta/E_n$.
    \end{enumerate}
\end{enumerate}
\mn
For $\zeta<2^{\aleph_0}$ let
\[
\clS_\zeta \defeq \big\{a \in [\kappa]^{\aleph_0} : \otp(a) \text{ is a limit ordinal and } t_\zeta = 
\{|C_{\sup}(a) \cap \beta| : \beta \in a\}\big\}
\]
and 
\[
\clS'_\zeta = \big\{ a \in \clS_t : \otp(a) = \otp(\gA^\zeta) \big\},
\]

\mn
and for $a\in \clS'_\zeta$ let $N_a$ be the model isomorphic to 
$\gA^\zeta$ by the function $f_a$, where $\Dom(f_a) = a$, 
$f_a(\gamma) = \otp(\gamma \cap a)$.

Let $\clS$ be the union of $\clS'_\zeta$ for $\zeta<2^{\aleph_0}$. Clearly
$\zeta \ne \xi\ \Rightarrow\ \clS_\zeta\cap \clS_\xi= \varnothing$, and so 
$\clS'_\zeta \cap \clS'_\xi = \varnothing$.  
Hence $N_a$ is well defined for $a\in \clS$.
 
Let $K_n$ be the set of pairs $(\gA,\bar \alpha)$ such that $\gA$ is 
a $\tau_1$-model with universe a countable subset of $\kappa$
with no last member, and $\bar \alpha$ is an increasing	sequence 
of ordinals $ < \kappa$ of length $n$ such that for all $k < n$ we have 
$\alpha_k < \sup(\gA)$,
$[\alpha_k,\alpha_{k+1}) \cap \gA \ne \varnothing$, and 
$\gA \rest \alpha_k \prec \gA$.
So clearly there is a function $\cd_n : K_n \rightarrow \clP(s_n)$
such that for $\zeta < 2^{\aleph_0}$, $\cd_n (\gA,\bar\alpha) =
t^\zeta \cap s_n$ \underline{iff} the pairs $(\gA,\bar\alpha),
(\gA^\zeta,\bar \alpha^\zeta\rest n)$ are isomorphic.

Let $M$ be a $\tau_1$-model with universe $\kappa$.
Now\footnote{See \cite[1.16=L1.15]{Sh:E62} or history in the introduction of
\S3, and the proof of \ref{6.11A}.}
we can find a full subtree $\clT$ of ${}^{\omega>}(\aleph_2)$ (i.e. it
is non-empty, closed under initial segments, and each member has $\aleph_2$
immediate successors) and elementary submodels $N_\eta$ of $M$ for 
$\eta \in \clT$ such that:
\begin{enumerate}
    \item  $\rang(\eta)\subseteq N_\eta$
\sn
    \item  If $\eta$ is an initial segment of $\rho$ then $N_\eta$ 
    is a submodel $N_\rho$. Moreover, $N_\eta \cap \aleph_2$ is an 
    initial segment of $N_\rho$.
\end{enumerate}
\mn
Now let $E$ be the set of $\delta < \kappa = \aleph_2$ satisfying the 
following condition: if $\rho \in \clT \cap {}^{\omega>}\delta$ then
$N_\rho \cap \kappa$ is a bounded subset of $\delta$, and $\delta$ is 
a limit ordinal. Let $E_1$ be the set of $\delta \in E$ such that if 
$\rho \in \clT\cap {}^{\omega>}\delta$ then for every $\beta < \delta$, 
there is $\gamma$ such that $\beta < \gamma < \delta$ and
$\rho \caret \LL\gamma\RR \in \clT$. So
by the choice of $\LL C_\delta : \delta \in S\RR$,
for some $\delta \in S$ we have $C_\delta\subset E_1$.

Let $\LL \alpha_{\delta,k} : k < \omega\RR$ list
$C_\delta$ in increasing order.

Now we choose, by induction on $n$, a quadruple
$(\eta_n,s^*_n,\alpha_n,k_n)$ such that:
\mn
\begin{enumerate}
    \item[$(*)$] 
    \begin{enumerate}
        \item $\eta_n \in \clT$ has length $n$ (so $\eta_0$ is necessarily $\LL\ \RR$).

        \item If $n = m+1$ then $\eta_n$ is a successor of $\eta_m$.

        \item $s^*_n$ is $\cd_n\!\big((N_{\eta_n},\LL \alpha_\ell : \ell < n\RR) \big)$ if the pair $(N_{\eta_n},\LL \alpha_\ell : \ell < n\RR)$ belongs to $K_n$ and is $s_n$ otherwise (so $s_n^* \subseteq s_n$ is infinite). 

        \item $\alpha_n = \sup(N_{\eta_n})+1$

        \item $k_n = \min\{k : N_{\eta_n} \subseteq \alpha_{\delta,k}\}$ and $k_0 = 0$. 

        \item if $n=m+1$ \underline{then}
        \begin{enumerate}
            \item[$(\alpha)$] $\min(N_{\eta_n} \setminus N_{\eta_m}) > \alpha_{\delta,k_n-1}$

            \item[$(\beta)$] $k_m < k_n$

            \item[$(\gamma)$] $k_n \in \bigcup\{s^*_\ell : \ell < n\}$

            \item[$(\delta)$] {If} $n = (n_1 + n_2)^2 + n_2 < (n_1 + n_2 + 1)^2$ (so $n_1,n_2$ are uniquely determined by $n$ and $n_2 < n$) then $k_n \in s_{n_2}^*$.

            \item[$(\eps)$] $k_n$ is minimal under those restrictions.
        \end{enumerate}

    \end{enumerate}
\end{enumerate}
\mn
There is no problem to carry the induction.
In the end, let $\eta = \bigcup\limits_n \eta_n\in \lim(\clT)$, so 
we get a $\tau_1$-model $N_\eta \defeq \bigcup\{N_{\eta_n} : n < \omega\}$
and an increasing sequence $\LL \alpha_n : n < \omega\RR$ of ordinals
with limit $\sup(\gA)$.
Now by the choice of $\LL(\gA^\zeta,\bar \alpha^\zeta) : \zeta < 2^{\aleph_0}\RR$, clearly for some $\zeta$ we have
$(N_\eta,\bar{\alpha})$ isomorphic to $(\gA^\zeta,\bar \alpha^\zeta)$,  
so necessarily $(N_\eta \rest \alpha_n,\bar \alpha\rest n)$
belongs to $K_n$ and $\cd_n (N_{\eta\rest n},\LL \alpha_\ell:
\ell<n\RR)=s^*_n$.

Also, clearly ${\rm sup}(N_\eta) = \delta$ and $\{k_n : n < \omega\} =
\{|C_\delta \cap \beta| : \beta \in N_{\eta}\}$. 

Letting $a$ be the universe of $N_\eta$, it follows that $a \in \clS_\zeta$,
so $N_a$ is well defined and isomorphic to $\gA^\zeta$ (hence to $N_\eta$).
Using $<^M$ we get $N_a = N_\eta$.  But $N_\eta\prec M$, so
$\LL N_a : a \in \clS\RR$ is really a
diamond sequence. (Well, for $\tau_1$-models rather then $\tau$-models,
but this does no harm and will even help us for $\kappa > \aleph_2$.)

Second, we consider the case $\kappa>\aleph_2$. Given a countable vocabulary $\tau$, let $\tau_1 = \tau \cup \{<\}$ (pedantically, assuming $< \notin \tau$) and let $\LL N_c : c \in [\aleph_2]^{\aleph_0}\RR$ be as was proved above with $\kappa = \aleph_2$.
For each $c \in [\kappa]^{\aleph_0}$, if $\otp(c) =
\otp(c\cap\omega_2,<^{N_{c\cap\omega_2}})$, let $g_c$
be the unique isomorphism from $(c\cap\omega_2, <^{N_{c\cap \omega_2}})$
onto $(c,<),<$ the usual order,
and let $M_c$ be the $\tau$-model with universe $c$ such that $g$ is an
isomorphism from $N_{c\cap\omega_2}\rest  \tau $
onto $M_c$.
Clearly it is an isomorphism and the $M_c$-s form a diamond sequence. 

[Why? For notational simplicity $\tau$ has predicates only 
(and, of course, $< \notin \tau$).
Let $M_0=M$ be a $\tau$-model with universe $\kappa$, let $M_1$ be 
an elementary submodel of $M$ of cardinality $\aleph_2$ such that 
$\omega_2\subseteq M_1$, let $h$ be a one-to-one function from $M_1$ onto
$\omega_2$, $M_2$ be a $\tau$-model with universe $\omega_2$ such that
$h$ is an isomorphism from $M_1$ onto $M_2$, and let
$M_3$ be the $\tau_1$-model expanding $M_2$ such that
$$<^{M_3} = \big\{\big(h(\alpha),h(\beta)\big) : \alpha < \beta \text{ are from } M_1\big\}.$$

So for some $a \in \clS \subseteq [\kappa ]^{\aleph_0}$ we have
$N_a\prec M_3$ and 
$$h(\alpha) = \beta \in N_a \wedge \alpha < \omega_2 \Rightarrow \alpha \in a.$$ (Note that the set of $a$-s satisfying 
this contains a club of $[\aleph_2]^{\aleph_0}$.)

Let $c = \{\alpha : h(\alpha) \in a\}$, so clearly $c \cap \omega_2=a$ 
and $M_c\prec M_1$ hence $M_c\prec M$, so we are done.]
\end{PROOF}

\begin{discussion}\label{4.3new}
Some concluding remarks:

\sn
1)  We can use other cardinals, but it is natural if we deal with
$D_{\kappa,<\theta ,\aleph_0}$ (see below).

\sn
2)   The context is very near to \S3, but the stress is different.
\end{discussion}

\begin{definition}
\label{4.4new}
Let $\kappa\ge \theta \ge \sigma$ and $\theta$ be uncountable regular. If
$\theta=\mu^+$ we may write $\mu$ instead of $<\theta$.

\sn
1) Let $D = D_1 = D^1_{\kappa,<\theta,\aleph_0}$ be the filter on
$[\kappa]^{<\theta}$ generated by $\{A^1_x : x \in \clH(\chi)\}$, where
\begin{align*}
    A^1_x = \big\{N\cap\kappa : &\ N = \textstyle\bigcup\limits_{n<\omega} N_n
    \text{ is an elementary submodel of }(\clH(\chi),\in), \\
  &\ N_n \text{ is increasing, } N_n \in N_{n+1},\ \|N_n\| < \theta, 
  \text{ and } N_n \cap \theta \in \theta\big\}. 
\end{align*}
2)  Let $D= D_2 = D^2_{\kappa,<\theta,\sigma}$ be the filter on
$[\kappa]^{<\theta}$ generated by $\{A^2_x: x \in \clH(\chi)\}$, where 
\begin{align*}
A^2_x = \big\{N \cap \kappa : &\ N = \textstyle\bigcup\limits_{\zeta<\sigma} N_\zeta \text{ is an elementary submodel of }
(\clH(\chi),\in),\\
  &\ N_\zeta \text{ increasing, } \LL N_\eps : \eps \le \zeta\RR \in N_{\zeta+1},
\text{ and } N_\eps \cap \theta \in \theta \big\}.
\end{align*}

\mn
3) For a filter $D$ on $[\kappa]^{<\theta}$, let $\diamondsuit_D$
mean the following: fixing any countable vocabulary $\tau$ there are $S\in D$ and $N=\LL
N_a: a\in S\RR$, each $N_a$ a $\tau$-model with universe $a$, such that
for every $\tau$-model $M$ with universe $\lambda$ we have
$\{a\in S: N_a\subseteq M\}\neq \varnothing \mod D$.

\sn
4) If $D$ is a filter on $[\kappa]^{<\theta}$ and $S \in D^+$, then 
$$D \rest S \defeq \big\{X \subseteq [\kappa]^{<\theta} : X \cup \big([\kappa]^{<\theta} \setminus S \big) \in D\big\}.$$ 
\end{definition}

\begin{claim}\label{4.5new}
Assume $\theta \le \sigma$ and $\kappa > \sigma^+$, and let 
$D = D_{\kappa,\theta,\aleph_0}$.

\sn
1) $[\kappa]^{\theta}$ can be partitioned to $\sigma^{\aleph_0}$
(pairwise disjoint) $D$-positive sets.

\sn
2) Assume in addition that $\sigma^{\aleph_0}\ge 2^\theta$. Then
\begin{enumerate}
    \item[$(\alpha)$]  We can find $A_\alpha \subseteq [\kappa]^\theta$ 
    for $\alpha < \lambda \defeq 2^{\kappa^\theta}$ such that each is $D$-positive but they are pairwise disjoint $\mod D$.
\sn
    \item[$(\beta)$]  If $\lambda = \kappa^\theta$ and $\tau$ is a countable vocabulary \underline{then} $\diamondsuit_{\lambda,\theta,\aleph_0}$. Moreover, there exist $S^*\subseteq [\lambda]^\theta$ and a function $N^*$ with domain $S^*$ such that
    \begin{enumerate}
        \item[$(a)$]  For distinct $a,b$ from $S^*$ we have 
        $a \cap \kappa \ne b \cap\kappa$.
\sn
        \item[$(b)$] For $a\in S^*$ we have that $N^*(a) = N^*_a$ is a 
        $\tau$-model with universe $a$.
\sn
        \item[$(c)$]  For a $\tau$-model $M$ with universe $\lambda$, 
        the set $\{a : N^*_a = M \rest a\}$ is stationary.
    \end{enumerate}
\end{enumerate}
\end{claim}

\begin{PROOF}{\ref{4.5new}}
Similar to earlier ones: part (1) like Claim \ref{4.1new} case (a), part
(2) like the proof of Claim \ref{4.2new}.
\end{PROOF}

\begin{claim}\label{4.6new}
1) If $\theta \le \kappa_0 \le \kappa_1$ and $\diamondsuit_{S_0}$
(i.e. $\diamondsuit_{D_{\kappa_0,\theta,\sigma} \rest S_0}$), where $S_0$ is 
a subset of $[\kappa_0]^\theta$ which is
$D_{\kappa_0,\theta,\sigma}$-positive and $S_1 \defeq \{a \in 
[\kappa_1]^\theta : a \cap \kappa_0 \in S_0\}$, then $\diamondsuit_{S_1}$
(i.e. $\diamondsuit_{D_{\kappa_1,\theta,\sigma} \rest S_1}$). 

\sn
2) In part (1), if in addition $\kappa_0=(\kappa_0)^\theta$ and
$\kappa_2= (\kappa_1)^\theta$ \underline{then} we can find $S_2\subseteq
[\kappa_2]^\theta$ such that:
\begin{enumerate}
    \item[$(a)$]  $a\in S_2 \Rightarrow a\cap\kappa_0\in S_0$
\sn
    \item[$(b)$] If $b \neq c \in S_2$ then $b\cap\kappa_1 \neq c\cap\kappa_1$.
\sn
    \item[$(c)$] $\diamondsuit_{S_2}$
\end{enumerate}
\mn
3) If $\kappa=\kappa^\theta$ then $\diamondsuit_{D_{\kappa,\theta,\sigma}}$.
\end{claim}

\begin{remark}
This works for other uniform definitions of normal filters.

Above, $\kappa^{\theta^{\sigma}} = \kappa$ can be replaced by 
``every tree with $\le\theta$ nodes has at most $\theta^*$ branches, and $\kappa^{\theta^*} = \kappa$.''
\end{remark}

\begin{PROOF}{\ref{4.6new}}
1) Easy.

\sn
2) Implicit in earlier proof, \ref{4.2new}.

\sn
3) See \cite{Sh:212}, \cite{Sh:247}
\end{PROOF}
\bigskip

\bibliographystyle{amsalpha}
\bibliography{shlhetal}

\end{document}